\documentclass[11pt]{article}
\usepackage[top=1in, bottom=1in, left=1.25in, right=1.25in]{geometry}

\usepackage{amsmath,amssymb,amsthm,xcolor,enumerate,bbm}
\definecolor{mycitecolor}{rgb}{0,0,0}
\usepackage[unicode,breaklinks=true,colorlinks=true,citecolor=mycitecolor]{hyperref}

\usepackage[]{authblk}
\usepackage{titletoc}
\titlecontents{section}[1.5em]{\vspace{-0.5pt}}
    {\thecontentslabel\enspace\enspace}
    {\hspace*{1.5em}}
   {\titlerule*[0.5pc]{.}\contentspage}
\titlecontents{subsection}[3em]{\vspace{-.5pt}}
    {\thecontentslabel\enspace\enspace}
    {}
    {\titlerule*[0.5pc]{.}\contentspage}
\titlecontents{subsubsection}[3em]{\vspace{-2.5pt}}
    {\thecontentslabel\enspace\enspace}
    {}
    {\titlerule*[0.5pc]{.}\contentspage}

\numberwithin{equation}{section}
\newtheorem{thm}{Theorem}[section]

\newtheorem{lem}[thm]{Lemma}
\newtheorem{prop}[thm]{Proposition}

\theoremstyle{remark}
\newtheorem{remark}[thm]{Remark}
\theoremstyle{definition}
\newtheorem{defn}[thm]{Definition}

\newcommand{\bke}[1]{\left ( #1 \right )}
\newcommand{\bkt}[1]{\left [ #1 \right ]}

\newcommand{\norm}[1]{\left \| #1 \right \|}
\newcommand{\aA}{\mathbb{A}}
\newcommand{\R}{\mathbb{R}}
\newcommand{\QQ}{\mathbb{Q}}
\renewcommand{\div}{\mathop{\rm div}\nolimits}

\newcommand{\pd}{\partial}
\newcommand\ga{\gamma}
\newcommand\De{\Delta}
\newcommand\de{\delta}
\newcommand{\nb}{\nabla}
\newcommand{\lec}{{\ \lesssim \ }}
\newcommand{\wkto}{\rightharpoonup}

\newcommand{\abs}[1]{\left | #1 \right |}
\newcommand\al{\alpha}
\newcommand\ep{\epsilon}
\newcommand\e {\varepsilon}
\newcommand\ze{\zeta}
\renewcommand\th{\theta}
\newcommand\wth{\overline \theta}
\newcommand\ka{\kappa}
\newcommand\la{\lambda}

\newcommand\Th{\Theta}
\newcommand{\ZZ}{\mathbb{Z}}
\newcommand{\NN}{\mathbb{N}}
\newcommand{\esssup} {\mathop{\mathrm{ess\,sup}}}
\newcommand{\supp} {\mathop{\mathrm{supp}}}
\newcommand{\pv} {\mathop{\mathrm{p.v.\!}}}
\newcommand{\tr} {\mathop{\mathrm{tr}}}

\newcommand{\na}{\nabla}
\newcommand{\td}{\tilde}

\newcommand{\EQ}[1]{\begin{equation} #1 \end{equation}}
\newcommand{\EQS}[1]{\begin{equation}\begin{split} #1 \end{split}\end{equation}}
\newcommand{\EQN}[1]{\begin{equation*}\begin{split} #1 \end{split}\end{equation*}}
\newcommand{\EN}[1]{\begin{enumerate} #1 \end{enumerate}}
\newcommand{\uloc}{{\mathrm{uloc}}}
\newcommand{\loc}{{\mathrm{loc}}}

\newif\ifxyzon

\begin{document}

\title{Discretely self-similar solutions to coupled fluid systems with locally $L^2$ data}

\author[1]{Chen-Chih Lai}
\author[2]{Tai-Peng Tsai}

\affil[1]{\footnotesize Department of Mathematical Sciences, Middle Tennessee State University,
Murfreesboro, TN 37132, USA \vspace{8pt}}

\affil[2]{\footnotesize Department of Mathematics, University of British Columbia, Vancouver, BC V6T 1Z2, Canada}

\date{}
\maketitle

\vspace{-10mm}
\begin{abstract}
In this paper we present a systematic way to construct 
discretely self-similar solutions of fluid systems, including the viscous resistive magnetohydrodynamics  equations, the viscoelastic Navier-Stokes equations with damping, and the Oberbeck-Boussinesq system, 
for locally $L^2$ initial data which are not necessarily in the uniformly local $L^2$ space.
 The solutions satisfy local energy inequalities, local a priori bounds depending only on the $L^2$ norm of the initial data in a fixed ball (and, for the Oberbeck--Boussinesq system, on the force), and a limiting formula for the pressure.

\end{abstract}

\vspace{-0.55cm}

\setcounter{tocdepth}{2}

\tableofcontents

\section{Introduction}

The goal of this paper is to construct discretely self-similar (DSS) solutions for DSS data 
in $L^2_\loc$, the locally $L^2$ class (without assuming uniformness), for fluid systems, i.e., coupled systems of the Navier-Stokes equations (NSE) for the fluid velocity $v$ and pressure $\pi$, and the equation for another physical quantity, say $w$. 
These fluid systems should satisfy scaling properties similar to that of NSE, so that it is meaningful to consider self-similar and discretely self-similar solutions. In particular, if $(v(x,t),\pi(x,t),w(x,t))$ is a solution to the system, so is the rescaled triplet of functions
\EQ{
v_\la(x,t)=\la v(\la x, \la^2 t),\quad 
\pi_\la(x,t)=\la^2 \pi(\la x, \la^2 t),\quad 
w_\la(x,t)=\la w(\la x, \la^2 t),
}
for any $\la>0$.
These functions are called \emph{self-similar} if $(v_\la,\pi_\la,w_\la)=(v,\pi,w)$ for any $\la>0$, and \emph{discretely self-similar} (DSS) if $(v_\la,\pi_\la,w_\la)=(v,\pi,w)$ for one particular $\la>1$.  We add the  adjective ``forward'' when these are solutions for positive time, and defined in the space-time domain
\EQ{
\QQ=\R^3 \times (0,\infty).
}
We can also describe time-independent functions as self-similar or DSS similarly by omitting the time dependence.
The initial values $v_0$ and $w_0$ of $v$ and $w$ at time $0$ are
assumed to be locally in $L^2$, $v_0,w_0\in L^2_\loc(\R^3)$, i.e., 
\[
M(x)=\int_{B_1(x)} |v_0(y)|^2 + |w_0(y)|^2 \, dy < \infty,\quad (\forall x\in\R^3),
\]
 where $B_r(x_0)=\{x\in\R^3:|x-x_0|<r\}$ and $B_r=B_r(0)$,  but the above integral need not be bounded \textit{uniformly} in $x$, i.e., we do not assume $v_0$ and $w_0$ are in $L^2_\uloc$, the uniformly local $L^2$ class with the condition $\sup_{x} M(x)<\infty$, and $M(x)$ may grow unbounded as $|x|\to \infty$.
An example is given in \cite[(1-15)]{BT-AnalPDE2019}.

For Navier-Stokes equations,
the construction of forward self-similar solutions for initial data of arbitrary size was initiated by  Jia--\v{S}ver\'ak \cite{JiaSve} for locally H\"older continuous 
$v_0\in C^\al_\loc(\R^3 \setminus \{0\})$. 
It has been generalized to self-similar and DSS solutions in a number of directions by a variety of authors 
 \cite{Tsai2013,Korobkov-Tsai,BT1-AHP,BT2,BT3,MR3916974}
including relaxing the initial data to $L^3_w$ (weak $L^3$) and Besov spaces, 
solutions which satisfy a rotationally corrected scaling invariance, and in the half space setting.

The construction of DSS solutions for $L^2_\loc$ data 
was known only for the Navier-Stokes equations \cite{Chae-Wolf,BT-AnalPDE2019,MR4090469}  and the MHD equations \cite{MR4205085} prior to this work.
It was first done for self-similar solutions with self-similar data by Lemari{\'e}-Rieusset  \cite[Theorem 16.3]{LR16}. As explained in the paragraph after  \cite[(1-13)]{BT-AnalPDE2019}, these data are automatically in $L^2_{\uloc}$ and one can construct  approximation solutions
using the uniform in space local energy bound 
in the framework of local Leray solutions for $L^2_{\uloc}$ data (see \cite{LR02,KS07}).

For DSS solutions for strictly DSS $L^2_\loc$ data, it was first done by
Chae--Wolf \cite{Chae-Wolf} in the class of ``local Leray
solutions with projected pressure", and satisfy the ``local energy inequality with projected pressure".
A second construction was given by Bradshaw--Tsai \cite{BT-AnalPDE2019}, which we will follow.  A third construction, by energy estimates in the weighted spaces $L^2_{w_\gamma}=L^2((1+|x|)^{-\gamma}dx)$ and a fixed-point argument in the class of DSS functions, was given by Fern\'andez-Dalgo and Lemari\'e-Rieusset \cite[Section 7]{MR4090469}; it also gives self-similar solutions, and it was adapted to MHD in \cite{MR4205085}.  
The construction of DSS solutions is much harder than self-similar solutions because, when the data are in $L^2_\loc$ and DSS, they  need not be in $L^2_{\uloc}$. 
There is no existence theory for general data in $L^2_\loc$.  Global weak solutions exist for data in weighted spaces that contain all DSS data in $L^2_\loc$ \cite[Theorem 1]{MR4090469}, \cite{MR4395595,BCT}, but these existence theorems do not by themselves give DSS solutions, since their solutions are obtained by approximations that do not preserve the DSS property. The DSS solutions in \cite[Section 7]{MR4090469} are obtained by a separate fixed-point argument in the class of DSS functions. 
 For general data in $L^2_\loc$, local energy estimates do not close: the nonlocal terms, in particular the pressure, depend on the solution far away, which is not controlled by the local norm of the data. For DSS solutions, however, the solution far away is determined by the solution near the origin at earlier times. For example, for a $\la$-DSS vector field $v$ and $k\in\ZZ$,
\[
\int_{B_{\la^k}}|v(x,t)|^2\,dx=\la^k\int_{B_1}|v(x,\la^{-2k}t)|^2\,dx.
\]
Hence the far-field contributions can be bounded by the local energy $\sup_{0<s\le t}\int_{B_1}|v(x,s)|^2dx$, and the local energy inequality closes as an integral inequality for this quantity. In this way 
we can still recover some a priori bound of the solutions using the DSS property. This is the key insight of \cite{Chae-Wolf}. 

Forward self-similar and DSS solutions have also been constructed in other space dimensions and for the fractional Navier--Stokes equations, with $-\De$ replaced by $(-\De)^\al$; see, e.g., \cite{LMZ-fractional}. We refer to the recent preprint \cite{GLX-2D-DSS} for large $\la$-DSS solutions in two dimensions, for any $\la>1$ and data in $L^2_\loc(\R^2\setminus\{0\})$.

Self-similar and DSS solutions have also been constructed for other fluid systems.  For the Magnetohydrodynamics (MHD) system, self-similar and DSS solutions were constructed by
Lai \cite{Lai-JMFM2019} for $L^3_w$ data (also for the viscoelastic Navier-Stokes equations with damping, veNSD) and by
Zhang and Zhang \cite{MR3997556} for Besov data
(for the generalized MHD system with fractional dissipation $(-\Delta)^\alpha$);
DSS solutions were constructed by Fern\'andez-Dalgo and Jarr\'in \cite{MR4205085}
for weighted $L^2$ data (following the approach of \cite{MR4090469}, which also gives self-similar and DSS solutions);
and Yang \cite{MR4973365} proved pointwise decay estimates for self-similar solutions.
It is also extended to
{O}berbeck-{B}oussinesq (OB) system by Brandolese and Karch \cite{MR4848568} for self-similar solutions with locally bounded data, and
Tsai \cite{Tsai-OB} for DSS solutions with $L^3_w$ data.

In this paper, we extend the approach of \cite{Chae-Wolf,BT-AnalPDE2019} for DSS data in $L^2_\loc$ to these coupled fluid systems.  For the Oberbeck-Boussinesq system and the viscoelastic Navier-Stokes equations with damping, the existence of DSS solutions for DSS data in $L^2_\loc$ appears to be new. For MHD, our result gives a construction different from that of \cite{MR4205085}; our solutions moreover satisfy the local energy inequality with initial data and a local $L^{3/2}$ bound for the pressure.

\medskip

In the following, we describe in detail the fluid systems and the corresponding existence theorems of DSS solutions with local $L^2$ DSS data. 

\subsection{The viscous resistive MHD equations}
The viscous resistive magnetohydrodynamics (MHD) equations comprise a coupled system of PDEs that govern the macroscopic dynamics of an incompressible, electrically conducting fluid. It involves a fluid velocity field $ v:\QQ \to \R^3$, a pressure function $\pi:\QQ\to \R$, and a magnetic field $ b:\QQ \to \R^3$.
The system reads:
\EQS{\label{eq-MHD}
\partial_tv - \Delta v + v\cdot\nabla v - b\cdot\nabla b + \nabla \pi&= 0, \\
\partial_tb - \Delta b + v\cdot\nabla b - b\cdot\nabla v &= 0,\\
\div v = \div b &=0.
}

Scaling invariance: If $(v,b,\pi)$ satisfies \eqref{eq-MHD}, then for any $\la>0$
\[
v^\la(x,t) = \la v(\la x,\la^2 t),\quad
b^\la(x,t) = \la b(\la x,\la^2 t),\quad
\pi^\la(x,t) = \la^2\pi(\la x,\la^2 t),
\]
is also a solution.

A solution of \eqref{eq-MHD} is called \emph{self-similar} (SS) if 
\EQ{\label{eq-MHD-SS-def}
v^\la = v,\quad
b^\la = b,\quad
\pi^\la = \pi
}
for all $\la>0$, and is called \emph{discretely self-similar} with factor $\la$ ($\la$-DSS) if \eqref{eq-MHD-SS-def} is valid for a given $\la>1$.
They correspond to time-independent and time-periodic profiles under \emph{similarity transform}: Let
\EQ{\label{eq-MHD-SS-transform}
v(x,t) = \frac1{\sqrt{t}} V(y,s),\quad
b(x,t) = \frac1{\sqrt{t}} B(y,s),\quad
\pi(x,t) = \frac1t P(y,s),
}
with the similarity variables
\EQ{\label{eq-self-similar-variables}
y = \frac{x}{\sqrt{t}},\quad
s = \log t.
}
The profile $(V,B)$ solves the time-dependent Leray system for the MHD equations
\begin{equation}\label{eq-MHD-Leray}
\begin{aligned}
LV + V\cdot\nb V - B\cdot\nb B+ \nb P&= 0, \\
LB + V\cdot\nb B - B\cdot\nb V &= 0,\\
\div V = \div B &=0,
\end{aligned}
\quad\text{ in }\R^3\times\R,
\end{equation}
in the similarity variables,
where $L$ is the linear operator given by
\EQ{\label{eq-operator-L}
L = \pd_s - \De - \frac12 - \frac12 y\cdot\nb.
}
If $(v,b)$ is self-similar, then $(V,B)$ is independent of $s$. 
If $(v,b)$ is discretely self-similar with factor $\la>1$, then $(V,B)$ is periodic in $s$ with period $2\log \la$.

\begin{thm}\label{thm-MHD}
Assume $v_0$, $b_0\in L^2_\loc(\R^3)$ are divergence-free $\la$-DSS vector fields for some $\la>1$.
Then there exists a $\la$-DSS distributional solution $(v,b)$ to \eqref{eq-MHD} and associated pressure $\pi$ satisfying the local energy inequality: 
 for all nonnegative $\phi\in C^\infty_c(\R^3\times[0,\infty))$ and all $t>0$, we have
\EQS{\label{eq-MHD-LEI}
\int &\left(|v|^2 + |b|^2\right)\phi\,(x,t)\, dx + 
2\int_0^t\!\!\int \left(|\nb v|^2+|\nb b|^2\right)\phi\, dxds \\
&\le \int \left(|v_0|^2 + |b_0|^2\right)\phi(x,0)\, dx
+ \int_0^t\!\!\int \left(|v|^2+|b|^2\right)\left(\pd_s\phi+\De\phi\right) dxds\\
&\quad + \int_0^t\!\!\int \left(|v|^2+|b|^2+2\pi\right) (v\cdot\nb\phi) dxds
 - 2\int_0^t\!\!\int (v\cdot b)(b\cdot\nb\phi)\, dxds,
}
and $(v,b)$ satisfies 
\[
\lim_{t\to0^+}\norm{v(t)-v_0}_{L^2(E)} + \norm{b(t)-b_0}_{L^2(E)}= 0
\]
for all compact subsets $E$ of $\R^3$. Moreover, for any $T>0$ and compact subset $E$ of $\R^3$, we have $v, b\in L^\infty(0,T;L^2(E))\cap L^2(0,T;H^1(E))$ and $\pi\in L^{3/2}(0,T;L^{3/2}(E))$, with
\EQ{\label{MHD-apriori}
\norm{v, b}_{L^\infty(0,T;L^2(E))\cap L^2(0,T;H^1(E))}+ \norm{\pi}_{L^{3/2}(0,T;L^{3/2}(E))} \le C(T,E,\la,\al_0)
}
where $\al_0=\norm{v_0}_{L^2(B_\la)}^2 + \norm{b_0}_{L^2(B_\la)}^2$.
Furthermore,  the pressure satisfies the formula
\EQ{\label{eq-pressure-MHD}
\pi(x,t) = -\frac13(|v|^2-|b|^2)(x,t) + \lim_{\de\to0} \int_{| x-y|>\de} K_{ij}(x-y)\bke{v_iv_j - b_ib_j}(y,t)\, dy
}
in $L^{3/2}_\loc(\R^3\times[0,\infty))$, where $K(x) = \nb^2(1/(4\pi|x|))$.

If furthermore $(v_0,b_0)$ is self-similar, then $(v,b)$ can be chosen self-similar.
\end{thm}

\textit{Comments on Theorem \ref{thm-MHD}:}
\EN{

\item For small self-similar data $(v_0,b_0)$ in suitable Besov, Lorentz, or pseudomeasure spaces, He and Xin \cite{MR2514362} constructed self-similar solutions, with uniqueness in the corresponding small-data solution classes.

\item When the data $(v_0,b_0)$ is DSS and in $L^3_w(\mathbb R^3)$, a DSS local Leray solution is constructed in \cite{Lai-JMFM2019}, where a self-similar local Leray solution is also constructed for self-similar data in the same $L^3_w$ class. This includes  self-similar/DSS data in $C^\al_\loc(\mathbb R^3\setminus\{0\})$, $0<\al<1$.

\item When the data $(v_0,b_0)$ is DSS and in the critical Besov space $\dot B^{-1+3/p}_{p,\infty}(\mathbb R^3)$, $3<p<6$, a DSS  distributional solution is constructed in \cite{MR3997556}, where a self-similar  distributional solution is also constructed for self-similar data in the same critical Besov space.
This is the case $\alpha=1$ of their results for the generalized MHD system with fractional dissipation $(-\Delta)^\alpha$.

\item When the data $(v_0,b_0)$ is DSS and in weighted spaces $L^2_{w_\gamma}(\R^3)= L^2((1+|x|)^{-\gamma}\,dx)$, $\gamma\in(1,2)$, a DSS suitable weak solution is constructed in \cite[Theorem 2]{MR4205085}, following the approach of \cite{MR4090469} for the Navier-Stokes equations.
Since DSS data in $L^2_\loc(\R^3)$ belong to $L^2_{w_\gamma}(\R^3)$, this includes all DSS data in $L^2_\loc(\R^3)$.
The contribution of Theorem \ref{thm-MHD} is to construct the solutions by a different method. At first sight, the a priori bound \eqref{MHD-apriori} may seem stronger than the bounds in \cite[Theorem 2]{MR4205085}: it holds on any compact set $E$ and its constant depends only on the local norm $\al_0$ of the data, whereas the bounds there are global in the weighted space, $v,b\in L^\infty(0,T;L^2_{w_\gamma})$ and $\nb v,\nb b\in L^2(0,T;L^2_{w_\gamma})$, with constants depending on $\norm{(v_0,b_0)}_{L^2_{w_\gamma}}$.
However, for DSS functions and $\gamma>1$ the two are equivalent. Indeed, with $A_k=\{x:\la^{k-1}\le|x|<\la^k\}$, the DSS scaling gives
$\int_{A_k}|v_0|^2\,dx=\la^k\int_{A_0}|v_0|^2\,dx$ and
$\int_0^T\!\int_{A_k}|\nb v|^2\,dx\,dt=\la^k\int_0^{\la^{-2k}T}\!\int_{A_0}|\nb v|^2\,dx\,dt$, and similarly for $b$ and for $\sup_t\int_{A_k}|v(t)|^2dx$, while $w_\gamma\sim\la^{-k\gamma}$ on $A_k$. Since $\sum_{k\ge1}\la^{k(1-\gamma)}<\infty$, the weighted norms are bounded by local ones, and conversely.
For the pressure, the bound in \cite{MR4205085} is in a weighted $L^\beta$ space with $\beta<3/2$, and a local $L^{3/2}$ bound needs an additional near/far splitting argument as in Lemma \ref{BT2019-lem-2.1}.

\item When the data $(v_0, b_0)$ is self-similar and in $L^\infty_\loc(\R^3\setminus\{0\})$, 
Yang \cite{MR4973365} constructed 
a self-similar solution $(v,b)\in BC_w([0,\infty); L^3_w(\R^3))$ and  proved its pointwise decay estimates.

}

\subsection{The viscoelastic Navier-Stokes equations with damping}
The viscoelastic Navier-Stokes equations with damping (veNSD) are a modification of the classical fluid dynamics equations that act as a simplified or regularized model for complex, non-Newtonian fluids (like polymer melts or Oldroyd fluids). It involves a fluid velocity field $ v:\QQ \to \R^3$, a pressure function $\pi:\QQ\to \R$, and an elastic deformation tensor $\mathbf{f}:\QQ\to \R^{3\times 3}$.
The system reads:
\EQS{\label{eq-veNSD}
\partial_tv - \Delta v + v\cdot \nabla v - \nabla\cdot({\bf f}\,{\bf f}^\top) + \nabla \pi&= 0,\\
\partial_t{\bf f} - \Delta{\bf f} + v\cdot\nabla{\bf f} - (\nabla v){\bf f} &= 0,\\
\div v = \div f_m&=0,\quad m=1,2,3,
}
where 
$f_m$ is the $m$-th column vector of ${\bf f}$. 
If ${\bf f} =(f_{ij})_{i,j=1}^3$, then $f_m = (f_{1m},f_{2m},f_{3m})^\top$.
We use the convention that
$(\div {\bf f} )_i = \pd_j { f}_{ji}$ and $(\nb v)_{ij} = \pd_j v_i$, so  $\div f_m=0$ for all $m$ is the same as $\div {\bf f}=0$.
Note ${\bf f}\,{\bf f}^\top$ is symmetric and 
\[
(\nb\cdot({\bf f}\,{\bf f}^\top))_i =  \pd_j ({f}_{jk}{f}_{ik}) = {f}_{jk}\pd_j {f}_{ik} = (f_k \cdot \nb f_k)_i.
\]

Scaling invariance: If $(v,{\bf f},\pi)$ satisfies \eqref{eq-veNSD}, then for any $\la>0$
\[
v^\la(x,t) = \la v(\la x,\la^2 t),\quad
{\bf f}^\la(x,t) = \la {\bf f}(\la x,\la^2 t),\quad
\pi^\la(x,t) = \la^2\pi(\la x,\la^2 t)
\]
is also a solution.

A solution of \eqref{eq-veNSD} is called \emph{self-similar} (SS) if 
\EQ{\label{eq-veNSD-SS-def}
v^\la = v,\quad
{\bf f}^\la = {\bf f},\quad
\pi^\la = \pi
}
for all $\la>0$, and is called \emph{discretely self-similar} with factor $\la$ ($\la$-DSS) if \eqref{eq-veNSD-SS-def} is valid for a given $\la>1$.
They correspond to time-independent and time-periodic profiles under similarity transform: Let
\EQ{\label{eq-veNSD-SS-transform}
v(x,t) = \frac1{\sqrt{t}} V(y,s),\quad
{\bf f}(x,t) = \frac1{\sqrt{t}} {\bf F}(y,s),\quad
\pi(x,t) = \frac1t P(y,s),
}
with the similarity variables \eqref{eq-self-similar-variables} and ${\bf F} = (F_1, F_2, F_3)$.

The profile $(V,{\bf F})$ solves the time-dependent Leray system for \eqref{eq-veNSD}
\begin{equation}\label{eq-veNSD-Leray}
\begin{aligned}
LV + V\cdot\nb V - \sum_{n=1}^3 F_n\cdot\nb F_n+ \nb P&= 0, \\
LF_m + V\cdot\nb F_m - F_m\cdot\nb V &= 0,\\
\div V = \div F_m &=0,
\end{aligned}
\quad\text{ in }\R^3\times\R,\quad m=1,2,3,
\end{equation}
in the similarity variables, where $L$ is the linear operator given by \eqref{eq-operator-L}.

If $(v,{\bf f})$ is self-similar, then $(V,{\bf F})$ is independent of $s$. 
If $(v,{\bf f})$ is discretely self-similar with factor $\la>1$, then $(V,{\bf F})$ is periodic in $s$ with period $2\log \la$.

\begin{thm}\label{thm-veNSD}
Assume $v_0$, $(f_m)_0\in L^2_\loc(\R^3)$, $m=1,2,3$, are divergence-free $\la$-DSS vector fields for some $\la>1$.
Then there exists a $\la$-DSS distributional solution $(v,{\bf f})$ to \eqref{eq-veNSD} and associated pressure $\pi$ satisfying the local energy inequality: 
 for all nonnegative $\phi\in C^\infty_c(\R^3\times[0,\infty))$ and all $t>0$, we have
\EQS{\label{eq-veNSD-LEI}
\int &\left(|v|^2 + |{\bf f}|^2\right)\phi\,(x,t)\, dx + 
2\int_0^t\!\!\int \left(|\nb v|^2+|\nb {\bf f}|^2\right)\phi\, dxds \\
&\le \int \left(|v_0|^2 + |{\bf f}_0|^2\right)\phi(x,0)\, dx
+ \int_0^t\!\!\int \left(|v|^2+|{\bf f}|^2\right)\left(\pd_s\phi+\De\phi\right) dxds\\
&\quad + \int_0^t\!\!\int \left(|v|^2+|{\bf f}|^2+2\pi\right) (v\cdot\nb\phi) dxds
 - 2\sum_{n=1}^3\int_0^t\!\!\int (v\cdot f_n)(f_n\cdot\nb\phi)\, dxds,
}
and $(v,{\bf f})$ satisfies 
\[
\lim_{t\to0^+}\norm{v(t)-v_0}_{L^2(E)} + \sum_{n=1}^3 \norm{f_n(t)-(f_n)_0}_{L^2(E)}= 0
\]
for all compact subsets $E$ of $\R^3$.
Moreover, for any $T>0$ and compact subset $E$ of $\R^3$, we have $v, f_m\in L^\infty(0,T;L^2(E))\cap L^2(0,T;H^1(E))$, $m=1,2,3$, and $\pi\in L^{3/2}(0,T;L^{3/2}(E))$, with
\EQ{\label{veNSD-apriori}
\norm{v, \bf f}_{L^\infty(0,T;L^2(E))\cap L^2(0,T;H^1(E))}+ \norm{\pi}_{L^{3/2}(0,T;L^{3/2}(E))} \le C(T,E,\la,\al_0)
}
where $\al_0=\norm{v_0}_{L^2(B_\la)}^2 + \norm{{\bf f}_0}_{L^2(B_\la)}^2$.
Furthermore,  the pressure satisfies the formula
\EQS{\label{eq-pressure-veNSD}
\pi(x,t) &= -\frac13\bke{|v|^2-\sum_{n=1}^3|f_n|^2}(x,t)\\
&\quad + \lim_{\de\to0} \int_{| x-y|>\de} K_{ij}(x-y)\left(v_i(y,t)v_j(y,t) - \sum_{n=1}^3 (f_n)_i(y,t)(f_n)_j(y,t)\right) dy
}
in $L^{3/2}_\loc(\R^3\times[0,\infty))$, where $K(x) = \nb^2(1/(4\pi|x|))$.

If furthermore $(v_0,{\bf f}_0)$ is self-similar, then $(v,{\bf f})$ can be chosen self-similar.
\end{thm}

\textit{Comments on Theorem \ref{thm-veNSD}:}
\EN{
\item When the data $(v_0,{\bf f}_0)$ is self-similar and in $C^\al_\loc(\mathbb R^3\setminus\{0\})$, $0<\al<1$, a self-similar solution is constructed by Lai, Lin, and Wang in \cite{MR3609234}.

\item When the data $(v_0, {\bf f}_0)$ is DSS and in $L^3_w(\mathbb R^3)$, a DSS local Leray solution is constructed in \cite{Lai-JMFM2019}, where a self-similar local Leray solution is also constructed for self-similar data in the same $L^3_w$ class. 
}

\subsection{Oberbeck-Boussinesq system}
The Oberbeck-Boussinesq (OB) system comprises a coupled system of PDEs that govern the motion of a fluid driven by buoyancy forces due to temperature variations. It acts as a foundational model for natural convection, atmospheric flows, and oceanic circulation by assuming the fluid density is constant everywhere except in the buoyancy term of the momentum equation. It involves a fluid velocity field \(v: \QQ \to \R^3\), a pressure function \(\pi: \QQ \to \R\), and a temperature field \(\theta: \QQ \to \R\). The system reads:
\EQS{\label{eq-OB}
\pd_t v-\De v +v\cdot\na v+\na\pi &= \th \nb G+f,\\
\div v&=0,\\
\pd_t \th-\De \th +v\cdot\na \th& = 0,
}
where $G(x) = |x|^{-1}$ is the gravity potential. It is chosen (instead of the usual $G(x)=gx_3$) to have a suitable scaling property.

Scaling invariance: If $(v,\th,\pi)$ satisfies \eqref{eq-OB} with force $f$, then for any $\la>0$
\[
v^\la(x,t) = \la v(\la x,\la^2 t),\quad
\th^\la(x,t) = \la\th(\la x,\la^2 t),\quad
\pi^\la(x,t) = \la^2\pi(\la x,\la^2 t)
\]
is also a solution of the same system with force $f^\la(x,t) = \la^3 f(\la x, \la^2 t)$.

A solution of \eqref{eq-OB} is called \emph{self-similar} (SS) if 
\EQ{\label{eq-OB-SS-def}
v^\la = v,\quad
\th^\la = \th,\quad
\pi^\la = \pi,\quad 
f^\la = f
}
for all $\la>0$, and is called \emph{discretely self-similar} with factor $\la$ ($\la$-DSS) if \eqref{eq-OB-SS-def} is valid for a given $\la>1$.
They correspond to time-independent and time-periodic profiles under similarity transform: Let
\EQS{\label{eq-OB-SS-transform}
v(x,t) &= \frac1{\sqrt{t}} V(y,s),\quad
\th(x,t) = \frac1{\sqrt{t}} \Th(y,s),\\
\pi(x,t) &= \frac1t P(y,s),\quad\ \
f(x,t) = \frac1{t^{3/2}} F(y,s),
}
with the similarity variables \eqref{eq-self-similar-variables}.
Note that $G(x) = \frac1{\sqrt{t}} G(y)$.
The profile $(V,\Th)$ solves the time-dependent Leray system for the Oberbeck-Boussinesq system
\begin{equation}\label{eq-OB-Leray}
\begin{aligned}
LV + V\cdot\nb V + \nb P&= \Th\nb G + F, \\
\div V &=0,\\
L\Th + V\cdot\nb\Th&= 0,
\end{aligned}
\quad\text{ in }\R^3\times\R,
\end{equation}
in the similarity variables, where $L$ is the linear operator defined in \eqref{eq-operator-L}.

If $(v,\th)$ is self-similar, then $(V,\Th)$ is independent of $s$. 
If $(v,\th)$ is discretely self-similar with factor $\la>1$, then $(V,\Th)$ is periodic in $s$ with period $2\log \la$.

\begin{thm}\label{thm-OB}
Assume $v_0$, $\th_0\in L^2_\loc(\R^3)$  are a divergence-free vector field and a function that are $\la$-DSS for some $\la>1$.
Assume $f=f_1+f_2$, $f_1,f_2$ are 
$(\la,3)$-DSS vector fields, $f_j(x,t) = \la^3f_j(\la x,\la^2t)$, $j=1,2$,
\[
 f_1 \in L^\infty(1,\la^2; L^{6/5,2}(\R^3)),  \quad f_2\in L^\infty(1,\la^2; L^2(\R^3)),\quad
\div f_2=0.
\]
Then there exists a $\la$-DSS distributional solution $(v,\th)$ to \eqref{eq-OB} and associated pressure $\pi$ satisfying the local energy inequalities: 
 for all nonnegative $\phi\in C^\infty_c(\R^3\times[0,\infty))$ and all $t>0$, we have
\EQS{\label{eq-OB-LEI}
\int \th^2\phi\,(x,t)\, dx &+ 2\int_0^t\!\!\int |\nb\th|^2\phi\, dxds\\
&\le \int \th_0^2\phi(x,0)\, dx + \int_0^t\!\!\int \th^2\left(\pd_s\phi + \De\phi + v\cdot\nb\phi\right) dxds,\\
\int \bke{|v|^2+\th^2}\phi\,(x,t)\, dx &+ 2\int_0^t\!\!\int \bke{|\nb v|^2+|\nb\th|^2}\phi\, dxds\\
&\le \int\bke{|v_0|^2+\th_0^2}\phi(x,0)\, dx + \int_0^t\!\!\int \bke{|v|^2+\th^2}\left(\pd_s\phi+\De\phi\right) dxds\\
&\quad + \int_0^t\!\!\int \bke{|v|^2+\th^2+2\pi} (v\cdot\nb\phi)\, dxds\\
&\quad + 2\int_0^t\!\!\int (\th\nb G + f)\cdot v\phi\, dxds,
}
and $(v,\th)$ satisfies 
\[
\lim_{t\to0^+}\norm{v(t)-v_0}_{L^2(E)} + \norm{\th(t)-\th_0}_{L^2(E)}= 0
\]
for all compact subsets $E$ of $\R^3$.
Moreover, for any $T>0$ and compact subset $E$ of $\R^3$, we have $v, \th\in L^\infty(0,T;L^2(E))\cap L^2(0,T;H^1(E))$ and $\pi\in L^{3/2}(0,T;L^{3/2}(E))$, with
\EQ{\label{OB-apriori}
\norm{v, \th}_{L^\infty(0,T;L^2(E))\cap L^2(0,T;H^1(E))}+ \norm{\pi}_{L^{3/2}(0,T;L^{3/2}(E))} \le C(T,E,\la,\al_0, M_1,M_2)
}
where $\al_0=\norm{v_0}_{L^2(B_\la)}^2 + \norm{\th_0}_{L^2(B_\la)}^2$, $M_1=\norm{f_1}_{L^\infty(1,\la^2;L^{6/5,2})}$ and $M_2=\norm{f_2}_{L^\infty(1,\la^2;L^2)}$.
Furthermore,  the pressure satisfies the formula
\EQS{\label{eq-pressure-OB}
\pi(x,t) &= -\frac13|v|^2(x,t) + \lim_{\de\to0} \int_{| x-y|>\de} K_{ij}(x-y)v_i(y,t)v_j(y,t)\, dy \\
&\quad - \int \nb\bke{\frac1{4\pi|x-y|}}\cdot(\th\nb G + f_1)(y,t)\, dy 
}
in $L^{3/2}_\loc(\R^3\times[0,\infty))$, where $K(x) = \nb^2(1/(4\pi|x|))$.

If furthermore $(v_0,\th_0)$ is self-similar, and  $f_1$ and $f_2$ are self-similar with degree 3, then $(v,{\th})$ can be chosen self-similar.
\end{thm}

\textit{Comments on Theorem \ref{thm-OB}:}
\EN{
\item When the data $(v_0,\th_0)$ is self-similar and in 
$L^\infty_\loc(\R^3 \setminus \{0\})$, and $f(x,t)=t^{-3/2} F(\frac x{\sqrt t})$ with $F \in H^{-1}(\R^3)$, a self-similar solution $(v,{\th})$ is constructed in Brandolese-Karch \cite{MR4848568}.

\item When the data $(v_0,\th_0)$ is DSS and in 
 $L^3_w(\R^3)$, and  $f(x,t) = \la^3f(\la x,\la^2t)$, with 
  $f \in L^\infty(1,\la^2; (L^{6/5,2}+L^2)(\R^3))$,
 a DSS solution $(v,{\th})$ satisfying  local energy inequalities (see Comment 4 below)  is constructed in Tsai \cite{Tsai-OB} (see also the correction \cite{Tsai-OB-corr}), where a self-similar solution is also constructed for self-similar data in the same $L^3_w$ class.

\item  Compared to \cite{MR4848568}, which treats self-similar forces $f(x,t)=t^{-3/2}F(x/\sqrt t)$ and only requires $F\in H^{-1}(\R^3)$, \cite{Tsai-OB} requires $f =f_1+f_2$ with
$f_1\in L^\infty(1,\la^2;L^{6/5,2}(\R^3))$ and $f_2\in L^\infty(1,\la^2; L^2(\R^3))$. In particular, $f$ is a function and the local energy inequality \eqref{eq-OB-LEI} is meaningful. 
 Since $L^{6/5,2}(\R^3)\subset\dot H^{-1}(\R^3)$, we have $f_1=\div \aA$ with $\aA_{jk}=-\pd_j(-\De)^{-1}(f_1)_k\in L^\infty(1,\la^2; L^2(\R^3))$, but we do not need this form. Theorem \ref{thm-OB} further assumes
$\div f_2=0$, so that $f_2$ is absent
in the  integral formula \eqref{eq-pressure-OB} for the pressure. 
 If $f_2$ remained in \eqref{eq-pressure-OB}, our bound of its contribution to the pressure near $t=0$ would diverge, since $\norm{f_2(t)}_{L^2}\sim t^{-3/4}$ by Lemma \ref{lem2.4} and $t^{-3/4}\notin L^{3/2}(0,1)$; see  the remark after \eqref{pnear.est}.
 For a general $f_2\in L^\infty(1,\la^2;L^2(\R^3))$, write $f_2=\mathbb Pf_2+\nb q$, where $\mathbb P$ is the  Helmholtz projection. Then $\mathbb Pf_2$ satisfies the assumptions of Theorem \ref{thm-OB}, and a solution for the force $f_1+\mathbb Pf_2$ gives a solution for the force $f$, with the same $(v,\th)$ and the pressure $\pi+q$. However, $q$ need not belong to $L^{3/2}_\loc$ near $t=0$. So the assumption $\div f_2=0$ only concerns the normalization of the pressure.  We apply the same idea to $f_1$ in the proofs of Lemma \ref{BT2019-lem-2.1-OB} and Proposition \ref{BT2019-prop-3.1-OB}: there the pressure part $-(-\De)^{-1}\div f_1$ does belong to $L^{3/2}_\loc$ near $t=0$, so we may assume $\div f_1=0$ without loss of generality.

\item The second inequality in \eqref{eq-OB-LEI}  has the form of the sum of the local energy inequality for $v$ and  the first inequality. Adding further multiples of the first inequality, it also holds with $\th^2$ and $|\nb\th|^2$ replaced by $\mu\th^2$ and $\mu|\nb\th|^2$ for any $\mu\ge1$.
We do not know whether the local energy inequality for $v$ alone holds for our solutions. The obstruction is the buoyancy term $2\iint\th\nb G\cdot v\phi$: since $|\nb G|=|x|^{-2}$, this term is critical with respect to Hardy's inequality, and the energy bounds do not prevent concentration at $x=0$ when passing to the limit. Lemma \ref{lem-buoyancy} shows that a possible defect is absorbed by the dissipation in the combined inequality.
 The same applies to the solutions of \cite{Tsai-OB}. There the two inequalities are stated separately, in the weak form \cite[(1.4)]{Tsai-OB}, and are obtained by passing to the limit in the local energy equalities of the approximations; this does not control the product of the differences between the approximations and their limits in the buoyancy term, which is the term $g_kh_k$ in the proof of Lemma \ref{lem-buoyancy}. So the local energy inequality for $v$ alone in \cite[(1.4)]{Tsai-OB} is not justified by that argument.  This is corrected in \cite[Section 2]{Tsai-OB-corr}, where \cite[(1.4)]{Tsai-OB} is replaced by the weak form of \eqref{eq-OB-LEI}. What can be proved, with the modified scheme of Remark \ref{approx-scheme}, is \eqref{eq-OB-LEI}, by Proposition \ref{BT2019-prop-3.1-OB}.

}

\subsection{Strategy of the proof}\label{S:strategy}
Because the proof of Theorem \ref{thm-veNSD} is very similar to Theorem \ref{thm-MHD}, we will skip it, and only present the proofs of Theorems \ref{thm-MHD} and \ref{thm-OB}.

 The statements for self-similar data in Theorems \ref{thm-MHD}--\ref{thm-OB} follow from the DSS case by the limiting argument of \cite[Section 5.1]{BT1-AHP}, with factors $\la_k\to1$; see Subsection \ref{S:SS}.

\paragraph{Approximation.}  DSS solutions for data in $L^{3}_w(\R^3)$ will be used for the purpose of approximations in this paper. As Bradshaw--Tsai \cite{BT-AnalPDE2019} used $L^{3}_w$-data solutions in \cite{BT1-AHP} for NSE, 
we will use $L^{3}_w$-data solutions in Lai  \cite{Lai-JMFM2019} for MHD and veNSD, and Tsai \cite{Tsai-OB} for OB. 
Since $L^{3}_w$-data belong to $L^2_\uloc$, these
solutions are in the framework of local Leray solutions for $L^2_{\uloc}$ data (see \cite{LR02,KS07}).

Our proof adapts the approximation argument of \cite{BT-AnalPDE2019} to the MHD equations and the OB system.
We first approximate the initial data in $L^2(B_1)$ by DSS data in $L^3_w$, preserving the divergence-free constraints on the vector fields.
For the vector fields, this approximation is provided by \cite[Lemma~4.1]{BT-AnalPDE2019}; the scalar temperature is treated by an analogous annular approximation without a divergence correction.
The existence results of \cite{Lai-JMFM2019,Tsai-OB} then provide a sequence of DSS solutions satisfying the corresponding local energy inequalities.  For MHD and veNSD, the mollified approximation scheme of \cite{Lai-JMFM2019} needs an additional pressure term in the equation for $b$ (resp.\ ${\bf f}$), which vanishes in the limit; see Remark \ref{rem-pib}.

\paragraph{Limiting pressure formula.} 
As in \cite{BT-AnalPDE2019}, the key is to derive a \textbf{limiting pressure formula} of the form
\EQ{\label{pi.form}
\pi (x,t)= \pi^1(x,t) + \lim_{\e \to 0} \int_{|y-x|>\e} K_{ij} (x-y)\mathbb F_{ij}(y,t)\,dy
}
in $L^{3/2}_\loc(\R^3 \times [0,\infty))$,
where $K(x) = \nb^2(1/(4\pi|x|))$ denotes the Hessian matrix of the Newtonian kernel, $\mathbb F=v\otimes v + \cdots $ is quadratic in $v$ and $w$, and $\pi^1=-\frac13\sum_i\mathbb F_{ii}$, plus terms from the forces for the Oberbeck-Boussinesq system; see \eqref{eq-pressure-OB}.
The integral in \eqref{pi.form} is not a  Calder\'on--Zygmund  singular integral because we do not  assume a global bound for $\mathbb F$.
 Rather, it is defined in $L^{3/2}_\loc$ for each fixed $\e >0$
using the DSS property of $\mathbb F$.

The limiting pressure formula should not be confused with the \textbf{local pressure decomposition formula} which plays an essential role in the theory of local Leray or local energy solutions for general initial data in $L^2_\uloc$ \cite{LR02,KS07,JiaSve, BT1-AHP,Lai-JMFM2019,KwTs,MR4395595,BCT,BT-SIMA21,BLT}, and is of the form
\begin{multline}\label{pressure-decomposition}
\pi (x,t) -c_{x_0}(t)= - \frac13 \sum_i \mathbb F_{ii}(x,t)+ \pv\int_{B_2(x_0)} K_{ij} (x-y)\mathbb F_{ij}(y,t)\,dy
\\
+ \int_{\R^3\setminus B_2(x_0)} \bkt{K_{ij} (x-y)-K_{ij} (x_0-y)}\mathbb F_{ij}(y,t)\,dy
\end{multline}
which is valid in $L^{3/2}(0,T;L^{3/2}(B_1(x_0)))$ for every $x_0$, for some $c_{x_0}\in L^{3/2}(0,T)$.
It allows us to bound the pressure locally using very weak global control of the velocity $v$ and the other quantity $w$. Due to the extra decay of $K(x-y)-K(x_0-y)$ in $y$, the last integral converges absolutely when $\mathbb F_{ij}$ are integrable with some bound. In contrast, for  general initial data  in $L^2_\loc$, which belong neither to $L^2_\uloc$ nor to the weighted spaces of \cite{MR4090469,MR4395595,BCT},  $\mathbb F_{ij}$ are generally unbounded and we cannot hope for the validity of \eqref{pressure-decomposition}. 
 For DSS solutions, however, \eqref{pressure-decomposition} follows from \eqref{pi.form}. Indeed, if $\mathbb F$ is $(\la,2)$-DSS (see Definition \ref{def-DSS}) and $\sup_{0<\tau<T}\norm{\mathbb F(\tau)}_{L^1(B_1)}<\infty$, then the DSS scaling gives $\int_{A_k}|\mathbb F(y,t)|\,dy\le \la^k\sup_{0<\tau\le t}\norm{\mathbb F(\tau)}_{L^1(B_1)}$ on the annuli $A_k=\{y:\la^{k-1}\le|y|<\la^k\}$, $k\ge1$. Since $|K(x_0-y)|\lesssim|y|^{-3}$ for $|y|$ large, the integral
\[
c_{x_0}(t)=\int_{|y-x_0|>2}K_{ij}(x_0-y)\,\mathbb F_{ij}(y,t)\,dy
\]
converges absolutely and $c_{x_0}\in L^\infty(0,T)$. Splitting the integral in \eqref{pi.form} at $|y-x_0|=2$ and subtracting $c_{x_0}(t)$ gives \eqref{pressure-decomposition}.
In particular, our solutions and those in \cite{BT-AnalPDE2019} for DSS initial data in $L^2_\loc$ satisfy \eqref{pressure-decomposition}; for the Oberbeck--Boussinesq system, the force term of \eqref{eq-pressure-OB} has to be added to the right side of \eqref{pressure-decomposition}. However, we do not use \eqref{pressure-decomposition} for the solutions from  \cite{BT1-AHP,Lai-JMFM2019,Tsai-OB} when we use them as approximations.

\paragraph{Comparison of methods.} 
The constructions in \cite{Lai-JMFM2019,Tsai-OB} work in similarity variables and obtain DSS solutions by constructing time-periodic solutions of the associated Leray systems, using energy estimates for perturbations of suitably revised background profiles and Galerkin approximations.
We use these solutions as approximants and obtain solutions for rougher initial data through local energy estimates and compactness.
The main difficulty is that DSS data in $L^2_\loc$ need not belong to $L^2_\uloc$, and the $L^3_w$ norms of the approximating data need not remain bounded.
Consequently, estimates depending on these stronger norms do not suffice for the limiting argument.
Following \cite{BT-AnalPDE2019}, we exploit discrete self-similarity to relate larger spatial scales to earlier times, thereby closing the local energy estimates and controlling the nonlocal pressure contributions.

The coupled systems require additional estimates. For MHD, the combined local energy estimate uses the cancellation between the velocity and magnetic coupling terms, while the pressure representation involves the quadratic tensor $v\otimes v-b\otimes b$.
For the OB system, the buoyancy term and its contribution to the pressure require separate control: we use a localized Hardy inequality, with the temperature energy  suitably weighted, and a separate far-field bound for the buoyancy part of the pressure.  Because the buoyancy term is critical, we obtain the local energy inequality for $v$ only in combination with that for $\th$; see Comment 4 after Theorem \ref{thm-OB}. 

\subsection{Organization and notation}\label{S:notation}
 In Section \ref{sec:DSS}, we collect elementary properties of DSS functions.
In Section \ref{sec:limiting-pressure}, we establish the stability of the limiting pressure formula under the local convergence arising in our approximation scheme.
In Section \ref{S:L3w-data}, we justify the limiting pressure formula for the approximating solutions and derive local energy and pressure estimates.  We also prove there a lemma for passing to the limit in the buoyancy term of the Oberbeck-Boussinesq system. 
The bounds depend only on the $L^2(B_\la)$ norms of the initial data (and, for OB, on the force), independently of their $L^3_w$ norms.
These uniform estimates provide the compactness needed to pass to a subsequential limit.
In Section \ref{sec:construction}, we show that a DSS solution in a neighborhood of the origin extends to a DSS solution on $\QQ$, and complete the construction; in Subsection \ref{S:SS}, we treat self-similar data.

We use the following notation.
\begin{itemize}
\item $\QQ=\R^3\times(0,\infty)$; $B_r(x_0)=\{x\in\R^3:|x-x_0|<r\}$ and $B_r=B_r(0)$; for $\la>1$ and $k\in\ZZ$, $A_k=\{x\in\R^3:\la^{k-1}\le|x|<\la^k\}$.
\item $L^{q,r}(\R^3)$ denotes the Lorentz spaces, and $L^3_w=L^{3,\infty}$ is the weak $L^3$ space. $L^2_\uloc$ is the space of uniformly locally square integrable functions, with norm $\sup_{x\in\R^3}\norm{f}_{L^2(B_1(x))}$.
\item $K(x)=\nb^2\big(1/(4\pi|x|)\big)$, i.e., $K_{ij}(x)=\pd_i\pd_j\frac1{4\pi|x|}$, and repeated indices are summed.
\item A function is $(\la,a)$-DSS if it is discretely self-similar with factor $\la$ and degree $a$, see Definition \ref{def-DSS}; we simply say $\la$-DSS when $a=1$.
\item $(f,g)=\int_{\R^3}f\cdot g\,dx$, and $C$ denotes a constant that may change from line to line.
\end{itemize}

\section{Discretely self-similar functions}\label{sec:DSS}

\begin{defn}\label{def-DSS}
Fix $\la>1$ and $a \in \R$. 
A scalar, vector or matrix function $g(x,t)$ defined on 
\EQ{
\QQ=\R^3\times (0,\infty)
}
 is called \textit{discretely self-similar} (DSS) with factor $\la$ and degree $a$, or $(\la,a)$-DSS, if 
\EQ{\label{DSS.defn}
g(x,t) = \la^{-a} g(\la^{-1} x,\la^{-2}t), \quad \forall (x,t)\in \QQ.
}
We simply say $\la$-DSS when $a=1$.
The same definition applies to those $g(x)$ defined on $\R^3$ that is independent of $t$.

We say $g$ is \textit{self-similar} of degree $a$ if \eqref{DSS.defn} is valid for all $\la>0$. In this case we can choose $\la=\sqrt t$ to get
\[
g(x,t) = \frac 1{t^{a/2}}\, g\bke{\frac x{\sqrt t}, 1}.
\]
\end{defn}

For example, for DSS solutions of the Navier-Stokes equations, 
the velocity $v$ has degree $a=1$ and the pressure $\pi$ has degree $a=2$.

When $g$ is $(\la,a)$-DSS,
by iterating the equality \eqref{DSS.defn} we get
$g(x,t) = \la^{-ak} g(\la^{-k} x,\la^{-2k}t)$ for any $k\in\ZZ$. So, for any $q\in[1,\infty)$, we have
\EQS{\label{eq-BT2019-3.6}
\int_{B_{\la^k}} |g(x,t)|^q\, dx &= \la^{(3-aq)k} \int_{B_1} |g(\td x, \la^{-2k}t)|^q\, d\td x,\\
\int_{B_{\la^k}} |g_0(x)|^q\, dx &= \la^{(3-aq)k} \int_{B_1} |g_0(\td x)|^q\, d\td x.
}
In the above equality, we can replace ball $B_r$ by annulus $B_{\la r} \setminus B_r$.

\begin{lem}\label{lem2.2}
$(i)$  If two functions $h(x,t)$, $g(x,t)$ on $\QQ$ are $(\la,a)$-DSS and $(\la,b)$-DSS, respectively, for some $\la>1$ and $a,b\in\R$, then their product $(hg)(x,t)$  is $(\la,a+b)$-DSS.

$(ii)$ Let $\eta\in C^\infty_c(\R^3)$ and $\eta_\ep(y) = (1/\ep^3)\eta(y/\ep)$ for $\ep>0$.
If a function $g(x,t)$ on $\QQ$ is $(\la,a)$-DSS, then so is
\[
(\eta_{\e \sqrt t}* g)(x,t)
=\int_{\R^3} \eta_{\e \sqrt t} (x-y)\, g(y,t)\,dy.
\]

\end{lem}

 Part (i) is immediate. For (ii), let $\tilde g(x,t)=(\eta_{\e\sqrt t}*g)(x,t)$. Since
\[
\eta_{\e\la^{-1}\sqrt t}(\la^{-1}w)=\la^3\eta_{\e\sqrt t}(w)\quad\text{and}\quad g(\la^{-1}z,\la^{-2}t)=\la^ag(z,t),
\]
the change of variables $y=\la^{-1}z$ gives $\tilde g(\la^{-1}x,\la^{-2}t)=\la^a\tilde g(x,t)$. For our applications, we may take  $h(x,t)=|x|^{-a}$.

\begin{lem}\label{lem2.3}
Let $T>0$, $q \ge 1$. Let $A_0=\{y\in\R^3:\la^{-1}\le|y|<1\}$.
If a function $g(x,t)$ on $\QQ$ is $(\la,a)$-DSS with $a>3+\frac 2q$, then
\[
\int_0 ^T\bke{\int_{|y|>\la^2} |g(y,t)| dy}^q dt \le C \int_0^T \int_{A_0} |g(y,t)|^q dy\, dt.
\]
\end{lem}
This lemma will be useful in bounding the far field part of the pressure. 

\begin{proof}
Let $A_k=\{x\in\R^3:\la^{k-1}\le|x|<\la^k\}$ for $k\in\ZZ$.
By Minkowski inequality,
\EQN{
&\bke{\int_0 ^T\bke{\int_{|y|>\la^2} |g(y,t)| dy}^q dt}^{1/q}
\le\int_{|y|>\la^2} \bke{ \int_0 ^T |g(y,t)|^q dt}^{1/q} dy
\\
&\qquad\qquad=\sum_{k \ge 3} \int_{A_k}\bke{\int_0 ^T \bke{\la^{-ak} |g(\la^{-k} y,\la^{-2k}t)|}^q dt}^{1/q} dy
\\
&\qquad\qquad=\sum_{k \ge 3} \la^{k(-a+3+2/q)} \int_{A_0} \bigg( \int_0 ^{T\la^{-2k}} |g(z,\tau)|^q d\tau\bigg)^{1/q} dz.
}
The sum converges if $a>3+2/q$ and is bounded by $C (\int_0 ^T\int_{A_0}  |g(z,\tau)|^q dz d\tau)^{1/q}$.
\end{proof}

\begin{lem}\label{lem2.4}
Let $0<q<\infty$ and $0<r\le \infty$.
If a function $g(x,t)$ on $\QQ$ is $(\la,a)$-DSS and $g \in L^\infty(1,\la^2; L^{q,r}(\R^3))$, then
\[
\norm{g(\cdot,t)}_{L^{q,r}(\R^3)} \le C \norm{g}_{L^\infty(1,\la^2; L^{q,r}(\R^3))} t^{\frac{-a+3/q}2} ,
\quad \text{for a.e.\ } t>0.
\]
\end{lem}

\begin{proof}
Choose $k \in \ZZ$ so that $\tau=\la^{-2k}t \in [1,\la^2)$. From
$g(x,t) = \la^{-ak} g(\la^{-k} x,\tau)$, we have
\[
\norm{g(\cdot,t)}_{L^{q,r}(\R^3)} =  \la^{-ak} (\la^{-k})^{-\frac 3q} \norm{g(\cdot,\tau)}_{L^{q,r}(\R^3)} 
= (\tau/t)^{\frac 12(a-\frac 3q)} \norm{g(\cdot,\tau)}_{L^{q,r}(\R^3)} .
\]
The lemma follows with $C=\max(1, \la ^{a-\frac 3q} )$.
\end{proof}

\section{A limiting pressure formula for DSS solutions}\label{sec:limiting-pressure}

In this section we will prove that, under certain conditions, the limiting pressure distribution of an approximation scheme for either MHD or OB inherits 
the structure of the approximate pressure distributions. These
results will be applied in Sections \ref{S:L3w-data} and \ref{S:L2loc-data}.

\begin{lem}[MHD]\label{BT2019-lem-2.1}
Fix $\la>1$ and $T>0$. 
Assume $v_k, \td v_k, b_k$ and $\td b_k$ are divergence-free, $\la$-DSS vector fields and that there   exist distributions $\pi_k$ and $\pi^b_k$  so that the following conditions are satisfied:\begin{itemize}
\item $v_k, \td v_k, b_k, \td b_k$,  $\pi_k$ and $\pi^b_k$  solve the system
\EQN{
\pd_tv_k - \De v_k + \td v_k\cdot\nb v_k - \td b_k\cdot\nb b_k + \nb\pi_k &= 0,\quad
(x,t)\in\R^3\times[0,T],\\
\pd_tb_k - \De b_k + \td v_k\cdot\nb b_k - \td b_k\cdot\nb v_k  +\nb\pi^b_k  &= 0,\quad
(x,t)\in\R^3\times[0,T].
}

\item $v_k, \td v_k, b_k$, and $\td b_k$ are uniformly bounded in $L^\infty(0,T;L^2(B_1))\cap L^2(0,T; H^1(B_1))$ over all $k\in\NN$.

\item The distribution $\pi_k$ satisfies in $L^{3/2}((0,T)\times B_\la)$ the formula
\EQN{
\pi_k(x,t) &= -\frac13(\td v_k\cdot v_k - \td b_k\cdot b_k)(x,t)\\
&\quad + \lim_{\de\to0} \int_{| x-y|>\de} K_{ij}(x-y) \bkt{(\td v_k)_i(v_k)_j - (\td b_k)_i(b_k)_j}(y,t)\, dy
}

\item There exists a $\la$-DSS solution $(v,b)$ in $L^\infty(0,T;L^2(B_1))\cap L^2(0,T;H^1(B_1))$ with pressure $\pi$ in $L^{3/2}(0,T; L^{3/2}(B_\la))$ so that
\EQN{
\text{$v_k, \td v_k\to v$ and $b_k, \td b_k\to b$ weakly}& \text{ in $L^2(0,T; H^1(B_1))$,}\\
\text{$v_k, \td v_k\to v$ and $b_k, \td b_k\to b$}& \text{ in $L^2(0,T; L^2(B_1))$,}\\
\text{$\pi_k\to \pi$ weakly}& \text{ in $L^{3/2}(0,T; L^{3/2}(B_\la))$.}
}
\end{itemize}
Then, 
 the pressure $\pi$ satisfies in $L^{3/2}((0,T)\times B_\la)$ the formula
\EQ{\label{eq-BT2019-2.2}
\pi(x,t) = -\frac13(|v|^2 - |b|^2)(x,t) + \lim_{\de\to0} \int_{| x-y|>\de} K_{ij}(x-y) \bke{v_iv_j - b_ib_j}(y,t)\, dy.
}

\end{lem}

\begin{proof}
The proof of the lemma is an adaption of the proof of \cite[Lemma 2.1]{BT-AnalPDE2019} for the Navier-Stokes equations to the MHD equations.

Note that since $v_k, \td v_k, b_k$, and $\td b_k$ are all uniformly bounded in $L^\infty(0,T; L^2(B_1))\cap L^2(0,T; H^1(B_1))$,  hence in $L^{10/3}(B_1\times(0,T))$ by the Sobolev embedding, the convergence in $L^2(0,T; L^2(B_1))$ and interpolation  imply 
\[
v_k,\td v_k\to v\ \text{ and }\ b_k,\td b_k\to b\quad\text{ in }L^3(0,T;L^3(B_1)).
\]
It also shows that $v_k, \td v_k, b_k$, and $\td b_k$ are all uniformly bounded in $L^3(0,T;L^3(B_1))$.

Let
\EQN{
\pi_k^1(x,t) &=  -\frac13\bkt{\td v_k\cdot v_k - \td b_k\cdot b_k}(x,t),\\
\pi_k^2(x,t) &= \lim_{\de\to0} \int_{|y|<\la^2,\,|x-y|>\de} K_{ij}(x-y) \bkt{(\td v_k)_i(v_k)_j - (\td b_k)_i(b_k)_j}(y,t)\, dy,\\
\pi_k^3(x,t) &= \int_{|y|\ge\la^2} K_{ij}(x-y) \bkt{(\td v_k)_i(v_k)_j - (\td b_k)_i(b_k)_j}(y,t)\, dy.
}
Also let
\EQN{
\pi^1(x,t) &=  -\frac13\bkt{|v|^2 - |b|^2}(x,t),\\
\pi^2(x,t) &= \lim_{\de\to0} \int_{|y|<\la^2,\,|x-y|>\de} K_{ij}(x-y) \bke{v_iv_j - b_ib_j}(y,t)\, dy,\\
\pi^3(x,t) &= \int_{|y|\ge\la^2} K_{ij}(x-y) \bke{v_iv_j - b_ib_j}(y,t)\, dy.
}
That $\pi_k^3$ and $\pi^3$ are defined in $L^{3/2}((0,T)\times B_\la)$ will
follow from the same argument to bound $\pi_k^3-\pi^3$ below.

Since,  by DSS rescaling, $v_k,\td v_k\to v$ and $b_k,\td b_k\to b$ in $L^3(0,T;L^3(B_\la))$, we have
\[
\pi_k^1\to\pi^1\quad\text{ in }L^{3/2}(0,T;L^{3/2}(B_\la)).
\]

Suppressing the dependence on $k$, let
\EQN{
h_{ij}(y,t) &= (\td v_k)_i(v_k)_j - v_iv_j - (\td b_k)_i(b_k)_j + b_ib_j\\
&= \bkt{(\td v_k)_i\bke{(v_k)_j-v_j} + \bke{(\td v_k)_i-v_i}v_j - (\td b_k)_i\bke{(b_k)_j-b_j} - \bke{(\td b_k)_i-b_i}b_j}(y,t).
}
Using the Calder\'on--Zygmund theory we have
\EQS{\label{eq-BT2019-2.3}
\int_0^T\int_{B_\la}& |\pi_k^2(x,t) - \pi^2(x,t)|^{3/2}\, dxdt\\
&\le C\int_0^T\int_{B_{\la^2}} |h_{ij}(x,t)|^{3/2}\, dxdt\\
&\le C \bke{\int_0^T\int_{B_{\la^2}} |\td v_k|^3 + |v|^3 + |\td b_k|^3  + |b|^3 \, dxdt}^{1/2} \cdot
\bke{\int_0^T\int_{B_{\la^2}} F_k\, dxdt}^{1/2},
}
where $F_k = |\td v_k-v|^3 + |v_k-v|^3 + |\td b_k-b|^3  + |b_k-b|^3$.
Rescaling gives
\EQN{
\int_0^T\int_{B_{\la^2}} F_k (x,t)\, dxdt &= \la^4\int_0^{T\la^{-4}}\int_{B_1} F_k(z,\tau)\, dzd\tau
}
for the obvious choice of $z$ and $\tau$.
Since the right-hand side of the equation above vanishes as $k\to\infty$, 
we conclude that $\pi_k^2$ converges to $\pi^2$ in $L^{3/2}(0,T;L^{3/2}(B_\la))$.

For the convergence of $\pi_k^3$ to $\pi^3$, let
\[
p_k(x,t) = \pi_k^3(x,t) - \pi^3(x,t) = \int_{|y|\ge\la^2} K_{ij}(x-y)h_{ij}(y,t)\, dy.
\]
For $x\in B_\la$, we have $|x-y|\sim |y|$ for $|y|\ge\la^2$, and hence
$|p_k(x,t)|\le C\int_{|y|\ge\la^2} |y|^{-3} |h(y,t)|\, dy$. The integrand is $\la$-DSS with degree 5 by Lemma \ref{lem2.2}. By Lemma \ref{lem2.3},
\EQN{
&\int_0^T\!\!\int_{B_\la} |p_k(x,t)|^{3/2}\,dxdt 
\le C\int_0^T\bigg(\int_{|y|\ge\la^2}  |y|^{-3} |h(y,t)|\bigg)^{3/2}dt 
\\
&\le C\int_0^T\int_{A_0} \bke{ |y|^{-3} |h(y,t)| }^{3/2}dy\,dt 
\le C\int_0^T\int_{B_1} \bke{ |h(y,t)| }^{3/2}dy\,dt.
}
Therefore, this term is bounded as \eqref{eq-BT2019-2.3}.

We have now shown that $\pi_k$ converges weakly in $L^{3/2}(0,T; L^{3/2}(B_\la))$ to both $\pi^1 + \pi^2 + \pi^3$ and $\pi$, implying that $\pi = \pi^1 + \pi^2 + \pi^3$ as distributions.
In other words, $\pi$ satisfies \eqref{eq-BT2019-2.2} in $L^{3/2}((0,T)\times B_\la)$.
\end{proof}

\begin{lem}[OB]\label{BT2019-lem-2.1-OB}
Fix $\la>1$ and $T>0$. 
Assume $f=f_1+f_2$, $f_1,f_2$ are 
$(\la,3)$-DSS vector fields, $f_j(x,t) = \la^3f_j(\la x,\la^2t)$, $j=1,2$,
\[
 f_1 \in L^\infty(1,\la^2; L^{6/5,2}(\R^3)),  \quad f_2\in L^\infty(1,\la^2; L^2(\R^3)),\quad
\div f_2=0.
\]
Assume $v_k$ and $\td v_k$ are sequences of divergence-free, $\la$-DSS vector fields, and $\th_k$ is a sequence of $\la$-DSS functions, and that there exists a distribution $\pi_k$ so that the following conditions are satisfied:
\begin{itemize}
\item $v_k, \td v_k, \th_k$, and $\pi_k$ solve the system
\EQN{
\pd_tv_k - \De v_k + \td v_k\cdot\nb v_k + \nb\pi_k &= \th_k\nb G + f,\quad
(x,t)\in\R^3\times[0,T],\\
\pd_t\th_k - \De \th_k + \td v_k\cdot\nb \th_k &= 0,\quad
(x,t)\in\R^3\times[0,T].
}

\item $v_k, \td v_k$, and $\th_k$ are uniformly bounded in $L^\infty(0,T;L^2(B_1))\cap L^2(0,T; H^1(B_1))$ over all $k\in\NN$.

\item  The distribution $\pi_k$ satisfies in $L^{3/2}((0,T)\times B_\la)$ the formula
\EQS{\label{eq-BT2019-2.1-OB}
\pi_k(x,t) &= -\frac13\bkt{\td v_k\cdot v_k}(x,t)
 + \lim_{\de\to0} \int_{| x-y|>\de} K_{ij}(x-y) (\td v_k)_i(y,t)(v_k)_j (y,t)\, dy\\
 &\quad - \int \nb_x\bke{\frac1{4\pi|x-y|}}\cdot(\th_k\nb G + f_1)(y,t)\, dy.
}

\item There exists a $\la$-DSS solution $(v,\th)$ in $L^\infty(0,T;L^2(B_1))\cap L^2(0,T;H^1(B_1))$ with pressure $\pi$ in $L^{3/2}(0,T; L^{3/2}( B_\la))$ so that
\EQN{
\text{$v_k, \td v_k\to v$ and $\th_k\to \th$ weakly}& \text{ in $L^2(0,T; H^1(B_1))$,}\\
\text{$v_k, \td v_k\to v$ and $\th_k\to \th$}& \text{ in $L^2(0,T; L^2(B_1))$,}\\
\text{$\pi_k\to \pi$ weakly}& \text{ in $L^{3/2}(0,T; L^{3/2}(B_\la))$.}
}
\end{itemize}
Then the pressure $\pi$ satisfies in $L^{3/2}((0,T)\times B_\la)$ the formula
\EQS{\label{eq-BT2019-2.2-OB}
\pi(x,t) &= -\frac13|v|^2(x,t) + \lim_{\de\to0} \int_{| x-y|>\de} K_{ij}(x-y)v_i(y,t)v_j(y,t)\, dy \\
&\quad - \int \nb_x\bke{\frac1{4\pi|x-y|}}\cdot(\th\nb G + f_1)(y,t)\, dy .
}

\end{lem}

\begin{proof} 
\emph{Reduction to $\div f_1=0$.} Let $\mathbb P$ be the Helmholtz projection and
\[
q_1=-(-\De)^{-1}\div f_1=-\int \nb_x\bke{\frac1{4\pi|x-y|}}\cdot f_1(y,t)\, dy,
\]
which is the term with $f_1$ in \eqref{eq-BT2019-2.1-OB} and \eqref{eq-BT2019-2.2-OB}, so that $f_1=\mathbb Pf_1+\nb q_1$.
Since the Riesz transforms are bounded on $L^{6/5,2}(\R^3)$ and commute with dilations, $\mathbb Pf_1$ is a divergence-free $(\la,3)$-DSS vector field with $\norm{\mathbb Pf_1}_{L^\infty(1,\la^2;L^{6/5,2})}\le CM_1$, where $M_1=\norm{f_1}_{L^\infty(1,\la^2;L^{6/5,2})}$.
By O'Neil's inequality $L^{3/2,\infty}*L^{6/5,2}\subset L^2$ and Lemma \ref{lem2.4},
\[
\norm{q_1(t)}_{L^2(\R^3)}\le C\norm{f_1(t)}_{L^{6/5,2}}\le C M_1 t^{-1/4},
\]
and $t^{-1/4}\in L^{3/2}(0,T)$, so $q_1\in L^{3/2}(0,T;L^{3/2}(B_\la))$.
Hence $(v_k,\td v_k,\th_k,\pi_k-q_1)$ and $(v,\th,\pi-q_1)$ satisfy the assumptions of the lemma with $f_1$ replaced by $\mathbb Pf_1$, whose term in \eqref{eq-BT2019-2.1-OB} vanishes, and \eqref{eq-BT2019-2.2-OB} for $\pi$ is equivalent to the same formula for $\pi-q_1$ with $\mathbb Pf_1$ in place of $f_1$.
Moreover, since $\div v=0$, we have $\nb q_1\cdot v\phi+q_1v\cdot\nb\phi=\div(q_1v\phi)$, so the local energy inequalities \eqref{eq-OB-LEI} for $(v,\th,\pi)$ with force $f$ are equivalent to those for $(v,\th,\pi-q_1)$ with force $\mathbb Pf_1+f_2$.
Therefore we may assume $\div f_1=0$, so that the term with $f_1$ is absent in \eqref{eq-BT2019-2.1-OB} and \eqref{eq-BT2019-2.2-OB}. This reduction will also be used in the proof of Proposition \ref{BT2019-prop-3.1-OB}.

 By the same argument in the beginning of the proof of Lemma \ref{BT2019-lem-2.1}, we have
\[
v_k,\td v_k\to v\text{ and }\th_k\to \th\quad\text{ in }L^q(0,T;L^q(B_1)),
\]
for any $q<10/3$. We take $q=3$ in Lemma \ref{BT2019-lem-2.1}, but now we need $3<q<10/3$, as $q=3$ corresponds to $r=3/2$ below, which is critical for $|x|^{-2}\in L^r_\loc$.

The convergences of the first two terms of \eqref{eq-BT2019-2.1-OB} to the corresponding terms in \eqref{eq-BT2019-2.2-OB} follow from the same argument as in the proof of Lemma \ref{BT2019-lem-2.1} (see also the proof of \cite[Lemma 2.1]{BT-AnalPDE2019}).
It suffices to show that
\[
\int \nb_x\bke{\frac1{4\pi|x-y|}} \cdot(\th_k\nb G)(y,t)\, dy \to 
\int \nb_x\bke{\frac1{4\pi|x-y|}} \cdot(\th\nb G)(y,t)\, dy
\]
in $L^{3/2}(0,T; L^{3/2}(B_\la))$.
We split their difference into local and nonlocal parts as follows:
\EQN{
p_{\rm near}(x,t) + p_{\rm far}(x,t)
&:= 
\int_{|y|\le\la^2} \nb_x\bke{\frac1{4\pi|x-y|}}\cdot(\wth_k\nb G)(y,t)\, dy\\
&\quad + \int_{|y|>\la^2} \nb_x\bke{\frac1{4\pi|x-y|}}\cdot(\wth_k\nb G)(y,t)\, dy,
}
where $\wth_k = \th_k -\th$.

To estimate $p_{\rm near}$, choose $r \in ( 60/41,3/2)$ by $\frac 2r=\frac53-\frac1q$.
By Young's convolution inequality and H\"older's inequality,
\EQN{
\norm{p_{\rm near}(\cdot,t)}_{L^{3/2}(B_{\la})} 
&\lec \norm{|\cdot|^{-2}}_{L^r(B_{\la+\la^2})} \norm{\wth_k(\cdot,t)\nb G}_{L^{3r/(5r-3)}(B_{\la^2})} \\
&\lec 1\cdot \norm{\wth_k(\cdot,t)}_{L^q(B_{\la^2})} \norm{|\cdot|^{-2}}_{L^r(B_{\la^2})}
=C\norm{\wth_k(\cdot,t)}_{L^q(B_{\la^2})}.
}
It follows that
\EQS{\label{eq3.5}
\int_0^T \norm{p_{\rm near}(\cdot,s)}_{L^{3/2}(B_\la)}^{3/2}\, ds
&\le C(\la) \int_0^T \norm{\wth_k(\cdot,s)}_{L^q(B_{\la^2})}^{3/2}\, ds\\
&\le C(\la,T) \bke{\int_0^T \norm{\wth_k(\cdot,s)}_{L^q(B_{\la^2})}^q ds}^{\frac3{2q}} ,
}
which goes to $0$ as $k\to\infty$.

Now we estimate the term $p_{\rm far}$. 
For $x\in B_\la$, we have $|x-y|\sim |y|$ for $|y|\ge\la^2$, and hence
$|p_{\rm far}(x,t)|\le C\int_{|y|\ge\la^2} |y|^{-4} |\wth_k(y,t)|\, dy$. The integrand is $\la$-DSS with degree 5 by Lemma \ref{lem2.2}. By Lemma \ref{lem2.3},
\EQN{
&\int_0^T\!\!\int_{B_\la} |p_{\rm far}(x,t)|^{3/2}\,dxdt 
\le C\int_0^T\bigg(\int_{|y|\ge\la^2}  |y|^{-4}  |\wth_k(y,t)| \bigg)^{3/2}dt 
\\
&\le C\int_0^T\int_{A_0} \bke{ |y|^{-4}  |\wth_k(y,t)| }^{3/2}dy\,dt 
\le C\int_0^T\int_{B_1} \bke{  |\wth_k(y,t)| }^{3/2}dy\,dt,
}
which goes to $0$ as $k\to\infty$.
This completes the proof of the lemma.
\end{proof}

\section{Properties of DSS solutions with data in $L^3_w$}\label{S:L3w-data}

In this section we will prove the local energy estimates for $\la$-DSS solutions of MHD and OB with weak $L^3$ initial data constructed in 
\cite{Lai-JMFM2019,Tsai-OB} 
solely in terms of the $L^2$ norm of the initial data in $B_\la$ (and, for OB, of the force). 
These estimates heavily rely on the DSS property, and are not available for general solutions with non-DSS weak $L^3$ data. They will provide the a priori bound when we construct solutions in  Section \ref{sec:construction}  for initial data in local $L^2$ by approximating the data with data in weak $L^3$. In the proof we must work with a quantity that is continuous in time. Hence, we need to work at the level of a
mollified approximation scheme.

\begin{prop}[MHD]\label{BT2019-prop-3.1}
Fix $\la>1$. Assume $v_0, b_0\in L^3_w(\R^3)$ are $\la$-DSS and divergence-free, and $(v,b)$ is a local Leray solution of the MHD equations \eqref{eq-MHD} evolving from $(v_0, b_0)$ constructed  as in \cite{Lai-JMFM2019}, with the modifications of Remarks \ref{approx-scheme} and \ref{rem-pib}, and $\pi$ is its associated pressure. 
Let $\al_0=\norm{v_0}_{L^2(B_\la)}^2 + \norm{b_0}_{L^2(B_\la)}^2$.
Then, there exist positive $T=T(\al_0,\la)$ and $C(\al_0,\la)$ independent of $\norm{v_0}_{L^2_\uloc}$, $\norm{b_0}_{L^2_\uloc}$, $\norm{v_0}_{L^3_w}$, and $\norm{b_0}_{L^3_w}$ so that 
\EQ{
\esssup_{0\le t\le T} \int_{B_1} \left(|v(x,t)|^2 + |b(x,t)|^2\right) dx + \int_0^T\int_{B_1} \left(|\nb v|^2 + |\nb b|^2\right) dxdt < C(\al_0, \la),
}
and 
\EQ{\label{eq-BT2019-3.2}
\int_0^T\int_{B_\la} |\pi(x,t)|^{3/2}\, dxdt < C(\al_0, \la).
}
Moreover,  the pressure  satisfies the formula
\EQ{\label{eq-BT2019-3.3}
\pi(x,t) = -\frac13(|v|^2-|b|^2)(x,t) + \lim_{\de\to0} \int_{| x-y|>\de} K_{ij}(x-y)
\bke{v_iv_j - b_ib_j}(y,t)\, dy
}
in $L^{3/2}(B_\la\times(0,T))$.
 Furthermore, the second equation of \eqref{eq-MHD} holds in the sense of distributions without a pressure term, and $(v,b)$ satisfies the local energy inequality \eqref{eq-MHD-LEI} with initial data $(v_0,b_0)$. 
\end{prop}

\begin{proof}
The proof of the proposition is an adaption of the proof of \cite[Proposition 3.1]{BT-AnalPDE2019} for the Navier-Stokes equations to the MHD equations.

Since $(v,b)$ is  constructed as in \cite{Lai-JMFM2019}, its image $(V,B)$ in the similarity variables \eqref{eq-self-similar-variables} solves the time-periodic Leray equations \eqref{eq-MHD-Leray} and is the limit of a mollified approximation scheme \cite[(2.16)]{Lai-JMFM2019}.
In particular, for each $\ep>0$, there exists a time-periodic solution $(V_\ep, B_\ep)$ to the problem (regularized version of \eqref{eq-MHD-Leray}) in $\R^3 \times \R$ 
\EQS{\label{eq4.4}
\pd_sV_\ep - \De V_\ep - \frac12 V_\ep - \frac12 y\cdot\na V_\ep + (\eta_\ep * V_\ep)\cdot\na V_\ep - (\eta_\ep * B_\ep)\cdot\na B_\ep + \nb P_\ep&= 0,\\
\pd_sB_\ep - \De B_\ep - \frac12 B_\ep - \frac12 y\cdot\na B_\ep + (\eta_\ep * V_\ep)\cdot\na B_\ep - (\eta_\ep * B_\ep)\cdot\na V_\ep + \nb P^B_\ep &= 0,
}
where $\eta_\ep(y) = (1/\ep^3)\eta(y/\ep)$ and $\eta$ is a fixed mollifier satisfying 
\EQ{\label{mollifier}
\eta \in C^\infty_c(\R^3), \quad \eta\ge0, \quad \text{radial}, \quad \eta(y)=0 \text{ if } |y|\ge1, \quad \int _{\R^3}\eta(y)\, dy = 1.
}
Moreover, a subsequence $V_{\ep_j},B_{\ep_j}$ and $P_{\ep_j}$ converges weakly in the energy class and $L^{5/3}_{x,t}$ to a weak solution $(V,B,P)$ of \eqref{eq-MHD-Leray}
as ${\ep_j} \to 0$, see \cite[(2.37)--(2.38)]{Lai-JMFM2019}.

A remark is in order: the approximation scheme in \cite{Lai-JMFM2019} needs to be slightly changed to get \eqref{eq4.4} so that we have a simpler local energy equality
\eqref{eq-MHD-LEE-epsilon} below. See Remark \ref{approx-scheme}.

 A pressure term $\nb P^B_\ep$ is needed in the second equation of \eqref{eq4.4}: because of the mollification, $(\eta_\ep * V_\ep)\cdot\na B_\ep - (\eta_\ep * B_\ep)\cdot\na V_\ep$ is not divergence free, and in \cite[(2.17)]{Lai-JMFM2019} the equation for $B_\ep$ is tested only with divergence-free fields. We will show that the corresponding term disappears in the limit $\ep\to0$.

Applying \eqref{eq-MHD-SS-transform} and \eqref{eq-self-similar-variables}, we obtain a $\la$-DSS solution $(v_\ep, b_\ep,\pi_\ep,\pi^b_\ep)$ to the system
\EQS{\label{eq-BT2019-3.5}
\pd_tv_\ep(x,t) - \De v_\ep(x,t) + (\eta_{\ep\sqrt{t}} * v_\ep)\cdot\nb v_\ep(x,t) - (\eta_{\ep\sqrt{t}} * b_\ep)\cdot\nb b_\ep(x,t) + \nb\pi_\ep(x,t) &= 0,\\[1mm]
\pd_tb_\ep(x,t) - \De b_\ep(x,t) + (\eta_{\ep\sqrt{t}} * v_\ep)\cdot\nb b_\ep(x,t) - (\eta_{\ep\sqrt{t}} * b_\ep)\cdot\nb v_\ep(x,t) + \nb\pi^b_\ep(x,t)&= 0,
}
and a DSS solution $(v, b,\pi)$ to the MHD system \eqref{eq-MHD}. 

By the convergence properties of $(V_{\ep_j},B_{\ep_j},P_{\ep_j})$ to $(V,B,P)$
\cite[(2.37)--(2.38)]{Lai-JMFM2019} and discretely self-similar scaling (to extend the estimates down to $t=0$), it follows that for all $T>0$ and all compact sets $E$ of $\R^3$, (up to the subsequence $\ep_j$)
\EQS{\label{eq4.7}
v_\ep\wkto v,\ b_\ep\wkto b\ \text{ weakly }&\text{in $L^2(0,T;H^1(E))$},\\
v_\ep\to v,\ b_\ep\to b\ \text{ strongly }&\text{in $L^2(0,T;L^2(E))$},\\
\pi_\ep\wkto \pi\ \text{ weakly }&\text{in $L^{5/3}(0,T;L^{5/3}(E))$},\\
v_\ep(s)\wkto v(s),\ b_\ep(s)\wkto b(s)\ \text{ weakly }&\text{in $L^2(E)$ for all $s\in[0,T]$}.
}
Note also that $v_\ep(t)\to v_0$ and $b_\ep(t)\to b_0$  as $t\to0^+$  in $L^2_\loc$; i.e., the mollification does not affect the initial data.
Furthermore, because each $(v_\ep, b_\ep)$ is smooth on $\R^3\times(0,\infty)$ and right continuous in $L^2_\loc$ at $t=0$, it follows that
\[
\al_\ep(t) = \int_{B_1} \left(|v_\ep(x,t)|^2 + |b_\ep(x,t)|^2\right) dx
\]
and 
\[
\td\al_\ep(t) = \sup_{0\le\tau\le t} \al_\ep(\tau)
\]
are continuous in $t$. 
Our goal is to establish local-in-time a priori bounds for $\al_\ep(t)$ that are independent of $\ep$. 
Note that $(v_\ep, b_\ep)$ satisfies the local energy equality; i.e., 
\EQS{\label{eq-MHD-LEE-epsilon}
\int &\left(|v_\ep|^2 + |b_\ep|^2\right)\phi(x,t)\, dx + 2\int_0^t\int \left(|\nb v_\ep|^2+|\nb b_\ep|^2\right)\phi\, dxds\\
&= \int \left(|v_0|^2 + |b_0|^2\right)\phi(x,0)\, dx
 + \int_0^t \int \left(|v_\ep|^2+|b_\ep|^2\right)\left(\pd_s\phi+\De\phi\right) dxds\\
&\quad + \int_0^t\int \left(|v_\ep|^2+|b_\ep|^2\right) \left((\eta_{\ep\sqrt{s}} * v_\ep)\cdot\nb\phi\right) dxds
 + 2 \int_0^t\int \pi_\ep (v_\ep\cdot\nb\phi) dxds\\
&\quad - 2\int_0^t\int (v_\ep\cdot b_\ep)\left((\eta_{\ep\sqrt{s}} * b_\ep) \cdot\nb\phi\right) dxds
 + 2 \int_0^t\int \pi^b_\ep (b_\ep\cdot\nb\phi) dxds
}
for any nonnegative $\phi\in C^\infty_c(\R^3\times[0,\infty))$.
Fix $\chi\in C^\infty(\R)$ with $\chi(t) = 1$ if $t\le 1$, $0\le\chi(t)\le1$ for $t\in[1,\la]$, and $\chi(t) = 0$ if $t\ge\la$.
We now fix $\phi$ in \eqref{eq-MHD-LEE-epsilon} as 
\[
\phi(x,t) = \chi^2(|x|)\cdot\chi(t).
\]

We will estimate the terms on the right-hand side of \eqref{eq-MHD-LEE-epsilon} for $0<t\le 1$, and we can treat $\phi$ as $t$-independent from now on.
The first term is bounded by $\al_0$. Using the scaling properties \eqref{eq-BT2019-3.6},
 the second term is bounded by 
\[
C\int_0^t\int_{B_\la} \left(|v_\ep|^2 + |b_\ep|^2\right) dxds
\le C\la^3 \int_0^{t/\la^2} \int_{B_1} \left(|v_\ep|^2 + |b_\ep|^2\right) dxds
\le C(\la) \int_0^t \td\al_\ep(s)\, ds.
\]
For the third term, we begin by using Young's inequality to obtain
\EQN{
\int_0^t\int &\left(|v_\ep|^2 + |b_\ep|^2\right) \left((\eta_{\ep\sqrt{s}} * v_\ep)\cdot\nb\phi\right) dxds\\
&\le C \int_0^t\int_{B_\la} \left(|v_\ep|^3 + |b_\ep|^3\right) dxds + C \int_0^t\int_{B_\la} |(\eta_{\ep\sqrt{s}} * v_\ep)|^3\, dxds.
}
Rescaling the unmollified term and making a change of variables yields the estimate
\EQS{
\int_0^t\int_{B_\la}& \left(|v_\ep|^3 + |b_\ep|^3\right) dxds\\
&\le C(\la) \int_0^{t/\la^2} \int_{B_1} \left(|v_\ep|^3 + |b_\ep|^3\right) dyd\tau
\le C(\la) \int_0^t \int \left(|v_\ep|^3 + |b_\ep|^3\right)\phi^{3/2} dxds.
}
For the term involving the mollifier, note that $\supp\eta\subset B_1$ by \eqref{mollifier}.
By taking $\ep$ sufficiently small we can ensure that $\supp\eta_{\ep\sqrt{s}}\subset B_{\la-1}$ whenever $s<1$.
Note that $\la^k+(\la-1)\le\la^{k+1}$ for all $k\ge0$.
Thus, for $x\in B_\la$,
\[
|(\eta_{\ep\sqrt{s}} * v_\ep)(x,s)| 
\le 
(\eta_{\ep\sqrt{s}} * (\chi_{B_{\la^2}}|v_\ep|))(x,s)
\]
whenever $\ep$ is sufficiently small and $s<1$.
Therefore, under the same assumptions and after rescaling we see that, for any $1<q<\infty$,
\EQ{\label{eq-BT2019-3.8}
\norm{(\eta_{\ep\sqrt{s}} * v_\ep)(s)}_{L^q(B_\la)}
\le C(q,\eta) \norm{v_\ep(s)}_{L^q(B_{\la^2})}
\le C(q,\eta,\la) \norm{v_\ep(\la^{-4}s)}_{L^q(B_1)},
}
where $C$ is independent of $s$ and $\ep$.

Using standard inequalities and \eqref{eq-BT2019-3.8} with $q=3$ thus leads to the estimate
\EQ{\label{eq-BT2019-3.9}
\int_0^t\int \left(|v_\ep|^2 + |b_\ep|^2\right) \left((\eta_{\ep\sqrt{s}} * v_\ep)\cdot\nb\phi\right) dxds
\le C(\eta,\la)\int_0^t \int \left(|v_\ep|^3 + |b_\ep|^3\right)\phi^{3/2} dxds.
}

By the Gagliardo--Nirenberg inequality and rescaling \eqref{eq-BT2019-3.6}, we have, for any $s>0$, that 
\EQN{
\norm{\phi^{1/2}v_\ep(s)}_{L^3} 
&\le C\norm{\nb(\phi^{1/2}v_\ep)}_{L^2}^{1/2} \norm{\phi^{1/2}v_\ep}_{L^2}^{1/2}(s)\\
&\le C(\la) \left(\td\al_\ep(s)^{1/2} + \norm{\phi^{1/2}\nb v_\ep(s)}_{L^2}\right)^{1/2} \left(\td\al_\ep(s)\right)^{1/4}.
}
Hence, for any $\ga>0$,
\[
\norm{\phi^{1/2}v_\ep(s)}_{L^3}^3 
\le C(\la) \left(\ga^{-3}\td\al_\ep(s)^3 + \ga\td\al_\ep(s) + \ga\norm{\phi^{1/2}\nb v_\ep}_{L^2}^2\right).
\]
Similarly, we have, for any $\ga>0$,
\[
\norm{\phi^{1/2}b_\ep(s)}_{L^3}^3 
\le C(\la) \left(\ga^{-3}\td\al_\ep(s)^3 + \ga\td\al_\ep(s) + \ga\norm{\phi^{1/2}\nb b_\ep}_{L^2}^2\right).
\]
Thus,
\EQS{\label{eq-BT2019-3.10}
\int_0^t&\int \left(|v_\ep|^2 + |b_\ep|^2\right) \left((\eta_{\ep\sqrt{s}} * v_\ep)\cdot\nb\phi\right) dxds\\
&\quad \le C(\ga) \int_0^t \left(\td\al_\ep(s)^3 + \td\al_\ep(s)\right) ds 
+ C\ga\int_0^t\int \left(|\nb v_\ep|^2 + |\nb b_\ep|^2\right)\phi\, dxds.
}
Provided $\ga$ is small enough, the gradient term can be absorbed into the left-hand side of \eqref{eq-MHD-LEE-epsilon}.

We next estimate the pressure term.
Let $w_\ep = v_\ep -  v_L$ and $c_\ep = b_\ep -  b_L$, where $v_L(x,t) = e^{t\De}v_0$ and $b_L(x,t) = e^{t\De}b_0$.
We have
\EQN{
\pd_tw_\ep - \De w_\ep + \nb\pi_\ep = h,&\qquad\div w_\ep = 0,
}
where $h_i = -\pd_jH_{ji}$ 
with
\EQN{
H &= (\eta_{\ep\sqrt{t}} * v_\ep)\otimes v_\ep - (\eta_{\ep\sqrt{t}} * b_\ep)\otimes b_\ep\\[1mm]
&= (\eta_{\ep\sqrt{t}} * w_\ep + \eta_{\ep\sqrt{t}} * v_L)\otimes (w_\ep + v_L) - (\eta_{\ep\sqrt{t}} * c_\ep + \eta_{\ep\sqrt{t}} * b_L)\otimes (c_\ep + b_L).
}
For $0<t_1<t_2<\infty$, we have
\EQN{
&v_L, b_L\in C([t_1,t_2]; L^4(\R^3)\cap L^\infty(\R^3)),\\
&w_\ep, c_\ep\in L^\infty(t_1,t_2;L^2(\R^3))\cap L^2(t_1,t_2;L^6(\R^3))\subset L^4(t_1,t_2;L^3(\R^3)).
}
(The properties of $w_\ep$ and $c_\ep$ follow from the construction in \cite{Lai-JMFM2019}.) 
By Young's convolution inequality,
\EQS{
&\norm{H}_{L^2(t_1,t_2;L^2)} 
\\
&\lec \norm{\eta_{\ep\sqrt{t}}}_{L^\infty(t_1,t_2;L^{\frac65}\cap L^1)} \left(\big\| |w_\ep| + |c_\ep|\big\|_{L^4(t_1,t_2;L^3(\R^3))} ^2 + \big\| |v_L|+|b_L|\big\|_{L^4(\R^3\times[t_1,t_2])}^2 \right).
}
Since $h\in L^2(t_1,t_2;H^{-1})$, \cite[Lemma A.2]{CKN-CPAM1982} implies $w_\ep\in C([t_1,t_2];L^2)$. 

Consider the following nonstationary Stokes system with forcing $h$:
\[
\pd_tV - \De V + \nb P = h,\quad
\div V = 0,
\]
with initial data $V_0=w_\ep(t_1)\in L^2(\R^3)$. 
It is well known that if $h\in  L^2(t_1, t_2;H^{-1})$ and $V_0\in L^2$, then there exists a unique $V\in C_w([t_1,t_2];L^2(\R^3))\cap L^2(t_1,t_2; H^1(\R^3))$ and unique $\na P$ solving the nonstationary Stokes system above (see \cite[p.\,1107--1108]{BT1-AHP}, where $h\in L^\infty(t_1,t_2;H^{-1})$ is assumed; the same energy argument applies to $h\in L^2(t_1,t_2;H^{-1})$).
Letting $V=w_\ep$ and $P=\pi_\ep$, this implies that $w_\ep$ and $\nb \pi_\ep$ are unique.
Up to a function $\pi_*(t)$ independent of $x$,
\EQ{\label{eq-BT2019-3.11}
\pi_\ep(x,t) - \pi_*(t)
  = -\frac13 \tr  H(x,t)
 + \lim_{\de\to0} \int_{| x-y|>\de} K_{ij}(x-y) H_{ij} (y,t)\, dy,
}
where 
$K_{ij}(x) = \pd_i\pd_j\frac1{4\pi|x|}$.
The right-hand side of \eqref{eq-BT2019-3.11} is defined in $L^2(t_1,t_2; L^2(\R^3))$.
Since the only appearance of $\pi_\ep$ in \eqref{eq-BT2019-3.5} is $\nb\pi_\ep$, we can redefine $\pi_\ep$ to equal to $\pi_\ep - \pi_*(t)$ and, therefore, drop $\pi_*(t)$ from \eqref{eq-BT2019-3.11}.

The pressure $\pi_\ep$ given by \eqref{eq-BT2019-3.11} is already bounded in $L^2(t_1,t_2; L^2(\R^3))$ for any $0<t_1<t_2<\infty$ but the bound depends on $t_1, t_2$ and $\ep$.
We now bound it in $L^{3/2}(0,T; L^{3/2}(B_\la))$.
Bounding the first term from \eqref{eq-BT2019-3.11} is simple given H\"older's inequality, \eqref{eq-BT2019-3.8}, and \eqref{eq-BT2019-3.9}.
In particular, similar to \eqref{eq-BT2019-3.10}, we have for any $\ga>0$
\[
\int_0^t \norm{\tfrac13 \tr H(\cdot,s)}_{L^{3/2}(B_\la)}^{3/2} ds
\le  \text{RHS of } \eqref{eq-BT2019-3.10}.
\]
To bound the  principal  value integral in \eqref{eq-BT2019-3.11}, we split the integral into the sum of local and nonlocal parts, i.e., sum of
\EQN{
\pi_{\rm near}(x,t)&= \lim_{\de\to0} \int_{|y|<\la^2,\,|x-y|>\de} K_{ij}(x-y) H_{ij}(y,t)\, dy\\
\pi_{\rm far}(x,t)&= \int_{|y|>\la^2} K_{ij}(x-y) H_{ij}(y,t)\, dy.
}
To bound $\pi_{\rm near}$ note that, by the Calder\'on--Zygmund theory,
\[
\norm{\pi_{\rm near}(\cdot,t)}_{L^{3/2}(B_\la)}
\le C \norm{(\eta_{\ep\sqrt{t}} * v_\ep)\otimes v_\ep(\cdot,t)}_{L^{3/2}(B_{\la^2})} + C\norm{(\eta_{\ep\sqrt{t}} * b_\ep)\otimes b_\ep(\cdot,t)}_{L^{3/2}(B_{\la^2})},
\]
and, arguing as above using \eqref{eq-BT2019-3.8} but with $B_{\la^2}$ in place of $B_\la$, it follows that 
\[
\int_0^t \norm{\pi_{\rm near}(\cdot,s)}_{L^{3/2}(B_\la)}^{3/2}\, ds\le  \text{RHS of } \eqref{eq-BT2019-3.10}.
\]

Now we bound the term $\pi_{\rm far}$. For $x\in B_\la$, we have $|x-y|\sim |y|$ for $|y|\ge\la^2$, and hence
$|\pi_{\rm far}(x,s)|\le C\int_{|y|\ge\la^2} |y|^{-3} |H(y,s)|\, dy$.
 By Lemma \ref{lem2.2}, $H$ is $\la$-DSS with degree 2. Let $A_k=\{x:\la^{k-1}\le|x|<\la^k\}$. Changing variables $y=\la^k z$ on $A_k$, we get for $x\in B_\la$ and $0<s<1$,
\EQN{
|\pi_{\rm far}(x,s)| &\le C\sum_{k\ge3}\int_{A_k}|y|^{-3}|H(y,s)|\,dy
= C\sum_{k\ge3}\la^{-2k}\int_{A_0}|z|^{-3}|H(z,\la^{-2k}s)|\,dz
\\
&\le C(\la)\sum_{k\ge3}\la^{-2k}\int_{B_1}\bke{|\eta_{\ep\sqrt{\tau}}*v_\ep|\,|v_\ep|+|\eta_{\ep\sqrt{\tau}}*b_\ep|\,|b_\ep|}(z,\tau)\Big|_{\tau=\la^{-2k}s}\,dz
\\
&\le C(\la)\sum_{k\ge3}\la^{-2k}\bke{\norm{v_\ep(\la^{-2k}s)}_{L^2(B_\la)}^2+\norm{b_\ep(\la^{-2k}s)}_{L^2(B_\la)}^2}
\le C(\la)\,\td\al_\ep(s),
}
where we used $\supp\eta_{\ep\sqrt\tau}\subset B_{\la-1}$, Young's inequality, the rescaling
$\norm{v_\ep(\tau)}_{L^2(B_\la)}^2=\la\norm{v_\ep(\la^{-2}\tau)}_{L^2(B_1)}^2$ (and the same for $b_\ep$), and that $\td\al_\ep$ is nondecreasing.

Therefore,
\[
\int_0^t \norm{\pi_{\rm far}(\cdot,s)}_{L^{3/2}(B_\la)}^{3/2}\, ds
\le C(\la) \int_0^t \td\al_\ep(s)^{3/2}\, ds.
\]
After using H\"older's inequality, \eqref{eq-BT2019-3.9}, the bounds above, and $\al^{3/2}\le\al+\al^3$ for $\al>0$, it is clear that 
\EQ{\label{eq4.14}
\int_0^t \norm{\pi_\ep(\cdot,s)}_{L^{3/2}(B_\la)}^{3/2}\, ds + 2\int_0^t \int |\pi_\ep(v_\ep\cdot\nb\phi)|\, dxds
\le \text{RHS of } \eqref{eq-BT2019-3.10}.
}

 The pressure $\pi^b_\ep$ is treated in the same way. The function $c_\ep=b_\ep-b_L$ satisfies
$\pd_tc_\ep - \De c_\ep + \nb\pi^b_\ep = h^b$, $\div c_\ep=0$, where $h^b_i = -\pd_jH^b_{ji}$ with
\[
H^b = (\eta_{\ep\sqrt{t}} * v_\ep)\otimes b_\ep - (\eta_{\ep\sqrt{t}} * b_\ep)\otimes v_\ep .
\]
By the same uniqueness argument as for \eqref{eq-BT2019-3.11}, after subtracting a function of $t$ only,
\EQ{\label{eq-pib-formula}
\pi^b_\ep(x,t) = -\frac13 \tr H^b(x,t) + \lim_{\de\to0} \int_{|x-y|>\de} K_{ij}(x-y) H^b_{ij}(y,t)\, dy .
}
Since $H^b$ is bilinear in $(v_\ep,b_\ep)$ and their mollifications, exactly as $H$, the argument above also gives
\[
\int_0^t \norm{\pi^b_\ep(\cdot,s)}_{L^{3/2}(B_\la)}^{3/2}\, ds + 2\int_0^t \int |\pi^b_\ep(b_\ep\cdot\nb\phi)|\, dxds
\le \text{RHS of } \eqref{eq-BT2019-3.10}.
\]

For the  second-to-last term on the right-hand side of \eqref{eq-MHD-LEE-epsilon}, we use Young's inequality and follow the same argument deriving \eqref{eq-BT2019-3.10} to get
\[
- 2\int_0^t\int (v_\ep\cdot b_\ep)\left((\eta_{\ep\sqrt{s}} * b_\ep) \cdot\nb\phi\right) dxds
\le \text{RHS of } \eqref{eq-BT2019-3.10}.
\]

Combining the estimates above for the right-hand side of \eqref{eq-MHD-LEE-epsilon} and taking $\ga$ sufficiently small to absorb the gradient terms on the right-hand side, we obtain
\EQ{\label{eq-BT2019-3.12}
\al_\ep(t) + \int_0^t \int_{B_1} \bke{|\nb v_\ep|^2 + |\nb b_\ep|^2} dxds
\le \al_0 + C(\la,\eta,\ga) \int_0^t \left(\td\al_\ep(s)^3 + \td\al_\ep(s)\right) ds.
}
Therefore, 
\[
\td\al_\ep(t) \le \al_0 + C(\la,\eta,\ga) \int_0^t \left(\td\al_\ep(s)^3 + \td\al_\ep(s)\right) ds.
\]
By continuity of $\al_\ep(t)$, we have
\[
\td\al_\ep(t) \le 2\al_0\ \text{ for all $t<T$,}
\]
for some $T>0$. 
By a continuity argument, we may take $T=\min\{1,(C(2+8\al_0^2))^{-1}\}$; the restriction $T\le1$ is needed because the estimates above are for $0<t\le1$.

Using \eqref{eq4.7}$_4$ and letting $\ep\to0$ yields  (up to the subsequence $\ep_j$)
\[
(v,\chi_{B_1}v)(t) + (b,\chi_{B_1}b)(t) \le \liminf_{\ep\to0} (v_{\ep},\chi_{B_1}v_{\ep})_{L^2}(t) + \liminf_{{\ep}\to0} (b_{\ep},\chi_{B_1}b_{\ep})_{L^2}(t) \le 2\al_0
\]
for all $t\le T$.
Note that \eqref{eq-BT2019-3.12} gives uniform (in $\ep$) control of
\[
\int_0^T\int_{B_1}\bke{|\nb v_\ep|^2 + |\nb b_\ep|^2} dxdt \le C(\al_0,\la)
\]
for some constant $C(\al_0,\la)$.
From \cite{Lai-JMFM2019} we have that $v_\ep, b_\ep$ converges weakly to $v, b$, respectively, in $L^2(1/k, T; H^1(B_1))$ for every $k\in\NN$.
Hence,
\[
\int_{1/k}^T\int_{B_1} \bke{|\nb v|^2 + |\nb b|^2} dxdt
\le \sup_{\ep>0} \int_0^T\int_{B_1} \bke{|\nb v_\ep|^2 + |\nb b_\ep|^2} dxdt,
\]
and, letting $k\to\infty$, it follows that
\[
\int_0^T\int_{B_1} \bke{|\nb v|^2 + |\nb b|^2} dxdt \le C(\al_0,\la).
\]
Similarly, since $\pi_\ep\in L^{3/2}(0,T;L^{3/2}(B_\la))$ with uniformly bounded norms, it follows that 
\[
\pi_{\ep_j}\wkto \pi \quad\text{ weakly in }\quad L^{3/2}(0,T; L^{3/2}(B_\la)).
\]
 Here $\pi_\ep$ is normalized as in \eqref{eq-BT2019-3.11}; any other pressure associated with $(v,b)$ differs from $\pi$ by a function of $t$ only, which does not affect \eqref{eq-MHD-LEI}.
The hypotheses of Lemma \ref{BT2019-lem-2.1} on $\td v_k=\eta_{\ep_j\sqrt t}*v_{\ep_j}$ and $\td b_k=\eta_{\ep_j\sqrt t}*b_{\ep_j}$ follow from \eqref{eq4.7} and the bounds above, since $\norm{\eta_{\ep\sqrt t}*(v_\ep-v)}_{L^2(B_1)}\le\norm{v_\ep-v}_{L^2(B_{1+\ep})}$ for $t\le1$, $\eta_{\ep\sqrt t}*v\to v$ in $L^2(0,T;L^2(B_1))$, and the bounds on $B_\la$ follow by DSS scaling.
Applying Lemma \ref{BT2019-lem-2.1} yields the desired pressure representation \eqref{eq-BT2019-3.3} in $L^{3/2}(0,T; L^{3/2}(B_\la))$.

 It remains to show that the magnetic pressure disappears in the limit. By the bound on $\pi^b_\ep$ above, after passing to a further subsequence, $\pi^b_{\ep_j}\wkto \pi^b$ weakly in $L^{3/2}(0,T; L^{3/2}(B_\la))$ for some $\pi^b$. Let
\[
G = v\otimes b - b\otimes v .
\]
The proof of Lemma \ref{BT2019-lem-2.1}, with $h_{ij}$ replaced by $H^b_{ij}-G_{ij}$ and using \eqref{eq-pib-formula}, shows that
\[
\pi^b(x,t) = -\frac13 \tr G(x,t) + \lim_{\de\to0} \int_{|x-y|>\de} K_{ij}(x-y) G_{ij}(y,t)\, dy
\quad\text{in } L^{3/2}((0,T)\times B_\la).
\]
Since $G$ is antisymmetric and $K_{ij}=K_{ji}$, we have $\tr G=0$ and $K_{ij}G_{ij}=0$, so $\pi^b=0$. By DSS scaling the same holds on $(0,T')\times E$ for every $T'>0$ and compact $E$. Consequently, the limit $(v,b)$ satisfies the second equation of \eqref{eq-MHD} in the sense of distributions, tested with all smooth compactly supported vector fields (not only divergence-free ones). Moreover, $\pi^b_\ep\wkto0$ weakly in $L^{3/2}_{\loc}$ and $b_\ep\to b$ in $L^3_{\loc}$, so the last term of \eqref{eq-MHD-LEE-epsilon} vanishes in the limit. Letting $\ep\to0$ in \eqref{eq-MHD-LEE-epsilon}, we conclude that $(v,b)$ satisfies the local energy inequality \eqref{eq-MHD-LEI} with initial data $(v_0,b_0)$: the initial term does not depend on $\ep$, the other terms on the right converge by \eqref{eq4.7} and the strong convergence of $v_\ep$, $b_\ep$ and their mollifications in $L^3_\loc$, and the left side is lower semicontinuous for every $t$ by \eqref{eq4.7}$_4$. See also Remark \ref{rem-pib} for a different argument that the magnetic pressure vanishes. The same remark applies to the equation for ${\bf f}$ in the veNSD system \eqref{eq-veNSD}.
This concludes the proof of Proposition \ref{BT2019-prop-3.1}.
\end{proof}

\begin{remark}\label{approx-scheme}
This remark concerns the regularized system \eqref{eq4.4}.
 We note that \cite{Lai-JMFM2019} uses the similarity variables $y=x/\sqrt{2t}$ and $s=\log\sqrt{2t}$, in which $L=\pd_s-\De-1-y\cdot\nb$ and the period is $\log\la$; the equations of \cite{Lai-JMFM2019} quoted below become ours after this change of variables, and the mollifier $\eta_\ep$ there corresponds to $\eta_{\ep\sqrt{2t}}$ in physical variables.
In a nonlinear term $W\cdot \nabla V$ we call $W$ the drift. 
In \cite{Lai-JMFM2019}, we first construct a revised asymptotic profile $(W,D)$ in \cite[(2.9)--(2.10)]{Lai-JMFM2019},
and then decompose $V=W+ U$, $B=D+ A$. We consider the equations for the perturbation $(U,A)$ and their approximations $(U_\ep,A_\ep)$. In the equation \cite[(2.16)]{Lai-JMFM2019} for the approximation $(U_\ep,A_\ep)$, 
 only the parts $U_\ep$ and $A_\ep$ of the drifts of $\nb U_\ep$ and $\nb A_\ep$ are mollified, as in $(W+\eta_\ep*U_\ep)\cdot \nabla A_\ep$; the profile $(W,D)$ is not mollified, and neither are the drifts in the terms $U_\ep\cdot\nb W$, $A_\ep\cdot\nb D$, $U_\ep\cdot\nb D$ and $A_\ep\cdot\nb W$. As a result, $V_\ep=W+ U_\ep$, $B_\ep=D+ A_\ep$ do not satisfy \eqref{eq4.4}. This can be fixed by mollifying all drifts in  \cite[(2.16)]{Lai-JMFM2019} including  $(W,D)$. Then \eqref{eq4.4} is valid and the same proof goes through.
This form of \eqref{eq4.4} simplifies the local energy equality \eqref{eq-MHD-LEE-epsilon}.
This remark also applies to \cite[(3-4)]{BT-AnalPDE2019}
 when we apply \cite[(2.23)]{BT1-AHP}, and to our construction in Proposition \ref{BT2019-prop-3.1-OB}
  when we apply \cite[(3.2)--(3.4)]{Tsai-OB} to \eqref{eq4.18}.
\end{remark}

\begin{remark}\label{rem-pib}
 A pressure term $\nb P^B_\ep$ was missing from the magnetic equation of the mollified perturbed Leray system for the MHD equations \cite[(2.16)]{Lai-JMFM2019}. This term is necessary because the quadratic term is no longer divergence-free after mollification, so the equation for $A_\ep$, which is tested only against divergence-free fields, holds only modulo a gradient.

Nevertheless, the proofs in \cite{Lai-JMFM2019} remain valid with this correction because:
\EN{
\item The additional pressure term does not affect the Galerkin argument, as it vanishes when tested against divergence-free vector fields.

\item The boundedness and convergence of $P^B_\ep$ can be established using the same arguments as those for the fluid pressure $P_\ep$.

\item Since the limiting pressure satisfies $P^B \in L^{5/3}(\R^3\times(0, 2\log\la))$ and $\Delta P^B = 0$ in $y$ for a.e.\ $s$, it follows that $P^B = 0$ by the Liouville theorem for $L^p$ harmonic functions.
}
A different argument for the vanishing of the magnetic pressure, based on the antisymmetry of $v\otimes b-b\otimes v$, is given at the end of the proof of Proposition \ref{BT2019-prop-3.1}.

The same remark applies to the mollified perturbed Leray system for the viscoelastic Navier--Stokes equations with damping \cite[(2.55)]{Lai-JMFM2019}.
\end{remark}

The next lemma is used to pass to the limit in the local energy inequalities for OB. The buoyancy term is critical with respect to Hardy's inequality, and we can only control a possible defect by the dissipation of both $v$ and $\th$.

\begin{lem}[Buoyancy term]\label{lem-buoyancy}
Let $t>0$, $R>0$, and let $\phi\in C^\infty_c(\R^3\times[0,\infty))$ be nonnegative with $\supp\phi\subset B_R\times[0,\infty)$. Suppose that, as $k\to\infty$, $v_k\to v$ and $\th_k\to\th$ in $L^2(B_R\times(0,t))$, $\nb v_k\wkto\nb v$ and $\nb\th_k\wkto\nb\th$ weakly in $L^2(B_R\times(0,t))$,  and $\div v_k=\div v=0$. Then, for any $\mu\ge1$,
\EQN{
\limsup_{k\to\infty}\int_0^t\!\!\int&\bke{2\th_k\nb G\cdot v_k - 2|\nb v_k|^2-2\mu|\nb\th_k|^2}\phi\,dxds
\\
&\le \int_0^t\!\!\int\bke{2\th\nb G\cdot v - 2|\nb v|^2-2\mu|\nb\th|^2}\phi\,dxds.
}
\end{lem}

\begin{proof}
Let $g_k=\th_k-\th$, $h_k=v_k-v$, and $\psi=\phi^{1/2}$. By Glaeser's inequality $|\nb_x\phi|^2\le 2M\phi$, where $M=\norm{\nb_x^2\phi}_{L^\infty}>0$ (it follows from $0\le\phi(x-\tau\nb_x\phi(x,s),s)\le\phi(x,s)-\tau|\nb_x\phi(x,s)|^2+\frac M2\tau^2|\nb_x\phi(x,s)|^2$ with $\tau=1/M$, by Taylor's theorem), the function $\psi$ is Lipschitz in $x$  with constant $\sqrt{M/2}$, uniformly in $s$. All norms below are in $L^2(\R^3\times(0,t))$.

Since $\nb h_k\wkto0$ and $\nb g_k\wkto0$ weakly,
\EQN{
\int_0^t\!\!\int|\nb v_k|^2\phi &= \int_0^t\!\!\int|\nb v|^2\phi + \norm{\psi\nb h_k}^2 + o(1),\\
\int_0^t\!\!\int|\nb \th_k|^2\phi &= \int_0^t\!\!\int|\nb \th|^2\phi + \norm{\psi\nb g_k}^2 + o(1).
}
Next, $\th_k\nb G\cdot v_k-\th\nb G\cdot v=\nb G\cdot(g_kv+\th h_k+g_kh_k)$ and $|\nb G|=|x|^{-2}$.
By Hardy's inequality $\norm{u/|x|}_{L^2(\R^3)}\le 2\norm{\nb u}_{L^2(\R^3)}$, the functions $\psi g_k/|x|$ and $\psi h_k/|x|$ are bounded in $L^2(\R^3\times(0,t))$, and they converge to $0$ in $L^2((\R^3\setminus B_\de)\times(0,t))$ for every $\de>0$; hence they converge weakly to $0$ in $L^2(\R^3\times(0,t))$. Since $\psi\th/|x|$ and $\psi v/|x|$ belong to $L^2(\R^3\times(0,t))$, the terms with $g_kv$ and $\th h_k$ tend to $0$.  For the last term, since $\div h_k=0$ and $G=|x|^{-1}$, an integration by parts gives $\int\nb G\cdot g_kh_k\,\phi\,dx=-\int G\,h_k\cdot\nb(g_k\phi)\,dx$, where all integrals converge absolutely by Hardy's inequality. Hence
\[
\abs{2\int_0^t\!\!\int\nb G\cdot g_kh_k\,\phi}
\le 2\int_0^t\!\!\int \frac{|h_k||\nb g_k|\phi}{|x|} + 2\int_0^t\!\!\int \frac{|h_k||g_k||\nb\phi|}{|x|}
=: I_1 + I_2.
\]
Since $h_k\nb\psi$ tends to $0$ in $L^2$,
\[
I_1 \le 2\norm{\frac{\psi h_k}{|x|}}\norm{\psi\nb g_k} \le 4\norm{\nb(\psi h_k)}\norm{\psi\nb g_k} = 4\norm{\psi\nb h_k}\norm{\psi\nb g_k} + o(1),
\]
and, since $h_k\to0$ in $L^2(B_R\times(0,t))$ and $g_k\nb\phi$ is bounded in $L^2(0,t;H^1)$,
\[
I_2 \le 2\norm{h_k}_{L^2(B_R\times(0,t))}\norm{\frac{g_k\nb\phi}{|x|}} \le 4\norm{h_k}_{L^2(B_R\times(0,t))} \norm{\nb(g_k\nb\phi)} = o(1).
\]
Therefore, with $a_k=\norm{\psi\nb h_k}$ and $b_k=\norm{\psi\nb g_k}$, the left side of the lemma is bounded by its right side plus
$\limsup_{k}\bke{4a_kb_k-2a_k^2-2\mu b_k^2}\le0$, since $2a^2+2\mu b^2\ge4\sqrt{\mu}\,ab\ge4ab$ for $\mu\ge1$.
\end{proof}

\begin{prop}[OB]\label{BT2019-prop-3.1-OB}
Fix $\la>1$. Assume $v_0, \th_0\in L^3_w(\R^3)$  are  a $\la$-DSS divergence-free vector field and a $\la$-DSS function, respectively,  and $f=f_1+f_2$, where  $f_1,f_2$ are 
$(\la,3)$-DSS vector fields, $f_j(x,t) = \la^3f_j(\la x,\la^2t)$, $j=1,2$,
\[
 f_1 \in L^\infty(1,\la^2; L^{6/5,2}(\R^3)),  \quad f_2\in L^\infty(1,\la^2; L^2(\R^3)),\quad
\div f_2=0,
\]
and $(v,\th)$ is a $\la$-DSS solution of the Oberbeck-Boussinesq system \eqref{eq-OB} evolving from $(v_0, \th_0)$ with force $f$, constructed as in \cite{Tsai-OB} with the modification of Remark \ref{approx-scheme}, and $\pi$ is its associated pressure. 
Let $\al_0=\norm{v_0}_{L^2(B_\la)}^2 + \norm{\th_0}_{L^2(B_\la)}^2$.
Then, there exist positive $T=T(\al_0,\la, f)$ and $C=C(\al_0,\la,f)$ independent of $\norm{v_0}_{L^2_\uloc}$, $\norm{\th_0}_{L^2_\uloc}$, $\norm{v_0}_{L^3_w}$, and $\norm{\th_0}_{L^3_w}$ so that 
\EQ{
\esssup_{0\le t\le T} \int_{B_1} \left(|v(x,t)|^2 + |\th(x,t)|^2\right) dx + \int_0^T\int_{B_1} \left(|\nb v|^2 + |\nb \th|^2\right) dxdt < C,
}
and 
\EQ{\label{eq-BT2019-3.2-OB}
\int_0^T\int_{B_\la} |\pi(x,t)|^{3/2}\, dxdt < C.
}
Moreover,  the pressure  satisfies the formula
\EQS{\label{eq-BT2019-3.3-OB}
\pi(x,t) &= -\frac13|v|^2(x,t) + \lim_{\de\to0} \int_{| x-y|>\de} K_{ij}(x-y)v_i(y,t)v_j(y,t)\, dy \\
&\quad - \int \nb\bke{\frac1{4\pi|x-y|}}\cdot(\th\nb G + f_1)(y,t)\, dy
}
in $L^{3/2}(B_\la\times(0,T))$.
 Furthermore, $(v,\th)$ satisfies the local energy inequalities \eqref{eq-OB-LEI} with initial data $(v_0,\th_0)$. 
\end{prop}

In the  integral formula \eqref{eq-BT2019-3.3-OB} for the pressure, $f_2$ is absent because $\div f_2=0$.

\begin{proof}
The proof of the proposition is an adaption of the proof of \cite[Proposition 3.1]{BT-AnalPDE2019} for the Navier-Stokes equations to the Oberbeck-Boussinesq system.

 By the reduction in the proof of Lemma \ref{BT2019-lem-2.1-OB}, we may assume $\div f_1=0$. Indeed, the force enters the approximation scheme of \cite{Tsai-OB} only through its pairing with divergence-free fields, so $(v,\th)$ is also the solution constructed for the force $\mathbb Pf_1+f_2$, with pressure $\pi-q_1$. The bounds for $(v,\th)$, the formula \eqref{eq-BT2019-3.3-OB} and the local energy inequalities are not affected, \eqref{eq-BT2019-3.2-OB} for $\pi$ follows from that for $\pi-q_1$ and $q_1\in L^{3/2}(0,T;L^{3/2}(B_\la))$, and $M_1$ is changed only to $CM_1$.

Since $(v,\th)$ is a solution of \eqref{eq-OB} from \cite{Tsai-OB}, its image in the similarity variables \eqref{eq-self-similar-variables} solves the time-periodic Leray equations \eqref{eq-OB-Leray} and is the limit of a mollified approximation scheme \cite[(3.2)--(3.4)]{Tsai-OB}.
In particular, for each $\ep>0$, there exists a time-periodic solution $(V_\ep, \Th_\ep, P_\ep)$ to the problem (regularized version of \eqref{eq-OB-Leray}) in $\R^3 \times \R$
\begin{align}
\nonumber
\left(\pd_sV_\ep - \De V_\ep - \frac12 V_\ep - \frac12 y\cdot\na V_\ep + (\eta_\ep * V_\ep)\cdot\na V_\ep + \nb P_\ep\right)(y,s) &= (\Th_\ep\nb G + F)(y,s),\\
\label{eq4.18}
\left(\pd_s\Th_\ep - \De \Th_\ep - \frac12 \Th_\ep - \frac12 y\cdot\na \Th_\ep + (\eta_\ep * V_\ep)\cdot\na \Th_\ep\right)(y,s) &= 0,
\end{align}
where $\eta_\ep(y) = (1/\ep^3)\eta(y/\ep)$ and $\eta$ is a fixed mollifier satisfying \eqref{mollifier}. Moreover, a subsequence of $(V_\ep,\Th_\ep,P_\ep)$ converges weakly in the energy class locally to $(V,\Th,P)$ as $\ep \to 0$.

A remark is in order: similar to \eqref{eq4.4}, 
the approximation scheme in \cite[(3.2)--(3.4)]{Tsai-OB} needs to be slightly changed to get \eqref{eq4.18}, so that we have simpler local energy equalities \eqref{eq-OB-LEE-epsilon} and \eqref{eq-OB-LEE-epsilon2} below. See Remark \ref{approx-scheme}.

Applying the similarity transform \eqref{eq-OB-SS-transform} and \eqref{eq-self-similar-variables}, we obtain a $\la$-DSS solution $(v_\ep, \th_\ep,\pi_\ep)$ to the system
\EQS{\label{eq-BT2019-3.5-OB}
\pd_tv_\ep(x,t) - \De v_\ep(x,t) + (\eta_{\ep\sqrt{t}} * v_\ep)\cdot\nb v_\ep(x,t) + \nb\pi_\ep(x,t) &= (\th_\ep\nb G + f)(x,t),\\[1mm]
\pd_t\th_\ep(x,t) - \De \th_\ep(x,t) + (\eta_{\ep\sqrt{t}} * v_\ep)\cdot\nb \th_\ep(x,t) &= 0,
}
and a DSS solution $(v, \th,\pi)$ to the OB system \eqref{eq-OB}. 

By the convergence properties of $(V_{\ep_j},\Th_{\ep_j},P_{\ep_j})$ to $(V,\Th,P)$
\cite[display before (3.28)]{Tsai-OB}, where the convergence of the pressure holds on compact sets \cite[Section 1]{Tsai-OB-corr}, and discretely self-similar scaling (to extend the estimates down to $t=0$; for the pressure we use the revised \cite[(3.27)]{Tsai-OB} in \cite[Section 1]{Tsai-OB-corr}, whose force term vanishes since $\div f_1=\div f_2=0$), it follows that for all $T>0$ and all compact sets $E$ of $\R^3$, (up to the subsequence $\ep_j$)
\EQS{
v_\ep\to v,\ \th_\ep\to \th\ \text{ weakly }&\text{in $L^2(0,T;H^1(E))$},\\
v_\ep\to v,\ \th_\ep\to \th\ \text{ strongly }&\text{in $L^2(0,T;L^2(E))$},\\
\pi_\ep\wkto \pi\ \text{ weakly }&\text{in $L^{5/3}(0,T;L^{5/3}(E))$},\\
v_\ep(s)\to v(s),\ \th_\ep(s)\to \th(s)\ \text{ weakly }&\text{in $L^2(E)$ for all $s\in[0,T]$}.
}

 We first show the time continuity \eqref{wep-cont} of the approximations, which is needed below.  Let $w_\ep = v_\ep - v_L$ and $\tilde \theta_\ep=\th_\ep-\th_L$, where $v_L(\cdot,t) = e^{t\De}v_0$ and $\th_L(\cdot,t)=e^{t\De}\th_0$.
For $0<t_1<t_2<\infty$, we have  (the properties of $w_\ep$ and $\tilde \theta_\ep$ follow from the construction in \cite{Tsai-OB}, cf.\ \cite[(1.6)]{Tsai-OB})
\EQN{
v_L,\th_L&\in C([t_1,t_2]; L^4(\R^3)\cap L^\infty(\R^3)),\\
w_\ep,\tilde \theta_\ep&\in L^\infty(t_1,t_2;L^2(\R^3))\cap L^2(t_1,t_2;H^1(\R^3))\subset L^4(t_1,t_2;L^3(\R^3)),\\
f &\in L^\infty(t_1,t_2; (L^{6/5,2}+L^2)(\R^3)) \subset L^2(t_1,t_2; H^{-1}(\R^3)).
}
We have
\EQ{\label{eq-wep}
\pd_tw_\ep - \De w_\ep + \nb\pi_\ep = h,\qquad\div w_\ep = 0,
}
where $h_i = -\pd_jH_{ji} + \th_\ep\pd_iG + (f)_i$ 
with
\EQN{
H &= (\eta_{\ep\sqrt{t}} * v_\ep)\otimes v_\ep\\
&= (\eta_{\ep\sqrt{t}} * w_\ep + \eta_{\ep\sqrt{t}} * v_L)\otimes (w_\ep + v_L).
}
By Young's convolution inequality,
\[
\norm{H}_{L^2(t_1,t_2;L^2)} \lec \norm{\eta_{\ep\sqrt{t}}}_{L^\infty(t_1,t_2;L^{\frac65}\cap L^1)} \left(\norm{w_\ep}_{L^4(t_1,t_2;L^3(\R^3))}^2  + \norm{v_L}_{L^4(\R^3\times[t_1,t_2])}^2\right).
\]
 We write $\th_\ep\pd_iG=\tilde \theta_\ep\pd_iG+\th_L\pd_iG$. By Hardy's inequality, for $\psi\in \dot H^1$ with $\norm{\nb\psi}_{L^2}=1$,
\[
\abs{\int \tilde \theta_\ep\pd_iG\psi\, dx}\le \norm{\frac{\tilde \theta_\ep}{|\cdot|}}_{L^2} \norm{\frac{\psi}{|\cdot|}}_{L^2}
\le 4 \norm{\nb\tilde \theta_\ep}_{L^2} \norm{\nb\psi}_{L^2}
= 4 \norm{\nb\tilde \theta_\ep}_{L^2}.
\]
Since $|\pd_iG|\le|x|^{-2}$, by H\"older's inequality,
\[
\norm{\th_L\pd_iG}_{L^{6/5}} \le \norm{|\cdot|^{-2}}_{L^{6/5}(B_1)}\norm{\th_L}_{L^\infty} + \norm{|\cdot|^{-2}}_{L^{12/7}(\R^3\setminus B_1)}\norm{\th_L}_{L^4},
\]
which is finite, and $L^{6/5}(\R^3)\subset \dot H^{-1}(\R^3)$ by Sobolev embedding. Hence $\th_\ep\pd_iG\in L^2(t_1,t_2; \dot H^{-1})\subset L^2(t_1,t_2; H^{-1})$.
Since $h\in L^2(t_1,t_2;H^{-1})$, and similarly $\pd_t\tilde \theta_\ep-\De\tilde \theta_\ep=-\div((\eta_{\ep\sqrt{t}} * v_\ep)\th_\ep)\in L^2(t_1,t_2;H^{-1})$  by Young's convolution inequality as for $H$, \cite[Lemma A.2]{CKN-CPAM1982} implies
\EQ{\label{wep-cont}
 w_\ep,\ \tilde \theta_\ep\in C([t_1,t_2];L^2(\R^3))\quad\text{for all }0<t_1<t_2<\infty.
}

Note also that $v_\ep(t)\to v_0$ and $\th_\ep(t)\to \th_0$ as $t\to0^+$ in $L^2_\loc$; i.e., the mollification does not affect the initial data.
Together with \eqref{wep-cont} for $w_\ep=v_\ep-v_L$ and $\tilde \theta_\ep=\th_\ep-\th_L$, it follows that, for a constant $C_1\ge1$ to be determined (see \eqref{C1.def}),
\[
\al_\ep(t) = \int_{B_1} \left(|v_\ep(x,t)|^2 + C_1\th_\ep(x,t)^2\right) dx
\]
and 
\[
\td\al_\ep(t) = \sup_{0\le\tau\le t} \al_\ep(\tau)
\]
are continuous in $t$. 
Our goal is to establish local-in-time a priori bounds for $\al_\ep(t)$ that are independent of $\ep$. 
Note that $(v_\ep, \th_\ep)$ satisfies the local energy equalities (the approximations need not be smooth, because of the force $f$ and the singular term $\th_\ep\nb G$; the equalities follow by testing the equations with $2v_\ep\phi$ and $2\th_\ep\phi$, which is justified by the Lions--Magenes lemma since $\pd_tv_\ep,\pd_t\th_\ep\in L^2_\loc((0,\infty);H^{-1}_\loc)$, all terms being absolutely integrable, by Hardy's inequality for the buoyancy term, and by the right continuity at $t=0$); i.e.,
\begin{align}
\nonumber
\int& |v_\ep|^2\phi(x,t)\, dx + 2\int_0^t \int |\nb v_\ep|^2\phi\, dxds\\
\nonumber
&= \int |v_0|^2\phi(x,0)\, dx + \int_0^t \int |v_\ep|^2\left(\pd_s\phi+\De\phi\right) dxds
 + \int_0^t \int |v_\ep|^2 \bke{(\eta_{\ep\sqrt{s}} * v_\ep)\cdot\nb\phi} dxds\\
\label{eq-OB-LEE-epsilon}
&\quad + 2\int_0^t \int \pi_\ep (v_\ep\cdot\nb\phi) dxds + 2\int_0^t \int (\th_\ep\nb G + f)\cdot v_\ep\phi\, dxds,
\end{align}
and
\EQS{\label{eq-OB-LEE-epsilon2}
\int& \th_\ep^2\phi(x,t)\, dx + 2\int_0^t \int |\nb\th_\ep|^2\phi\, dxds\\
&= \int \th_0^2\phi(x,0)\, dx + \int_0^t \int \th_\ep^2(\pd_s\phi + \De\phi)\, dxds
 + \int_0^t \int \th_\ep^2 \bke{(\eta_{\ep\sqrt{s}} * v_\ep)\cdot\nb\phi} dxds
}
for any nonnegative $\phi\in C^\infty_c(\R^3\times[0,\infty))$.
Fix $\chi\in C^\infty(\R)$ with $\chi(t) = 1$ if $t\le 1$, $0\le\chi(t)\le1$ for $t\in[1,\la]$, and $\chi(t) = 0$ if $t\ge\la$.
We now choose $\phi$ in \eqref{eq-OB-LEE-epsilon}--\eqref{eq-OB-LEE-epsilon2} as 
\[
\phi(x,t) = \chi^2(|x|)\cdot\chi(t).
\]

We will multiply the second equation \eqref{eq-OB-LEE-epsilon2} by $C_1$, add it to the first equation \eqref{eq-OB-LEE-epsilon}, and estimate the terms on the right-hand side to get  an a priori  bound. The  first three terms of the right-hand side of \eqref{eq-OB-LEE-epsilon} and the full right-hand side of \eqref{eq-OB-LEE-epsilon2} can be controlled as in the proof of Proposition \ref{BT2019-prop-3.1} (see also the proof of \cite[Proposition 3.1]{BT-AnalPDE2019}).
More precisely, we have
\[
\int (|v_0|^2 + C_1\th_0^2)\phi(x,0)\, dx\le \al_1:=(1+C_1)\al_0,
\]
\EQN{
\int_0^t \int (|v_\ep|^2 + C_1\th_\ep^2)\left(\pd_s\phi+\De\phi\right) dxds\le C(\la) \int_0^t \td\al_\ep(s)\, ds,
}
and by the same derivation of \eqref{eq-BT2019-3.10}, for any $\gamma>0$,
\EQS{\label{eq-BT2019-3.10-OB}
\int_0^t&\int (|v_\ep|^2 + C_1\th_\ep^2) \left((\eta_{\ep\sqrt{s}} * v_\ep)\cdot\nb\phi\right) dxds\\
&\quad \le C(\ga) \int_0^t \left(\td\al_\ep(s)^3 + \td\al_\ep(s)\right) ds 
+ C\ga\int_0^t\int \left(|\nb v_\ep|^2 + |\nb \th_\ep|^2\right)\phi\, dxds.
}

It remains to estimate the pressure and forcing terms in the last line of \eqref{eq-OB-LEE-epsilon}.
We first estimate the pressure term.
By \eqref{eq-wep}, \eqref{wep-cont}  and the same uniqueness argument  before \eqref{eq-BT2019-3.11},
up to a function $\pi_*(t)$ independent of $x$,
\EQS{\label{eq-BT2019-3.11-OB}
\pi_\ep(x,t) - \pi_*(t)
 &= -\frac13 \left((\eta_{\ep\sqrt{t}} * v_\ep)\cdot v_\ep\right)(x,t)\\
&\quad + \lim_{\de\to0} \int_{| x-y|>\de} K_{ij}(x-y) (\eta_{\ep\sqrt{t}} * v_\ep)_i(y,t) (v_\ep)_j(y,t)\, dy\\
&\quad - \int \nb\bke{\frac1{4\pi|x-y|}}\cdot(\th_\ep\nb G)(y,t)\, dy,
}
where $K_{ij}(x) = \pd_i\pd_j\frac1{4\pi|x|}$. In the last integral  $f_1$ and $f_2$ are absent because $\div f_1=\div f_2=0$.
The right-hand side of \eqref{eq-BT2019-3.11-OB} is defined in $L^2(t_1,t_2; L^2(\R^3))$; for the last integral we use $\th_\ep\nb G\in L^2(t_1,t_2;\dot H^{-1})$ (shown before \eqref{wep-cont}) and that $\nb(-\De)^{-1}$ maps $\dot H^{-1}$ to $L^2$.
Since the only appearance of $\pi_\ep$ in \eqref{eq-BT2019-3.5-OB} is $\nb\pi_\ep$, we can redefine $\pi_\ep$ to equal to $\pi_\ep - \pi_*(t)$ and, therefore, drop $\pi_*(t)$ from \eqref{eq-BT2019-3.11-OB}.

The pressure $\pi_\ep$ given by \eqref{eq-BT2019-3.11-OB} is already bounded in $L^2(t_1,t_2; L^2(\R^3))$ for any $0<t_1<t_2<\infty$ but the bound depends on $t_1, t_2$ and $\ep$.
We now bound it in $L^{3/2}(0,T; L^{3/2}(B_\la))$.
The first and second terms of \eqref{eq-BT2019-3.11-OB} are bounded by the same argument as in the proof of Proposition \ref{BT2019-prop-3.1}. In  particular, we have for any $\ga>0$ (see \eqref{eq4.14})
\EQN{
\int_0^t&\int_{B_\la} \bigg|\tfrac13(\eta_{\ep\sqrt{s}} * v_\ep)\cdot v_\ep(x,s)- \lim_{\de\to0} \int_{| x-y|>\de} K_{ij}(x-y) (\eta_{\ep\sqrt{ s}} * v_\ep)_i(v_\ep)_j(y,s)\, dy \bigg|^{\frac32}dx\, ds\\
&\qquad\qquad\qquad\qquad\le C(\ga) \int_0^t \left(\td\al_\ep(s)^3 + \td\al_\ep(s)\right) ds
+ \ga \int_0^t\int |\nb v_\ep|^2\phi\, dx\,ds.
}
For the last term in \eqref{eq-BT2019-3.11-OB}, we split it into local and nonlocal parts as
\EQN{
p_{\rm near}(x,t) + p_{\rm far}(x,t)
&:= 
\int_{|y|\le\la^2} \nb\bke{\frac1{4\pi|x-y|}}\cdot(\th_\ep\nb G)(y,t)\, dy\\
&\quad + \int_{|y|>\la^2} \nb\bke{\frac1{4\pi|x-y|}}\cdot(\th_\ep\nb G)(y,t)\, dy.
}

The bound of $p_{\rm near}$ is similar to \eqref{eq3.5}: Fix $3<q<10/3$ and choose $r \in ( 60/41,3/2)$ by $\frac 2r=\frac53-\frac1q$. By Young's convolution inequality and H\"older's inequality,
\EQN{
\norm{p_{\rm near}(\cdot,t)}_{L^{3/2}(B_\la)} 
&\lec \norm{|\cdot|^{-2}}_{L^r(B_{\la+\la^2})} \norm{\th_\ep(\cdot,t)\nb G}_{L^{3r/(5r-3)}(B_{\la^2})}  \\
&\lec \norm{\th_\ep(\cdot,t)}_{L^q(B_{\la^2})} \norm{|\cdot|^{-2}}_{L^r(B_{\la^2})}
\lec \norm{\th_\ep(\cdot,t)}_{L^q(B_{\la^2})}.
}
 Note that, by the rescaling \eqref{eq-BT2019-3.6}, the Gagliardo--Nirenberg inequality and Young's inequality, as in the derivation of \eqref{eq-BT2019-3.10},
\[
\norm{\th_\ep(s)}_{L^q(B_{\la^2})}^{3/2} \le C\norm{\phi^{1/2}\th_\ep(\la^{-4}s)}_{L^q}^{3/2}\le C(\ga) (\td\al_\ep(s) + 1) + \ga\norm{\phi^{1/2}\nb\th_\ep(\la^{-4}s)}_{L^2}^2,
\]
since $3<q<10/3$ gives a power of $\td\al_\ep$ less than $1$; the factor $\la^4$ from $\int_0^t g(\la^{-4}s)\,ds\le\la^4\int_0^t g(s)\,ds$ is absorbed into $\ga$. 
It follows that
\EQS{\label{pnear.est}
&\int_0^{ t} \norm{p_{\rm near}(\cdot,s)}_{L^{3/2}(B_\la)}^{3/2} \, ds
\le C  \int_0^{ t} \norm{\th_\ep(\cdot,s)}_{L^q(B_{\la^2})}^{3/2} \, ds\\
&\qquad \le
 C \int_0^t \left(\td\al_\ep(s)+1\right) ds
+ \ga \int_0^t\int |\nb \th_\ep|^2\phi\, dxds,
}
where $C=C(\la,\ga,\eta)$.
 Note that, if we did not assume $\div f_2=0$, the last integral in \eqref{eq-BT2019-3.11-OB} would contain $f_2$, and its contribution to the bound of $p_{\rm near}$ would include
$\int_0^1 \norm{f_2(\cdot,s)}_{L^{2}}^{3/2}\,ds\le C \int_0^1(s^{-3/4})^{3/2}\,ds$, which diverges.

Now we bound the term $p_{\rm far}$. For $x\in B_\la$, we have $|x-y|\sim |y|$ for $|y|\ge\la^2$, and hence
\[
\int_0^t\!\!\int_{B_\la} |p_{\rm far}(x,s)|^{3/2}\,dxds 
\le C\int_0^t\bigg(\int_{|y|\ge\la^2}  |y|^{-4} |\theta_\ep(y,s)| dy\bigg)^{3/2}ds.
\]
By Lemma \ref{lem2.2}, the integrand $|y|^{-4} |\theta_\ep(y,s)|$ is $\la$-DSS with degree 5. By Lemma \ref{lem2.3},
\EQN{
\int_0^t\!\!\int_{B_\la} |p_{\rm far}(x,s)|^{3/2}\,dxds
&\le C\int_0^t\int_{A_0} \bke{ |y|^{-4} |\theta_\ep(y,s)| }^{3/2}dy\,ds
\le C\int_0^t \td\al_\ep(s)^{3/4}\, ds  .
}

We conclude
\EQ{\label{eq4.28}
\int_0^t\!\!\int_{B_\la} |p_{\rm far}(x,s)|^{3/2}\,dxds 
\le C\int_0^t \bke{\td\al_\ep(s)+1}\, ds.
}

After using H\"older's inequality, the same argument deriving \eqref{eq-BT2019-3.10}, the bounds above, and $\al^{3/2}\le\al + \al^3$ for $\al>0$, we get for $0<t<1$
\EQN{
\int_0^t &\norm{\pi_\ep(\cdot,s)}_{L^{3/2}(B_\la)}^{3/2}\, ds + 2\int_0^t \int |\pi_\ep(v_\ep\cdot\nb\phi)|\, dxds\\
&\le C(\la,\ga,\eta)\int_0^t \bke{\td\al_\ep(s)^3 + \td\al_\ep(s) +1} ds + \ga\int_0^t\int \bke{|\nb v_\ep|^2+|\nb \th_\ep|^2}\phi\, dxds.
}

 Finally, we consider the forcing terms in the last line of \eqref{eq-OB-LEE-epsilon}. For $0<t\le1$ we have $\phi(x,t)=\chi^2(|x|)$. Recall Hardy's inequality $\norm{g/|\cdot|}_{L^2}\le 2\norm{\nb g}_{L^2}$ for $g\in H^1(\R^3)$. Hence
\[
\abs{2\int_0^t \int \th_\ep\nb G \cdot v_\ep\phi\, dx\,ds}
\le 2\int_0^t \norm{\frac{\chi\th_\ep}{|\cdot|}}_{L^2}\norm{\frac{\chi v_\ep}{|\cdot|}}_{L^2} ds
\le 8\int_0^t \norm{\nb(\chi\th_\ep)}_{L^2}\norm{\nb(\chi v_\ep)}_{L^2} ds .
\]
We have $\norm{\nb(\chi g)}_{L^2}\le \norm{\phi^{1/2}\nb g}_{L^2}+C\norm{g}_{L^2(B_\la)}$, and, by rescaling and $C_1\ge1$, $\norm{v_\ep(s)}_{L^2(B_\la)}^2+\norm{\th_\ep(s)}_{L^2(B_\la)}^2\le \la\td\al_\ep(s)$. Therefore, for any $\ga>0$,
\EQN{
&\abs{2\int_0^t \int \th_\ep\nb G \cdot v_\ep\phi\, dx\,ds}
\\
&\le\int_0^t  \bke{C_0\norm{\phi^{1/2}\nb v_\ep}_{L^2} \norm{\phi^{1/2}\nb\th_\ep}_{L^2} +\ga\norm{\phi^{1/2}\nb v_\ep}_{L^2}^2+\ga\norm{\phi^{1/2}\nb \th_\ep}_{L^2}^2 +C(\ga)\td\al_\ep}(s)\, ds
\\
&\le \int_0^t \bke{ (\tfrac12+\ga)\norm{\phi^{1/2}\nb v_\ep}_{L^2}^2 + (\tfrac12 C_1+\ga)\norm{\phi^{1/2}\nb\th_\ep}_{L^2}^2 } ds + C(\ga) \int_0^t \td\al_\ep(s)\, ds,
}
where $C_0=8$ and
\EQ{\label{C1.def}
C_1:= C_0^2.
}

 Recall $ M_1=\norm{f_1}_{L^\infty(1,\la^2;L^{6/5,2})} $ and denote $M_2=\norm{f_2}_{L^\infty(1,\la^2;L^{2})} $.
Since $\norm{f_1}_{L^2(0,1;L^{6/5,2})} \le C M_1$ and $\norm{f_2}_{L^1(0,1;L^{2})} \le CM_2$ by Lemma \ref{lem2.4}, for $0<t\le1$ we get, by H\"older's inequality in Lorentz spaces, the Sobolev embedding $\norm{g}_{L^{6,2}}\le C\norm{\nb g}_{L^2}$ and $\phi^2\le\phi$,
\EQN{
\abs{2\int_0^t \int f_1\cdot v_\ep\phi\, dx\,ds}
&\le C\norm{f_1}_{L^2(0,1;L^{6/5,2})}\norm{\nb(\phi v_\ep)}_{L^2(0,t;L^2)}
\le \ga\int_0^t \norm{\nb(\phi v_\ep)}_{L^2}^2\,ds + C\ga^{-1} M_1^2
\\
&\le 2\ga\int_0^t\int|\nb v_\ep|^2\phi\,dx\,ds + C\int_0^t\td\al_\ep(s)\,ds + C\ga^{-1} M_1^2,
}
and
\[
\abs{2\int_0^t \int f_2\cdot v_\ep\phi\, dx\,ds}
\le C\norm{f_2}_{L^1(0,1;L^{2})}\sup_{0\le s\le t}\norm{v_\ep(s)}_{L^2(B_\la)}
\le CM_2\td\al_\ep(t)^{1/2}
\le \tfrac14\td\al_\ep(t)+CM_2^2.
\]
Note that the last bound involves $\td\al_\ep(t)$ itself, not its time integral; it is absorbed below after taking the supremum in time.

Multiplying the second equation \eqref{eq-OB-LEE-epsilon2} by $C_1$ given in \eqref{C1.def}, adding it to the first equation \eqref{eq-OB-LEE-epsilon}, combining the estimates above, and taking $\ga$ sufficiently small to absorb the gradient terms on the right-hand side (including those above \eqref{C1.def}), we obtain for $0<t'\le t\le1$
\EQ{\label{eq-OB-sup-step}
\al_\ep(t') + \int_0^{t'} \int_{B_1} \bke{|\nb v_\ep|^2 + C_1|\nb \th_\ep|^2} dxds
\le \alpha_1 + C_f + C \int_0^t \left(\td\al_\ep(s)^3 + \td\al_\ep(s)\right) ds + \tfrac14\td\al_\ep(t),
}
where $C_f =  C(1 +  M_1^2+M_2^2)$; here we used $t\le1$ to absorb $C\int_0^t 1\,ds$.
Taking the supremum over $t'\in[0,t]$ in \eqref{eq-OB-sup-step} gives $\td\al_\ep(t)\le \frac43\big(\alpha_1 + C_f + C \int_0^t (\td\al_\ep^3 + \td\al_\ep) ds\big)$. Inserting this into \eqref{eq-OB-sup-step} with $t'=t$, we obtain
\EQ{\label{eq-BT2019-3.12-OB}
\al_\ep(t) + \int_0^t \int_{B_1} \bke{|\nb v_\ep|^2 + C_1|\nb \th_\ep|^2} dxds
\le 2(\alpha_1 + C_f) + C \int_0^t \left(\td\al_\ep(s)^3 + \td\al_\ep(s)\right) ds .
}
Therefore, 
\[
\td\al_\ep(t) \le 2(\alpha_1 + C_f) + C \int_0^t \left(\td\al_\ep(s)^3 + \td\al_\ep(s)\right) ds.
\]
By continuity of $\al_\ep(t)$, we have
\[
\td\al_\ep(t) \le  4(\alpha_1 + C_f)\ \text{ for all $t<T$,}
\]
for some $T>0$.
By a continuity argument, we may take $T=\min\{1,(C(2+32(\alpha_1+C_f)^2))^{-1}\}$.

Letting $\ep\to0$ yields
\[
(v,\chi_{B_1}v)(t) + (\th,\chi_{B_1}\th)(t) \le \liminf_{\ep\to0} (v_\ep,\chi_{B_1}v_\ep)_{L^2}(t) + \liminf_{\ep\to0} (\th_\ep,\chi_{B_1}\th_\ep)_{L^2}(t) \le  4(\alpha_1+C_f)
\]
for all $t\le T$.
Note that \eqref{eq-BT2019-3.12-OB} gives uniform (in $\ep$) control of
\[
\int_0^T\int_{B_1}\bke{|\nb v_\ep|^2 + |\nb \th_\ep|^2} dxdt \le C(\al_0,\la, f)
\]
for some constant $C(\al_0,\la, f)$.
From \cite{Tsai-OB} we have that $v_\ep, \th_\ep$ converges weakly to $v, \th$, respectively, in $L^2(1/k, T; H^1(B_1))$ for every $k\in\NN$.
Hence,
\[
\int_{1/k}^T\int_{B_1} \bke{|\nb v|^2 + |\nb \th|^2} dxdt
\le \sup_{\ep>0} \int_0^T\int_{B_1} \bke{|\nb v_\ep|^2 + |\nb \th_\ep|^2} dxdt,
\]
and, letting $k\to\infty$, it follows that
\[
\int_0^T\int_{B_1} \bke{|\nb v|^2 + |\nb \th|^2} dxdt \le C(\al_0,\la, f).
\]
Similarly, since $\pi_\ep\in L^{3/2}(0,T;L^{3/2}( B_\la))$ with uniformly bounded norms, it follows that
$\pi_{\ep_j}\wkto\pi$ weakly in $L^{3/2}(0,T;L^{3/2}(B_\la))$.
 As in the proof of Proposition \ref{BT2019-prop-3.1}, $\pi$ is the limit of the normalized pressures $\pi_\ep$, and the hypotheses of Lemma \ref{BT2019-lem-2.1-OB} on $\td v_k=\eta_{\ep_j\sqrt t}*v_{\ep_j}$ hold. Applying Lemma \ref{BT2019-lem-2.1-OB} yields the desired pressure representation in $L^{3/2}(0,T; L^{3/2}( B_\la))$.

 Finally, $(v,\th)$ satisfies the local energy inequalities \eqref{eq-OB-LEI} with initial data $(v_0,\th_0)$. Indeed, for nonnegative $\phi\in C^\infty_c(\R^3\times[0,\infty))$, let $\ep=\ep_j\to0$ in \eqref{eq-OB-LEE-epsilon2} and in the sum of \eqref{eq-OB-LEE-epsilon} and \eqref{eq-OB-LEE-epsilon2}. The initial terms do not depend on $\ep$. The terms with $\pd_s\phi+\De\phi$, the cubic terms, the pressure term and the forcing term converge, by the strong convergence of $v_\ep$, $\th_\ep$ and $\eta_{\ep\sqrt s}*v_\ep$ in $L^3_\loc$, the weak convergence of $\pi_\ep$ in $L^{5/3}_\loc$, and the weak convergence of $v_\ep$ in $L^2_\loc H^1$ and of $v_\ep(s)$ in $L^2_\loc$ for every $s$. The terms at time $t$ on the left are lower semicontinuous by the weak convergence of $v_\ep(t)$ and $\th_\ep(t)$ for every $t$. The buoyancy term and the gradient terms in the second inequality are handled together by Lemma \ref{lem-buoyancy} with $\mu=1$. This concludes the proof.
\end{proof}

\section{DSS solutions with data in $L^2_\loc$}\label{sec:construction}
In this section, for both 
 the  magnetohydrodynamics equations (MHD) and the Oberbeck-Boussinesq system (OB), after showing that a DSS solution is determined by its value  in a space-time neighborhood of the origin, e.g., $B_1(0)\times (0,1)$, we will construct DSS solutions for DSS data that lie  in $L^2(B_1(0))$ by approximating them with more regular solutions with $L^3_w$  data and using the a priori bounds shown in the previous section.
We will skip the proof for veNSD because its structure is similar to MHD.

\subsection{DSS solutions in a neighborhood of the origin}
As discussed in \cite[Section 4B]{BT-AnalPDE2019}, a DSS solution of the Navier-Stokes equations in a space-time neighborhood of the origin can be extended to a DSS solution on $\R^3\times(0,\infty)$.
In this subsection, we show that the same statement is valid for the MHD equations and the Oberbeck-Boussinesq system.

\subsubsection{MHD equations}
A pair $(V,B)$ is a time-periodic solution with period $T_\la= 2\log\la$ to the Leray system \eqref{eq-MHD-Leray} if and only if $V,B$ are $T_\la$-periodic,  $\div V=\div B=0$ and 
\EQS{\label{eq-BT2019-4.1}
\int_{s'}^{s'+ T_\la} \bke{(V,\pd_sg) - (\nb V,\nb g) + (\frac12 V + \frac12 y\cdot\nb V - V\cdot\nb V + B\cdot\nb B, g)}ds &= 0,\\
\int_{s'}^{s'+ T_\la} \bke{(B,\pd_sh) - (\nb B,\nb h) + (\frac12 B + \frac12 y\cdot\nb B - V\cdot\nb B + B\cdot\nb V, h)}ds &= 0,
}
hold for a fixed $s'\in\R$ and
all $g\in\mathcal{D}_{ T_\la}$ and $h\in\mathcal{D}_{ T_\la}^0$, where 
\EQS{\label{DT.def}
\mathcal{D}_{ T_\la} &=\{ g  \in C^\infty(\R^3\times\R;\R^3):\ g\in \mathcal{D}_{ T_\la}^0, \ \div g=0\},\\[3pt]
\mathcal{D}_{ T_\la}^0 &=\{ h \in C^\infty(\R^3\times\R): \ h(y,s)=h(y,s+ T_\la), \  \exists R>0:\ \supp h(\cdot,s)\subset B_R\ \ \forall s\}.
}
Here $h \in \mathcal{D}_{ T_\la}^0$ is a vector function in \eqref{eq-BT2019-4.1} for MHD, and will be a scalar function in \eqref{eq-BT2019-4.1-OB} for OB.

Since there is a one-to-one correspondence between time-periodic solutions to \eqref{eq-MHD-Leray} and DSS solutions of \eqref{eq-MHD}, an equivalent characterization of DSS solutions is obtained by reformulating \eqref{eq-BT2019-4.1} in the physical variables.
Let 
\EQ{\label{zeg.phih.def}
\ze_g(x,t) = t^{-1}g(y,s),\quad g\in\mathcal{D}_{ T_\la};\quad
\phi_h(x,t) = t^{-1}h(y,s),\quad h\in\mathcal{D}_{ T_\la}^0.
}
Note that $\ze_g(x,t)=\la^2\ze_g(\la x,\la^2 t)$ and $\phi_h(x,t)=\la^2\phi_h(\la x,\la^2 t)$, i.e., they are $(\la,2)$-DSS.
Then,  a $\la$-DSS pair $(v,b)$ with $v,b\in L^\infty_\loc([0,\infty);L^2_\loc)\cap L^2_\loc([0,\infty);H^1_\loc)$  is a solution of \eqref{eq-MHD} if and only if 
\EQS{\label{eq-BT2019-4.2}
\int_t^{\la^2t} \bkt{(v, \pd_t\ze_g) - (\nb v,\na\ze_g) - (v\cdot\nb v - b\cdot\nb b,\ze_g)} d\tau &= 0,\\
\int_t^{\la^2t} \bkt{(b, \pd_t\phi_h) - (\nb b,\na\phi_h) - (v\cdot\nb b - b\cdot\nb v,\phi_h)} d\tau &= 0,
}
for any fixed $t>0$, 
for all $g\in\mathcal{D}_{ T_\la}$ and $h\in\mathcal{D}_{ T_\la}^0$, since \eqref{eq-BT2019-4.1} is just \eqref{eq-BT2019-4.2} in similarity variables.
Note that $(v,\ze_g)|_{\tau=\la^2 t} = (v,\ze_g)|_{\tau=t}$ and $(b,\phi_h)|_{\tau=\la^2 t} = (b,\phi_h)|_{\tau=t}$.

Our goal is to give a third and localized characterization
of DSS solutions.
Fix $k\in\ZZ$ and let $Q_k=B_{\la^k}(0)\times(0,\la^{2k})$.
Let $g\in\mathcal{D}_{ T_\la}$ and $h\in\mathcal{D}_{ T_\la}^0$ be given and $\ze_g$ and $\phi_h$ be as in \eqref{zeg.phih.def}.
Since $\ze_g$ and $\phi_h$ for $1\le t\le \la^2$ are
compactly supported in space, their spatial supports shrink as $t\to0^+$. For $m\in\ZZ$, let 
\EQS{\label{DQkm.def}
\mathcal{D}_{Q_k}^m &= \{ \zeta_g: \ g\in\mathcal{D}_{ T_\la} ,\
\supp\zeta_g |_{t\in[\la^{-2m},\la^{2-2m}]}\subset Q_k\}.
\\[3pt]
\mathcal{D}_{Q_k}^{0,m} &= \{ \phi_h: \ h\in\mathcal{D}_{ T_\la}^0 ,\
\supp\phi_h |_{t\in[\la^{-2m},\la^{2-2m}]}\subset Q_k\}.
}
It is easy to see that $\mathcal{D}_{Q_k}^{m} \subset \mathcal{D}_{Q_k}^{m+1}$ and $\mathcal{D}_{Q_k}^{0,m} \subset \mathcal{D}_{Q_k}^{0,m+1}$ for all $m$ and
\[
\bigcup_{m\in\ZZ}\mathcal{D}_{Q_k}^m = \{\ze_g:\ g\in\mathcal{D}_{ T_\la}\},\quad
\bigcup_{m\in\ZZ}\mathcal{D}_{Q_k}^{0,m} = \{\phi_h:\ h\in\mathcal{D}_{ T_\la}^0\}.
\]
For $\ze_g\in\mathcal{D}_{Q_k}^m$ and $\phi_h\in\mathcal{D}_{Q_k}^{0,m}$,
taking $t=\la^{-2m}$ in \eqref{eq-BT2019-4.2} gives
\EQS{\label{eq-BT2019-4.4}
\int_{\la^{-2m}}^{\la^{2-2m}} \bkt{(v, \pd_t\ze_g) - (\nb v,\na\ze_g) - (v\cdot\nb v - b\cdot\nb b,\ze_g)} d\tau &= 0,\\
\int_{\la^{-2m}}^{\la^{2-2m}} \bkt{(b, \pd_t\phi_h) - (\nb b,\na\phi_h) - (v\cdot\nb b - b\cdot\nb v,\phi_h)} d\tau &= 0.
}
Note that
the integrals are computed over a subset of $Q_k$.
Thus, if $(v,b)$ is a DSS solution to \eqref{eq-MHD}, then it satisfies \eqref{eq-BT2019-4.2} for every $\zeta_g$ and $\phi_h$ with $g \in \mathcal{D}_{ T_\la}$ and $h\in\mathcal{D}_{ T_\la}^0$, and hence it satisfies \eqref{eq-BT2019-4.4} for $m$ sufficiently large. Conversely, if $(v,b)$ is DSS, div-free, and it satisfies \eqref{eq-BT2019-4.4} for every $\zeta_g\in\mathcal{D}_{Q_k}^m$, $\phi_h\in\mathcal{D}_{Q_k}^{0,m}$ and every $m$, then it satisfies \eqref{eq-BT2019-4.2} for every $\zeta_g$ and $\phi_h$ with $g \in \mathcal{D}_{ T_\la}$ and $h\in\mathcal{D}_{ T_\la}^0$. 
This leads to the following local characterization of DSS solutions.

\begin{lem}\label{th5.1}
Let $(v,b):\QQ \to \R^3 \times \R^3$ be $\la$-DSS and $\div v=\div b=0$, with $v,b\in L^\infty_\loc([0,\infty);L^2_\loc(\R^3))\cap L^2_\loc([0,\infty);H^1_\loc(\R^3))$. Let $k \in \ZZ$. Then $(v,b)$
satisfies \eqref{eq-BT2019-4.2} for every $\zeta_g$ and $\phi_h$ with $g \in \mathcal{D}_{ T_\la}$ and $h\in\mathcal{D}_{ T_\la}^0$
 if and only if $(v,b)$ satisfies 
\eqref{eq-BT2019-4.4} for every $m\in\ZZ$, $\zeta_g\in\mathcal{D}_{Q_k}^m$ and $\phi_h\in\mathcal{D}_{Q_k}^{0,m}$.
\end{lem}

In particular, if $(v,b)$ is a DSS solution of \eqref{eq-MHD} in a neighborhood of the origin, then it can be extended to a DSS solution on $\R^3\times(0,\infty)$.

\subsubsection{Oberbeck-Boussinesq system}
A pair $(V,\Th)$ is a time-periodic solution with period $T_\la= 2\log\la$ to the Leray system \eqref{eq-OB-Leray} if and only if $V,\Th$ are $T_\la$-periodic,  $\div V=0$ and 
\EQS{\label{eq-BT2019-4.1-OB}
\int_{s'}^{s'+ T_\la} \bke{(V,\pd_sg) - (\nb V,\nb g) + (\frac12 V + \frac12 y\cdot\nb V - V\cdot\nb V + \Th\nb G + F, g)}ds &= 0,\\
\int_{s'}^{s'+ T_\la} \bke{(\Th,\pd_sh) - (\nb \Th,\nb h) + (\frac12 \Th + \frac12 y\cdot\nb \Th - V\cdot\nb \Th, h)}ds &= 0,
}
hold for a fixed $s'\in\R$ and all $g\in\mathcal{D}_{ T_\la}$ and $h\in\mathcal{D}_{ T_\la}^0$, where $\mathcal{D}_{ T_\la}$ and $\mathcal{D}_{ T_\la}^0$ are defined in \eqref{DT.def}, and $h\in\mathcal{D}_{ T_\la}^0$ is a scalar function in \eqref{eq-BT2019-4.1-OB} for OB.

Since there is a one-to-one correspondence between time-periodic solutions to \eqref{eq-OB-Leray} and DSS solutions of \eqref{eq-OB}, an equivalent characterization of DSS solutions is obtained by reformulating \eqref{eq-BT2019-4.1-OB} in the physical variables.
Let $\ze_g(x,t) $ and $\phi_h(x,t)$ be defined as in \eqref{zeg.phih.def} for $g\in\mathcal{D}_{ T_\la}$ and $h\in\mathcal{D}_{ T_\la}^0$.
Then,  a $\la$-DSS pair $(v,\th)$ with $v,\th\in L^\infty_\loc([0,\infty);L^2_\loc)\cap L^2_\loc([0,\infty);H^1_\loc)$  is a solution of \eqref{eq-OB} if and only if 
\EQS{\label{eq-BT2019-4.2-OB}
\int_t^{\la^2t} \bkt{(v, \pd_t\ze_g) - (\nb v,\na\ze_g) - (v\cdot\nb v - \th\nb G - f,\ze_g)} d\tau &= 0,\\
\int_t^{\la^2t} \bkt{(\th, \pd_t\phi_h) - (\nb \th,\na\phi_h) - (v\cdot\nb \th,\phi_h)} d\tau &= 0,
}
for any fixed $t>0$, for all $g\in\mathcal{D}_{ T_\la}$ and $h\in\mathcal{D}_{ T_\la}^0$, since \eqref{eq-BT2019-4.1-OB} is just \eqref{eq-BT2019-4.2-OB} in similarity variables.
Note that $(v,\ze_g)|_{\tau=\la^2 t} = (v,\ze_g)|_{\tau=t}$ and $(\th,\phi_h)|_{\tau=\la^2 t} = (\th,\phi_h)|_{\tau=t}$.

Our goal is to give a third and localized characterization
of DSS solutions.
Fix $k\in\ZZ$ and let $Q_k=B_{\la^k}(0)\times(0,\la^{2k})$.
Let $\ze_g$ and $\phi_h$ be defined as in \eqref{zeg.phih.def} for
$g\in\mathcal{D}_{ T_\la}$ and $h\in\mathcal{D}_{ T_\la}^0$.
Let $\mathcal{D}_{Q_k}^m $ and $\mathcal{D}_{Q_k}^{0,m}$ be defined as in \eqref{DQkm.def} for $m\in\ZZ$.

For $\ze_g\in\mathcal{D}_{Q_k}^m$ and $\phi_h\in\mathcal{D}_{Q_k}^{0,m}$,
taking $t=\la^{-2m}$ in \eqref{eq-BT2019-4.2-OB} gives 
\EQS{\label{eq-BT2019-4.4-OB}
\int_{\la^{-2m}}^{\la^{2-2m}} \bkt{(v, \pd_t\ze_g) - (\nb v,\na\ze_g) - (v\cdot\nb v - \th\nb G - f,\ze_g)} d\tau &= 0,\\
\int_{\la^{-2m}}^{\la^{2-2m}} \bkt{(\th, \pd_t\phi_h) - (\nb \th,\na\phi_h) - (v\cdot\nb \th,\phi_h)} d\tau &= 0.
}
Note that the integrals are computed over a subset of $Q_k$.
Thus, if $(v,\th)$ is a DSS solution to \eqref{eq-OB}, then it satisfies \eqref{eq-BT2019-4.2-OB} for every $\zeta_g$ and $\phi_h$ with $g \in \mathcal{D}_{ T_\la}$ and $h\in\mathcal{D}_{ T_\la}^0$, and hence it satisfies \eqref{eq-BT2019-4.4-OB} for $m$ sufficiently large. Conversely, if $(v,\th)$ is DSS, $\div v=0$, and it satisfies \eqref{eq-BT2019-4.4-OB} for every $\zeta_g\in\mathcal{D}_{Q_k}^m$, $\phi_h\in\mathcal{D}_{Q_k}^{0,m}$ and every $m$, then it satisfies \eqref{eq-BT2019-4.2-OB} for every $\zeta_g$ and $\phi_h$ with $g \in \mathcal{D}_{ T_\la}$ and $h\in\mathcal{D}_{ T_\la}^0$. 
This leads to the following local characterization of DSS solutions.

\begin{lem}\label{th5.2}
Let $(v,\th):\QQ \to \R^3 \times \R$ be $\la$-DSS and $\div v=0$, with $v,\th\in L^\infty_\loc([0,\infty);L^2_\loc(\R^3))\cap L^2_\loc([0,\infty);H^1_\loc(\R^3))$. Let $k \in \ZZ$. Then $(v,\th)$
satisfies \eqref{eq-BT2019-4.2-OB} for every $\zeta_g$ and $\phi_h$ with $g \in \mathcal{D}_{ T_\la}$ and $h\in\mathcal{D}_{ T_\la}^0$
 if and only if $(v,\th)$ satisfies 
\eqref{eq-BT2019-4.4-OB} for every $m\in\ZZ$, $\zeta_g\in\mathcal{D}_{Q_k}^m$ and $\phi_h\in\mathcal{D}_{Q_k}^{0,m}$.
\end{lem}

In particular, if $(v,\th)$ is a DSS solution of \eqref{eq-OB} in a neighborhood of the origin, then it can be extended to a DSS solution on $\R^3\times(0,\infty)$.

\subsection{Construction of DSS solutions}\label{S:L2loc-data}

Given the similarity between the structures of MHD \eqref{eq-MHD} and veNSD \eqref{eq-veNSD}, we focus on presenting the proofs of Theorem \ref{thm-MHD} for MHD.
The proof of Theorem \ref{thm-veNSD} for veNSD follows the same logic and is left to the readers.
We will prove Theorem \ref{thm-OB} for OB after the proof of Theorem \ref{thm-MHD}.

\begin{proof}[Proof of Theorem \ref{thm-MHD}]
The proof of the theorem is an adaption of the proof of \cite[Theorem 1.2]{BT-AnalPDE2019} for the Navier-Stokes equations to the MHD equations.

Fix $\la>1$ and assume $v_0, b_0\in L^2_\loc$ are divergence-free $\la$-DSS vector fields.
Let $\{v_0^{(k)}\}_{k\in\NN}$ and $\{b_0^{(k)}\}_{k\in\NN}$ be the sequences of divergence-free $\la$-DSS vector fields in $L^3_w(\R^3)$ from \cite[Lemma 4.1]{BT-AnalPDE2019} that approximate $v_0$ and $b_0$, respectively. 
Then, the values $\|v_0^{(k)}\|_{L^2(B_1)}$ and $\|b_0^{(k)}\|_{L^2(B_1)}$ are uniformly bounded in $k$ and 
\[
\|v_0^{(k)}-v_0\|_{L^2(B_1)} + \|b_0^{(k)}-b_0\|_{L^2(B_1)} \to0\quad \text{as}\quad k\to\infty.
\]
Since $v_0^{(k)}, b_0^{(k)}\in L^3_w(\R^3)$ and are $\la$-DSS, by \cite[Theorem 1.2]{Lai-JMFM2019}, whose proof applies with the modifications of Remarks \ref{approx-scheme} and \ref{rem-pib}, there exists a $\la$-DSS local Leray solution $(v_k, b_k)$ to \eqref{eq-MHD} and an associated pressure $\pi_k$ having initial data $(v_0^{(k)}, b_0^{(k)})$ for every $k\in\NN$.  By  Proposition \ref{BT2019-prop-3.1}, the equation for $b_k$ holds without a pressure term, tested with all smooth compactly supported vector fields.
By Proposition \ref{BT2019-prop-3.1}, $v_k$ and $b_k$ are uniformly bounded in $L^\infty(0,T; L^2(B_1))\cap L^2(0,T; H^1(B_1))$, and hence also in $L^{10/3}(0,T;L^{10/3}(B_1))$, for some $T>0$ which depends only on $\la$, $\|v_0^{(k)}\|_{L^2( B_\la)}$, and $\|b_0^{(k)}\|_{L^2( B_\la)}$, which are bounded uniformly in $k$ by the convergence in $L^2(B_1)$ and the DSS scaling.
As usual, see \cite{BT1-AHP, KS07, LR16}, there exist distributions $v$ and $b$ and a subsequence of $\{(v_k, b_k)\}$ (still indexed by $k$ for simplicity) so that $v_k$ and $b_k$ converge to $v$ and $b$, respectively, in weak star topology on $L^\infty(0,T;L^2(B_1))$, in the weak topology on $L^2(0,T;H^1(B_1))$, and in $L^2(0,T;L^2(B_1))$.
Since they are uniformly bounded in $L^{10/3}(0,T;L^{10/3}(B_1))$, they also converge in $L^q(0,T; L^q(B_1))$ for any $q<10/3$. The limit $v$ and $b$ also inherit the DSS and divergence-free properties in $B_1\times (0,T)$.
By the pressure estimate \eqref{eq-BT2019-3.2} in Proposition \ref{BT2019-prop-3.1}, $\pi_k$ are uniformly bounded in $L^{3/2}(0,T;L^{3/2}( B_\la))$ by $C(\la,\norm{v_0}_{L^2(B_\la)}, \norm{b_0}_{L^2(B_\la)})$ and, therefore, we may extract a subsequence which converges weakly to a distribution $\pi\in L^{3/2}(0,T; L^{3/2}( B_\la))$.

Fix $\ka\in\ZZ$ so that $\la^\ka<1$ and $\la^{2\ka}<T$.
Let $Q_\ka = B_{\la^\ka}\times(0,\la^{2\ka})\subset B_1\times(0,T)$.
Then $(v_k, b_k)$ satisfies \eqref{eq-MHD} on $Q_\ka$ and satisfies \eqref{eq-BT2019-4.4} for every $m\in\ZZ$, $\zeta_g\in\mathcal{D}_{Q_\ka}^m$ and $\phi_h\in\mathcal{D}_{Q_\ka}^{ 0,m}$.  Since $v_k\to v$ and $b_k\to b$ in $L^2(0,T;L^2(B_1))$ and weakly in $L^2(0,T;H^1(B_1))$, the quadratic terms converge and the limit $(v,b)$ satisfies the same identities.
Thus, $(v,b)$ can be extended to a DSS solution on $\R^3\times(0,\infty)$ (which we still denote by $(v,b)$).

 We next show that $(v,b)$ satisfies the local energy inequality \eqref{eq-MHD-LEI} and attains the initial data.

First, $v_k$ and $b_k$ are equicontinuous in time in the weak topology, uniformly in $k$. Indeed, for $\psi\in C^\infty_c(B_1;\R^3)$ and $0\le s<t\le T$, the equation for $v_k$ gives
\EQN{
\abs{(v_k(t)-v_k(s),\psi)} &= \abs{\int_s^t \bkt{-(\nb v_k,\nb\psi) + \int (v_k\otimes v_k-b_k\otimes b_k):\nb\psi\,dx + (\pi_k,\div\psi)}d\tau}\\
&\le C(\psi)\Big((t-s)^{\frac12}\norm{\nb v_k}_{L^2(B_1\times(0,T))}\\
&\qquad + (t-s)^{\frac13}\big(\norm{(v_k,b_k)}_{L^3(B_1\times(0,T))}^2 + \norm{\pi_k}_{L^{3/2}(B_1\times(0,T))}\big)\Big)\\
&\le C(\psi)(t-s)^{1/3},
}
with $v_k(0)=v_0^{(k)}$, where $C(\psi)$ is independent of $k$ by Proposition \ref{BT2019-prop-3.1} and $T\le1$. The same holds for $b_k$, whose equation has no pressure term. Since $v_k$ and $b_k$ are also uniformly bounded in $L^\infty(0,T;L^2(B_1))$, after passing to a further subsequence and redefining $v$ and $b$ on a set of times of measure zero, we have
$v_k(t)\wkto v(t)$ and $b_k(t)\wkto b(t)$ weakly in $L^2(B_1)$ for every $t\in[0,T]$, with $v(0)=v_0$ and $b(0)=b_0$. In particular, $v(t)\wkto v_0$ and $b(t)\wkto b_0$ weakly in $L^2(B_1)$ as $t\to0^+$.

Second, by  Proposition \ref{BT2019-prop-3.1}, each $(v_k,b_k)$ satisfies the local energy inequality \eqref{eq-MHD-LEI} with initial data $(v_0^{(k)},b_0^{(k)})$. For nonnegative $\phi\in C^\infty_c(B_1\times[0,T))$, we pass to the limit $k\to\infty$: the initial term converges since $v_0^{(k)}\to v_0$ and $b_0^{(k)}\to b_0$ in $L^2(B_1)$; the other terms on the right converge since $(v_k,b_k)\to(v,b)$ in $L^3(B_1\times(0,T))$ and $\pi_k\wkto\pi$ weakly in $L^{3/2}(B_1\times(0,T))$; and the left side is lower semicontinuous by the weak convergence of $v_k(t)$, $b_k(t)$ for every $t$ and of $\nb v_k$, $\nb b_k$. For general $\phi$, \eqref{eq-MHD-LEI} follows by the DSS scaling.

Third, let $K$ be a compact subset of $B_1$ and choose $\phi\in C^\infty_c(B_1\times[0,T))$, $0\le\phi\le1$, with $\phi=1$ on $K\times[0,T/2]$. The time integrals on the right side of \eqref{eq-MHD-LEI} tend to $0$ as $t\to0^+$, so
\[
\limsup_{t\to0^+}\int \left(|v|^2+|b|^2\right)\phi\,(x,t)\,dx\le \int \left(|v_0|^2+|b_0|^2\right)\phi(x,0)\,dx.
\]
Together with the weak convergence $v(t)\wkto v_0$ and $b(t)\wkto b_0$, expanding $\int(|v(t)-v_0|^2+|b(t)-b_0|^2)\phi(x,t)\,dx$ shows that it tends to $0$. Hence
\[
\lim_{t\to0^+} \bke{\norm{v(t)-v_0}_{L^2(K)} + \norm{b(t)-b_0}_{L^2(K)}} = 0.
\]
For a general compact subset $K$ of $\R^3$, we have $K' = \la^mK\subset B_1$ for some $m\in\ZZ$, and 
\EQN{
\int_K |v(x,t) - v_0(x)|^2\, dx &= \la^{-m} \int_{K'} |v(x',\la^{2m}t) - v_0(x')|^2\, dx',\\
\int_K |b(x,t) - b_0(x)|^2\, dx &= \la^{-m} \int_{K'} |b(x',\la^{2m}t) - b_0(x')|^2\, dx'.
}
It follows that $\lim_{t\to0^+} \norm{v(t) - v_0}_{L^2(K)} + \norm{b(t) - b_0}_{L^2(K)} = 0$ for every compact set $K\subset\R^3$.
A similar rescaling argument also implies that $v,b \in L^\infty(0,T';L^2(K))\cap L^2(0,T';H^1(K))$ and $\pi\in L^{3/2}(0,T';L^{3/2}(K))$ for any $T'>0$ and compact subset $K$ of $\R^3$, with the bound \eqref{MHD-apriori}.

Finally, note that $\pi_k$ satisfies formula \eqref{eq-BT2019-3.3}.
Applying Lemma \ref{BT2019-lem-2.1} to the sequence and limit above implies that $\pi$ satisfies the desired pressure formula in $L^{3/2}(0,T; L^{3/2}( B_\la))$.
Rescaling establishes the formula in $L^{3/2}_\loc(\R^3\times[0,\infty))$.
 By the same rescaling, $(v,b,\pi)$ satisfies \eqref{eq-MHD} in the sense of distributions on $\QQ$: a DSS rescaling maps the support of any test function in $\QQ$ into $B_1\times(0,T)$, where the equations hold as the limit of those for $(v_k,b_k,\pi_k)$.
\end{proof}

\begin{proof}[Proof of Theorem \ref{thm-OB}]
Fix $\la>1$ and assume $v_0, \th_0\in L^2_\loc$  are a divergence-free vector field and a function, both $\la$-DSS.
Let $\{v_0^{(k)}\}$ be the sequence of divergence-free $\la$-DSS  vector fields from \cite[Lemma 4.1]{BT-AnalPDE2019} that approximates $v_0$.
By adapting the proof of \cite[Lemma 4.1]{BT-AnalPDE2019} for divergence-free vector fields to scalar functions, we obtain a sequence of $\la$-DSS functions $\{\th_0^{(k)}\}\subset L^3_w(\R^3)$
that approximates $\theta_0$. We have
\[
\|v_0^{(k)}-v_0\|_{L^2(B_1)} + \|\th_0^{(k)}-\th_0\|_{L^2(B_1)} \to0\quad \text{as}\quad k\to\infty.
\]
Since $v_0^{(k)}, \th_0^{(k)}\in L^3_w$ and are $\la$-DSS, by \cite[Theorem 1.1]{Tsai-OB} (see also \cite{Tsai-OB-corr}), whose proof applies with the modification of Remark \ref{approx-scheme}, there exists a $\la$-DSS distributional solution $(v_k, \th_k)$ to \eqref{eq-OB} and an associated pressure $\pi_k$ having initial data $(v_0^{(k)}, \th_0^{(k)})$ for every $k\in\NN$; by Proposition \ref{BT2019-prop-3.1-OB}, $(v_k,\th_k)$ satisfies  the local energy inequalities  \eqref{eq-OB-LEI} with initial data $(v_0^{(k)}, \th_0^{(k)})$.
By Proposition \ref{BT2019-prop-3.1-OB}, $v_k$ and $\th_k$ are uniformly bounded in $L^\infty(0,T; L^2(B_1))\cap L^2(0,T; H^1(B_1))$ (hence also in $L^{10/3}(0,T;L^{10/3}(B_1))$) for some $T$ which depends only on $\la$, $\|v_0^{(k)}\|_{L^2( B_\la)}$, and $\|\th_0^{(k)}\|_{L^2( B_\la)}$, which are bounded uniformly in $k$.
As usual, see \cite{BT1-AHP, KS07, LR16}, there exist distributions $v$ and $\th$ and a subsequence of $\{(v_k, \th_k)\}$ (still indexed by $k$ for simplicity) so that $v_k$ and $\th_k$ converge to $v$ and $\th$, respectively, in weak star topology on $L^\infty(0,T;L^2(B_1))$, in the weak topology on $L^2(0,T;H^1(B_1))$, and in $L^2(0,T;L^2(B_1))$.
Since they are uniformly bounded in $L^{10/3}(0,T;L^{10/3}(B_1))$, they also converge in $L^q(0,T; L^q(B_1))$ for any $q<10/3$. 
The limit $v$ and $\th$ also inherit the DSS and divergence-free properties in $B_1\times (0,T)$.
By the pressure estimate \eqref{eq-BT2019-3.2-OB} in Proposition \ref{BT2019-prop-3.1-OB}, $\pi_k$ are uniformly bounded in $L^{3/2}(0,T;L^{3/2}( B_\la))$ by $C(\la,\norm{v_0}_{L^2(B_\la)}, \norm{\th_0}_{L^2(B_\la)}, f)$ and, therefore, we may extract a subsequence which converges weakly to a distribution $\pi\in L^{3/2}(0,T; L^{3/2}( B_\la))$.

Fix $\ka\in\ZZ$ so that $\la^\ka<1$ and $\la^{2\ka}<T$.
Then, $Q_\ka = B_{\la^\ka}\times(0,\la^{2\ka})\subset B_1\times(0,T)$.
Therefore $(v_k, \th_k)$ satisfies \eqref{eq-OB} on $Q_\ka$ and satisfies \eqref{eq-BT2019-4.4-OB} for every $m\in\ZZ$, 
$\zeta_g\in\mathcal{D}_{Q_\ka}^m$ and $\phi_h\in\mathcal{D}_{Q_\ka}^{0,m}$.
Passing to the limit as in the proof of Theorem \ref{thm-MHD}, the limit $(v,\th)$ satisfies the same identities; the buoyancy term $(\th_k\nb G,\ze_g)$ is linear in $\th_k$ and converges, since $\th_k\to\th$ in $L^q(B_1\times(0,T))$ for some $q>3$ and $|x|^{-2}\ze_g\in L^{q'}$.
Thus, $(v,\th)$ can be extended to a DSS solution on $\R^3\times(0,\infty)$ (which we still denote by $(v,\th)$).

 We next show that $(v,\th)$ satisfies \eqref{eq-OB-LEI} and attains the initial data. As in the proof of Theorem \ref{thm-MHD}, $v_k$ and $\th_k$ are equicontinuous in time in the weak topology, uniformly in $k$: for $\psi\in C^\infty_c(B_1;\R^3)$ and $0\le s<t\le T$,
\[
\abs{(v_k(t)-v_k(s),\psi)}\le C(\psi)(t-s)^{1/4}.
\]
Compared with the MHD case, we also use, for some fixed $3<q<10/3$,
\[
\int_s^t\abs{(\th_k\nb G,\psi)}\,d\tau\le (t-s)^{1-\frac1q}\norm{\th_k}_{L^q(B_1\times(0,T))}\norm{|\cdot|^{-2}\psi}_{L^{q'}}
\]
and $\int_s^t\abs{(f,\psi)}\,d\tau\le C(\psi)( M_1+M_2)(t-s)^{1/4}$ by Lemma \ref{lem2.4}, where $ M_1$ and $M_2$ are as in the proof of Proposition \ref{BT2019-prop-3.1-OB}. Similarly, $\abs{(\th_k(t)-\th_k(s),\psi)}\le C(\psi)(t-s)^{1/3}$ for $\psi\in C^\infty_c(B_1)$. Hence, after passing to a further subsequence and redefining $v$ and $\th$ on a set of times of measure zero, $v_k(t)\wkto v(t)$ and $\th_k(t)\wkto\th(t)$ weakly in $L^2(B_1)$ for every $t\in[0,T]$, with $v(0)=v_0$ and $\th(0)=\th_0$.

We pass to the limit $k\to\infty$ in \eqref{eq-OB-LEI} for $(v_k,\th_k)$, for nonnegative $\phi\in C^\infty_c(B_1\times[0,T))$, as at the end of the proof of Proposition \ref{BT2019-prop-3.1-OB}: the initial terms converge since $v_0^{(k)}\to v_0$ and $\th_0^{(k)}\to\th_0$ in $L^2(B_1)$; the other terms converge except the buoyancy term, which is handled together with the gradient terms by Lemma \ref{lem-buoyancy} with $\mu=1$; and the terms at time $t$ on the left are lower semicontinuous. For general $\phi$, \eqref{eq-OB-LEI} follows by the DSS scaling.

For the initial data, let $K$ be a compact subset of $B_1$ and choose $\phi\in C^\infty_c(B_1\times[0,T))$, $0\le\phi\le1$, with $\phi=1$ on $K\times[0,T/2]$. The time integrals in \eqref{eq-OB-LEI} tend to $0$ as $t\to0^+$; for the buoyancy term this follows from Hardy's inequality, and for the forcing term from Lemma \ref{lem2.4}. The first inequality in \eqref{eq-OB-LEI} and the weak convergence $\th(t)\wkto\th_0$ give $\int(\th(t)-\th_0)^2\phi(x,t)\,dx\to0$. Then the second inequality gives
$\limsup_{t\to0^+}\int|v|^2\phi\,(x,t)\,dx\le\int|v_0|^2\phi(x,0)\,dx$, and together with $v(t)\wkto v_0$ we get $\int|v(t)-v_0|^2\phi(x,t)\,dx\to0$. Hence
\[
\lim_{t\to0^+} \bke{\norm{v(t)-v_0}_{L^2(K)} + \norm{\th(t)-\th_0}_{L^2(K)}} = 0.
\]
For a general compact subset $K$ of $\R^3$, we have $K' = \la^mK\subset B_1$ for some $m\in\ZZ$, and 
\EQN{
\int_K |v(x,t) - v_0(x)|^2\, dx &= \la^{-m} \int_{K'} |v(x',\la^{2m}t) - v_0(x')|^2\, dx',\\
\int_K (\th(x,t) - \th_0(x))^2\, dx &= \la^{-m} \int_{K'} (\th(x',\la^{2m}t) - \th_0(x'))^2\, dx'.
}
It follows that $\lim_{t\to0^+} \norm{v(t) - v_0}_{L^2(K)} + \norm{\th(t) - \th_0}_{L^2(K)} = 0$ for every compact set $K\subset\R^3$.
A similar rescaling argument also implies that $v,\th \in L^\infty(0,T';L^2(K))\cap L^2(0,T';H^1(K))$ and $\pi\in L^{3/2}(0,T';L^{3/2}(K))$ for any $T'>0$ and compact subset $K$ of $\R^3$, with the bound \eqref{OB-apriori}.

Finally, note that $\pi_k$ satisfies formula \eqref{eq-BT2019-3.3-OB}.
Applying Lemma \ref{BT2019-lem-2.1-OB} to the sequence and limit above implies that $\pi$ satisfies the desired pressure formula in $L^{3/2}(0,T; L^{3/2}( B_\la))$.
Rescaling establishes the formula in $L^{3/2}_\loc(\R^3\times[0,\infty))$.
 By the same rescaling, $(v,\th,\pi)$ satisfies \eqref{eq-OB} in the sense of distributions on $\QQ$.
\end{proof}

\subsection{Self-similar solutions}\label{S:SS}
 In this subsection we prove the statements for self-similar data in Theorems \ref{thm-MHD}--\ref{thm-OB}, following \cite[Section 5.1]{BT1-AHP}, where self-similar solutions of the Navier--Stokes equations are obtained as limits of DSS solutions. We give the details for MHD; the proofs for veNSD and OB are the same.

Suppose that $(v_0,b_0)$ is self-similar. Then it is $\la$-DSS for every $\la>1$. Let
\[
\la_k=2^{2^{-k}},\qquad k\in\NN,
\]
so that $\la_{k+1}^2=\la_k$, $\la_k^{2^k}=2$ and $\la_k\to1$. By iterating \eqref{DSS.defn}, a $\la$-DSS function is also $\la^m$-DSS for every $m\in\NN$. In particular, every $\la_k$-DSS function is $2$-DSS.

Let $(v_k,b_k)$, with pressure $\pi_k$, be the $\la_k$-DSS solution given by Theorem \ref{thm-MHD} with factor $\la_k$. All functions used in its construction are $2$-DSS: the approximations of the data are $\la_k$-DSS, and the time-periodic solutions of the mollified scheme in the proof of Proposition \ref{BT2019-prop-3.1}, which have period $2\log\la_k$, also have period $2\log2$. Hence the proofs of Proposition \ref{BT2019-prop-3.1} and of Theorem \ref{thm-MHD} apply with $\la=2$, and give the bound \eqref{MHD-apriori} for $(v_k,b_k,\pi_k)$ with a constant $C(T,E,2,\al_0)$, $\al_0=\norm{v_0}_{L^2(B_2)}^2+\norm{b_0}_{L^2(B_2)}^2$, which is independent of $k$.

With these uniform bounds, we pass to the limit $k\to\infty$ as in the proof of Theorem \ref{thm-MHD}. Each $(v_k,b_k,\pi_k)$ is a $2$-DSS distributional solution of \eqref{eq-MHD}, in which the equation for $b_k$ holds without a pressure term, and it satisfies \eqref{eq-MHD-LEI} with the same initial data $(v_0,b_0)$ and the pressure formula \eqref{eq-pressure-MHD}. The weak equicontinuity, the passage to the limit in the equations and in \eqref{eq-MHD-LEI}, the convergence to the initial data, and Lemma \ref{BT2019-lem-2.1} with $\la=2$ apply verbatim. Since the bounds hold on $E\times(0,T')$ for every compact $E$ and every $T'>0$, a diagonal subsequence converges on every such set, and no extension argument is needed. This gives a $2$-DSS solution $(v,b,\pi)$ with all the properties stated in Theorem \ref{thm-MHD}.

It remains to show that $(v,b,\pi)$ is self-similar. For $l\ge k$, we have $\la_k=\la_l^{2^{l-k}}$, so $(v_l,b_l,\pi_l)$ is $\la_k$-DSS. The identity \eqref{DSS.defn} with a fixed factor is preserved under weak limits, so $(v,b,\pi)$ is $\la_k$-DSS for every $k$. Hence $v_\la=v$ for every $\la$ in the set $\{\la_k^j:k\in\NN,\ j\in\ZZ\}=\{2^{j2^{-k}}\}$, which is dense in $(0,\infty)$; here $v_\la(x,t)=\la v(\la x,\la^2t)$. Since $\la\mapsto v_\la$ is continuous from $(0,\infty)$ to $L^2_\loc(\QQ)$ (approximate $v$ by continuous functions), we get $v_\la=v$ for every $\la>0$. The same argument applies to $b$, and to $\pi$ in $L^{3/2}_\loc(\QQ)$ with $\pi_\la(x,t)=\la^2\pi(\la x,\la^2t)$. Thus $(v,b,\pi)$ is self-similar.

For OB, $f_1$ and $f_2$ are self-similar of degree $3$, hence $(\la_k,3)$-DSS for every $k$, and the constants $M_1$ and $M_2$ in Theorem \ref{thm-OB} and Proposition \ref{BT2019-prop-3.1-OB} are computed with $\la=2$, so they are independent of $k$ as well.

\section*{Acknowledgments}
CL gratefully acknowledges the support by the Department of Mathematical Sciences and the Office of Research and Sponsored Programs (ORSP) at Middle Tennessee State University.
The research of TT was partially supported by Natural Sciences and Engineering Research Council of Canada (NSERC) under grant RGPIN-2023-04534.

\addcontentsline{toc}{section}{References}
\small
\bibliographystyle{abbrv}
\bibliography{coupled-fluid-2026}

@article {MR3916974,
    AUTHOR = {Albritton, Dallas and Barker, Tobias},
     TITLE = {Global weak {B}esov solutions of the {N}avier-{S}tokes
              equations and applications},
   JOURNAL = {Arch. Ration. Mech. Anal.},
  FJOURNAL = {Archive for Rational Mechanics and Analysis},
    VOLUME = {232},
      YEAR = {2019},
    NUMBER = {1},
     PAGES = {197--263},
      ignore-ISSN = {0003-9527,1432-0673},
   MRCLASS = {35Q30 (35B44 35C06 35D30 76D03 76D05)},
  MRNUMBER = {3916974},
MRREVIEWER = {Tao\ Tao},
       DOI = {10.1007/s00205-018-1319-0},
       URL = {https://doi.org/10.1007/s00205-018-1319-0},
}

@article {BCT,
    AUTHOR = {Bradshaw, Zachary and Chernobai, Misha and Tsai, Tai-Peng},
     TITLE = {Global {N}avier-{S}tokes flows in intermediate spaces},
   JOURNAL = {J. Differential Equations},
  FJOURNAL = {Journal of Differential Equations},
    VOLUME = {429},
      YEAR = {2025},
     PAGES = {50--87},
      ignore-ISSN = {0022-0396,1090-2732},
   MRCLASS = {35Q30 (35D30 76D03 76D05)},
  MRNUMBER = {4866502},
       DOI = {10.1016/j.jde.2025.02.025},
       URL = {https://doi.org/10.1016/j.jde.2025.02.025},
}

@article {MR4395595,
    AUTHOR = {Bradshaw, Zachary and Kukavica, Igor and Tsai, Tai-Peng},
     TITLE = {Existence of global weak solutions to the {N}avier-{S}tokes
              equations in weighted spaces},
   JOURNAL = {Indiana Univ. Math. J.},
  FJOURNAL = {Indiana University Mathematics Journal},
    VOLUME = {71},
      YEAR = {2022},
    NUMBER = {1},
     PAGES = {191--212},
      ignore-ISSN = {0022-2518,1943-5258},
   MRCLASS = {35Q30},
  MRNUMBER = {4395595},
       DOI = {10.1512/iumj.2022.71.8789},
       URL = {https://doi.org/10.1512/iumj.2022.71.8789},
}

@article {BLT,
    AUTHOR = {Bradshaw, Zachary and Lai, Chen-Chih and Tsai, Tai-Peng},
     TITLE = {Mild solutions and spacetime integral bounds for {S}tokes and
              {N}avier-{S}tokes flows in {W}iener amalgam spaces},
   JOURNAL = {Math. Ann.},
  FJOURNAL = {Mathematische Annalen},
    VOLUME = {388},
      YEAR = {2024},
    NUMBER = {3},
     PAGES = {3053--3126},
      ISSN = {0025-5831,1432-1807},
   MRCLASS = {35Q30 (76D05 76D07)},
  MRNUMBER = {4705760},
MRREVIEWER = {Mengying\ Xiao},
       DOI = {10.1007/s00208-023-02607-9},
       URL = {https://doi-org.ezproxy.library.ubc.ca/10.1007/s00208-023-02607-9},
}

@article {BT1-AHP,
    AUTHOR = {Bradshaw, Zachary and Tsai, Tai-Peng},
     TITLE = {Forward discretely self-similar solutions of the
              {N}avier-{S}tokes equations {II}},
   JOURNAL = {Ann. Henri Poincar\'{e}},
  FJOURNAL = {Annales Henri Poincar\'{e}. A Journal of Theoretical and
              Mathematical Physics},
    VOLUME = {18},
      YEAR = {2017},
    NUMBER = {3},
     PAGES = {1095--1119},
      ISSN = {1424-0637},
   MRCLASS = {35Q30 (35C06 76D05)},
  MRNUMBER = {3611025},
MRREVIEWER = {Reinhard Farwig},
       DOI = {10.1007/s00023-016-0519-0},
       URL = {https://doi.org/10.1007/s00023-016-0519-0},
}

@article {BT2,
    AUTHOR = {Bradshaw, Zachary and Tsai, Tai-Peng},
     TITLE = {Rotationally corrected scaling invariant solutions to the
              {N}avier-{S}tokes equations},
   JOURNAL = {Comm. Partial Differential Equations},
  FJOURNAL = {Communications in Partial Differential Equations},
    VOLUME = {42},
      YEAR = {2017},
    NUMBER = {7},
     PAGES = {1065--1087},
      ISSN = {0360-5302},
   MRCLASS = {35Q30 (35C06 76D05 76U05)},
  MRNUMBER = {3691390},
       URL = {https://doi.org/10.1080/03605302.2017.1323922},
}

@article {BT3,
    AUTHOR = {Bradshaw, Zachary and Tsai, Tai-Peng},
     TITLE = {Discretely self-similar solutions to the {N}avier-{S}tokes
              equations with {B}esov space data},
   JOURNAL = {Arch. Ration. Mech. Anal.},
  FJOURNAL = {Archive for Rational Mechanics and Analysis},
    VOLUME = {229},
      YEAR = {2018},
    NUMBER = {1},
     PAGES = {53--77},
      ISSN = {0003-9527,1432-0673},
   MRCLASS = {35Q30 (35C06 76D05)},
  MRNUMBER = {3799090},
MRREVIEWER = {Sergey\ V.\ Ershkov},
       DOI = {10.1007/s00205-017-1213-1},
       URL = {https://doi.org/10.1007/s00205-017-1213-1},
}

@article{BT-AnalPDE2019,
  author={Bradshaw, Zachary and Tsai, Tai-Peng},
  title={Discretely self-similar solutions to the {N}avier--{S}tokes equations with data in {$L^2_{\rm loc}$} satisfying the local energy inequality},
  journal={Anal. PDE},
  fjournal={Analysis \& PDE},
  volume={12},
  number={8},
  pages={1943--1962},
  year={2019},
  publisher={Mathematical Sciences Publishers}
}

@article {BT-SIMA21,
    AUTHOR = {Bradshaw, Zachary and Tsai, Tai-Peng},
     TITLE = {Local energy solutions to the {N}avier-{S}tokes equations in
              {W}iener amalgam spaces},
   JOURNAL = {SIAM J. Math. Anal.},
  FJOURNAL = {SIAM Journal on Mathematical Analysis},
    VOLUME = {53},
      YEAR = {2021},
    NUMBER = {2},
     PAGES = {1993--2026},
      ISSN = {0036-1410,1095-7154},
   MRCLASS = {35Q30 (76D03 76D05)},
  MRNUMBER = {4242337},
MRREVIEWER = {Fuyi\ Xu},
       DOI = {10.1137/20M1361420},
       URL = {https://doi-org.ezproxy.library.ubc.ca/10.1137/20M1361420},
}

@article {MR4848568,
    AUTHOR = {Brandolese, Lorenzo and Karch, Grzegorz},
     TITLE = {Large self-similar solutions to {O}berbeck-{B}oussinesq system
              with {N}ewtonian gravitational field},
   JOURNAL = {J. Differential Equations},
  FJOURNAL = {Journal of Differential Equations},
    VOLUME = {424},
      YEAR = {2025},
     PAGES = {139--158},
      ISSN = {0022-0396,1090-2732},
   MRCLASS = {35Q30 (35C06 35J47 76D03)},
  MRNUMBER = {4848568},
MRREVIEWER = {Jiaqi\ Yang},
       DOI = {10.1016/j.jde.2024.12.041},
       URL = {https://doi-org.ezproxy.library.ubc.ca/10.1016/j.jde.2024.12.041},
}

@article{CKN-CPAM1982,
  title={Partial regularity of suitable weak solutions of the {N}avier-{S}tokes equations},
  author={Caffarelli, Luis and Kohn, Robert and Nirenberg, Louis},
  journal={Comm. Pure Appl. Math.},
  fjournal={Communications on p=Pure and Applied Mathematics},
  volume={35},
  number={6},
  pages={771--831},
  year={1982},
  publisher={Wiley Online Library}
}

@article{Chae-Wolf,
    AUTHOR = {Chae, Dongho and Wolf, J\"org},
     TITLE = {Existence of discretely self-similar solutions to the
              {N}avier-{S}tokes equations for initial value in
              {$L_{loc}^2(\Bbb R^3)$}},
   JOURNAL = {Ann. Inst. H. Poincar\'e{} C Anal. Non Lin\'eaire},
  FJOURNAL = {Annales de l'Institut Henri Poincar\'e{} C. Analyse Non
              Lin\'eaire},
    VOLUME = {35},
      YEAR = {2018},
    NUMBER = {4},
     PAGES = {1019--1039},
      ISSN = {0294-1449,1873-1430},
   MRCLASS = {76D03 (35C06 35Q30 76D05)},
  MRNUMBER = {3795025},
MRREVIEWER = {Isabelle\ Gruais},
       DOI = {10.1016/j.anihpc.2017.10.001},
       URL = {https://doi.org/10.1016/j.anihpc.2017.10.001},
}

@article {MR4205085,
    AUTHOR = {Fern\'andez-Dalgo, Pedro Gabriel and Jarr\'in, Oscar},
     TITLE = {Discretely self-similar solutions for 3{D} {MHD} equations and
              global weak solutions in weighted {$L^2$} spaces},
   JOURNAL = {J. Math. Fluid Mech.},
  FJOURNAL = {Journal of Mathematical Fluid Mechanics},
    VOLUME = {23},
      YEAR = {2021},
    NUMBER = {1},
     PAGES = {Paper No. 22, 30},
      ISSN = {1422-6928,1422-6952},
   MRCLASS = {35Q30 (76W05)},
  MRNUMBER = {4205085},
MRREVIEWER = {Dongfen\ Bian},
       DOI = {10.1007/s00021-020-00551-1},
       URL = {https://doi-org.ezproxy.mtsu.edu/10.1007/s00021-020-00551-1},
}

@article {MR4090469,
    AUTHOR = {Fern\'andez-Dalgo, Pedro Gabriel and Lemari\'e-Rieusset,
              Pierre Gilles},
     TITLE = {Weak solutions for {N}avier-{S}tokes equations with initial
              data in weighted {$L^2$} spaces},
   JOURNAL = {Arch. Ration. Mech. Anal.},
  FJOURNAL = {Archive for Rational Mechanics and Analysis},
    VOLUME = {237},
      YEAR = {2020},
    NUMBER = {1},
     PAGES = {347--382},
      ISSN = {0003-9527,1432-0673},
   MRCLASS = {35Q30 (35D30 76D03 76D05)},
  MRNUMBER = {4090469},
       DOI = {10.1007/s00205-020-01510-w},
       URL = {https://doi-org.ezproxy.library.ubc.ca/10.1007/s00205-020-01510-w},
}

@article {JiaSve,
    AUTHOR = {Jia, Hao and {\v{S}}ver{\'a}k, Vladim{\'{\i}}r},
     TITLE = {Local-in-space estimates near initial time for weak solutions
              of the {N}avier-{S}tokes equations and forward self-similar
              solutions},
   JOURNAL = {Invent. Math.},
  FJOURNAL = {Inventiones Mathematicae},
    VOLUME = {196},
      YEAR = {2014},
    NUMBER = {1},
     PAGES = {233--265},
      ISSN = {0020-9910},
   MRCLASS = {35Q30 (35A01 35D30 76D03 76D05)},
  MRNUMBER = {3179576},
MRREVIEWER = {Boris Muha},
       DOI = {10.1007/s00222-013-0468-x},
       URL = {http://dx.doi.org/10.1007/s00222-013-0468-x},
}

@incollection {KS07,
    AUTHOR = {Kikuchi, N. and Seregin, Gregory},
     TITLE = {Weak solutions to the {C}auchy problem for the
              {N}avier-{S}tokes equations satisfying the local energy
              inequality},
 BOOKTITLE = {Nonlinear equations and spectral theory},
    SERIES = {Amer. Math. Soc. Transl. Ser. 2},
    VOLUME = {220},
     PAGES = {141--164},
 PUBLISHER = {Amer. Math. Soc., Providence, RI},
      YEAR = {2007},
   MRCLASS = {35Q30 (35B65 35D05 76D03 76D05)},
  MRNUMBER = {2343610 (2008k:35359)},
MRREVIEWER = {J{\"u}rgen Socolowsky},
}

@article {Korobkov-Tsai,
    AUTHOR = {Korobkov, Mikhail V. and Tsai, Tai-Peng},
     TITLE = {Forward self-similar solutions of the {N}avier-{S}tokes
              equations in the half space},
   JOURNAL = {Anal. PDE},
  FJOURNAL = {Analysis {\&} PDE},
    VOLUME = {9},
      YEAR = {2016},
    NUMBER = {8},
     PAGES = {1811--1827},
      ignore-ISSN = {2157-5045},
   MRCLASS = {35Q30 (35C06 76D05)},
  MRNUMBER = {3599519},
       DOI = {10.2140/apde.2016.9.1811},
       URL = {http://dx.doi.org/10.2140/apde.2016.9.1811},
}

@article {KwTs,
    AUTHOR = {Kwon, Hyunju and Tsai, Tai-Peng},
     TITLE = {Global {N}avier-{S}tokes flows for non-decaying initial data
              with slowly decaying oscillation},
   JOURNAL = {Comm. Math. Phys.},
  FJOURNAL = {Communications in Mathematical Physics},
    VOLUME = {375},
      YEAR = {2020},
    NUMBER = {3},
     PAGES = {1665--1715},
      ignore-ISSN = {0010-3616,1432-0916},
   MRCLASS = {76D05 (35Q30)},
  MRNUMBER = {4091509},
       DOI = {10.1007/s00220-020-03695-3},
       URL = {https://doi.org/10.1007/s00220-020-03695-3},
}

@article {Lai-JMFM2019,
    AUTHOR = {Lai, Chen-Chih},
     TITLE = {Forward discretely self-similar solutions of the {MHD}
              equations and the viscoelastic {N}avier-{S}tokes equations
              with damping},
   JOURNAL = {J. Math. Fluid Mech.},
  FJOURNAL = {Journal of Mathematical Fluid Mechanics},
    VOLUME = {21},
      YEAR = {2019},
    NUMBER = {3},
     PAGES = {Paper No. 38, 28},
      ISSN = {1422-6928,1422-6952},
   MRCLASS = {76D03 (35Q35 76A10 76W05)},
  MRNUMBER = {3975833},
MRREVIEWER = {Dongfen\ Bian},
       DOI = {10.1007/s00021-019-0443-4},
       URL = {https://doi.org/10.1007/s00021-019-0443-4},
}

@book {LR02,
    AUTHOR = {Lemari\'e-Rieusset, Pierre Gilles},
     TITLE = {Recent developments in the {N}avier-{S}tokes problem},
    SERIES = {Chapman \& Hall/CRC Research Notes in Mathematics},
    VOLUME = {431},
 PUBLISHER = {Chapman \& Hall/CRC, Boca Raton, FL},
      YEAR = {2002},
     PAGES = {xiv+395},
      ISBN = {1-58488-220-4},
   MRCLASS = {35Q30 (46N20 76D03 76D05)},
  MRNUMBER = {1938147 (2004e:35178)},
MRREVIEWER = {Joel David Avrin},
       DOI = {10.1201/9781420035674},
       URL = {http://dx.doi.org/10.1201/9781420035674},
}

@book {LR16,
    AUTHOR = {Lemari\'e-Rieusset, Pierre Gilles},
     TITLE = {The {N}avier-{S}tokes problem in the 21st century},
 PUBLISHER = {CRC Press, Boca Raton, FL},
      YEAR = {2016},
     PAGES = {xxii+718},
      ISBN = {978-1-4665-6621-7},
   MRCLASS = {76D05 (35Q30 76D03)},
  MRNUMBER = {3469428},
MRREVIEWER = {Isabelle Gruais},
       DOI = {10.1201/b19556},
       URL = {http://dx.doi.org/10.1201/b19556},
}

@article {Tsai2013,
    AUTHOR = {Tsai, Tai-Peng},
     TITLE = {Forward discretely self-similar solutions of the
                  {N}avier-{S}tokes equations},
   JOURNAL = {Comm. Math. Phys.},
  FJOURNAL = {Communications in Mathematical Physics},
    VOLUME = {328},
      YEAR = {2014},
    NUMBER = {1},
     PAGES = {29--44},
      ISSN = {0010-3616},
   MRCLASS = {Preliminary Data},
  MRNUMBER = {3196979},
       DOI = {10.1007/s00220-014-1984-2},
       URL = {http://dx.doi.org/10.1007/s00220-014-1984-2},
}

@article{Tsai-OB,
  title={Large discretely self-similar solutions to {O}berbeck--{B}oussinesq system with {N}ewtonian gravitational field},
  author={Tsai, Tai-Peng},
  journal={Funkcialaj Ekvacioj},
  volume={69},
  number={1},
  pages={101--118},
  year={2026},
  doi={10.1619/fesi.69.101}
}

@preamble{"\providecommand{\noopsort}[1]{}"}

@unpublished{Tsai-OB-corr,
  title={Correction to: {L}arge discretely self-similar solutions to {O}berbeck--{B}oussinesq system with {N}ewtonian gravitational field},
  author={Tsai, Tai-Peng},
  note={preprint: \url{https://personal.math.ubc.ca/~ttsai/publications/Boussinesq-FE-errata.pdf}},
  month={November},
  year={{\noopsort{2027}}2026}
}

@article {MR4973365,
    AUTHOR = {Yang, Yifan},
     TITLE = {Forward self-similar solutions to the {MHD} equations:
              existence and pointwise estimates},
   JOURNAL = {Discrete Contin. Dyn. Syst.},
  FJOURNAL = {Discrete and Continuous Dynamical Systems},
    VOLUME = {47},
      YEAR = {2026},
     PAGES = {406--456},
      ISSN = {1078-0947,1553-5231},
   MRCLASS = {76W05 (35A01 35Q35)},
  MRNUMBER = {4973365},
       DOI = {10.3934/dcds.2025124},
       URL = {https://doi-org.ezproxy.library.ubc.ca/10.3934/dcds.2025124},
}

@article {MR2514362,
    AUTHOR = {He, Cheng and Xin, Zhouping},
     TITLE = {On the self-similar solutions of the magneto-hydro-dynamic
              equations},
   JOURNAL = {Acta Math. Sci. Ser. B (Engl. Ed.)},
  FJOURNAL = {Acta Mathematica Scientia. Series B. English Edition},
    VOLUME = {29},
      YEAR = {2009},
    NUMBER = {3},
     PAGES = {583--598},
      ISSN = {0252-9602,1572-9087},
   MRCLASS = {35Q35 (35C06 35Q60 76W05)},
  MRNUMBER = {2514362},
MRREVIEWER = {Elder\ Jes\'us\ Villamizar-Roa},
       DOI = {10.1016/S0252-9602(09)60055-2},
       URL = {https://doi-org.ezproxy.mtsu.edu/10.1016/S0252-9602(09)60055-2},
}

@article {MR3997556,
    AUTHOR = {Zhang, Jingjing and Zhang, Ting},
     TITLE = {Global existence of discretely self-similar solutions to the
              generalized {MHD} system in {B}esov space},
   JOURNAL = {J. Math. Phys.},
  FJOURNAL = {Journal of Mathematical Physics},
    VOLUME = {60},
      YEAR = {2019},
    NUMBER = {8},
     PAGES = {081515, 18},
      ISSN = {0022-2488,1089-7658},
   MRCLASS = {76W05 (35C06 35Q35)},
  MRNUMBER = {3997556},
MRREVIEWER = {Guo\ Chun\ Wu},
       DOI = {10.1063/1.5092787},
       URL = {https://doi-org.ezproxy.mtsu.edu/10.1063/1.5092787},
}

@article {MR3609234,
    AUTHOR = {Lai, Baishun and Lin, Junyu and Wang, Changyou},
     TITLE = {Forward self-similar solutions to the viscoelastic
              {N}avier-{S}tokes equation with damping},
   JOURNAL = {SIAM J. Math. Anal.},
  FJOURNAL = {SIAM Journal on Mathematical Analysis},
    VOLUME = {49},
      YEAR = {2017},
    NUMBER = {1},
     PAGES = {501--529},
      ISSN = {0036-1410,1095-7154},
   MRCLASS = {76A10 (35C06 35Q30 76D03 76D05)},
  MRNUMBER = {3609234},
MRREVIEWER = {Du\v san\ Zorica},
       DOI = {10.1137/16M1060340},
       URL = {https://doi-org.ezproxy.mtsu.edu/10.1137/16M1060340},
}

@misc{GLX-2D-DSS,
  title={Large forward discretely self-similar solutions to the two-dimensional {N}avier--{S}tokes equations with rough data},
  author={Gui, Changfeng and Liu, Hao and Xie, Chunjing},
  note={arXiv:2609.28157},
  year={2026}
}

@article{LMZ-fractional,
  title={Forward self-similar solutions of the fractional {N}avier--{S}tokes equations},
  author={Lai, Baishun and Miao, Changxing and Zheng, Xiaoxin},
  journal={Advances in Mathematics},
  volume={352},
  pages={981--1043},
  year={2019},
  doi={10.1016/j.aim.2019.06.021}
}
\end{document}